\documentclass[ss,ejs]{imsart}

\RequirePackage[OT1]{fontenc}
\RequirePackage{amsthm,amsmath}
\RequirePackage[numbers]{natbib}

\pubyear{2017}
\volume{0}
\issue{0}
\firstpage{1}
\lastpage{40}

\startlocaldefs
\numberwithin{equation}{section}
\theoremstyle{plain}

\endlocaldefs

\usepackage[bf,normalsize,tableposition=top]{caption}
\usepackage{subfig}
\usepackage{amsfonts}
\usepackage{color}
\usepackage{epsfig}
\usepackage{srcltx}
\usepackage{graphicx}
\usepackage[colorlinks=true]{hyperref}
\usepackage{multirow}
\usepackage{slashbox}
\usepackage[dvipsnames]{xcolor}
\usepackage{soul}
\def\itemSG{4}
\def\itemdrisk{5}

\def\inter{\mathop{\hbox{\rm int}}}

\def\Risk{\hbox{\rm Risk}}
\def\W{{\cal W}}
\def\NW{\hbox{\rm NW}}
\def\Hplus{\hbox{\scriptsize$\left[\begin{array}{c|c}
H&h\cr\hline h^T\cr\end{array}\right]$}}
\def\Z{{\cal Z}}

\def\Erf{{\mathop{\hbox{\small\rm Erf}}}}
\def\ErfInv{{\mathop{\hbox{\small\rm ErfInv}}}}
\usepackage{amssymb}

\def\Diag{{\hbox{\rm Diag}}}

\def\bR{{\mathbf{R}}}

\def\Opt{{\hbox{\rm Opt}}}

\def\A{{\cal A}}
\def\L{{\cal L}}
\def\U{{\cal U}}
\def\SG{{\cal SG}}
\def\T{{\cal T}}
\def\E{{\cal E}}
\def\P{{\cal P}}

\def\G{{\cal G}}
\def\H{{\cal H}}

\def\X{{\cal X}}

\def\N{{\cal N}}

\def\K{{\cal K}}
\def\bE{{\mathbf{E}}}
\def\bS{{\mathbf{S}}}
\def\Prob{\hbox{\rm Prob}}
\def\Erf{\hbox{\rm Erf}}
\newtheorem{corollary}{Corollary}[section]
\newtheorem{lemma}{Lemma}[section]
\newtheorem{proposition}{Proposition}[section]
\newtheorem{remark}{Remark}[section]
\newtheorem{assump}{Assumption}

\def\Det{{\hbox{\rm Det}}}
\def\qed{$\Box$}
\def\Card{\mathop{\hbox{\rm Card}}}

\def\Tr{{\hbox{\rm Tr}}}
\def\Erf{{\mathop{\hbox{\small\rm Erf}}}}
\def\ErfInv{{\mathop{\hbox{\small\rm ErfInv}}}}
\newcommand{\be}{\begin{eqnarray}}
\newcommand{\ee}[1]{\label{#1}\end{eqnarray}}

\newcommand{\I}{{\cal I}}
\newcommand{\Ii}{{\cal I}}
\newcommand{\D}{{\cal D}}

\newcommand{\eight}{ \mbox{\small$\frac{1}{8}$}}

\newcommand{\rf}[1]{~(\ref{#1})}
\newcommand{\half}{ \mbox{\small$\frac{1}{2}$}}

\newcommand{\ese}{\end{eqnarray*}}
\newcommand{\bse}{\begin{eqnarray*}}

\def\airho{{\rho^2}}
\definecolor{MyDarkBlue}{rgb}{0,0.08,0.45}
\definecolor{MyViolet}{rgb}{0.45,0.08,0.95}
\definecolor{MyBrown}{rgb}{0.45,0.08,0}

\definecolor{light-gray}{gray}{0.95}
\def\itemffinal{3}
\newcommand{\an}[2]{{\color{blue}\  #2}}
\newcommand{\aic}[2]{{\color{cyan}\  #2}}

\begin{document}

\begin{frontmatter}
\title{Change Detection via Affine and Quadratic Detectors
}
\runtitle{Change Detection via Affine and Quadratic Detectors}

\begin{aug}

\author{\fnms{Yang} \snm{Cao}\thanksref{t1}\ead[label=e1]{caoyang@gatech.edu}}

\address{H. Milton Stewart School of Industrial and Systems Engineering\\
Georgia Institute of Technology, Atlanta, Georgia 30332, USA\\
\printead{e1,e4,e5}}

\author{\fnms{Vincent} \snm{Guigues}\thanksref{t2}
\ead[label=e2]{vincent.guigues@fgv.br}}

\address{School of Applied Mathematics, FGV\\ 190 Praia de Botafogo, Rio de Janeiro, RJ, 22250-900, Brazil\\
\printead{e2}}

\author{\fnms{Anatoli} \snm{Juditsky}\thanksref{t3}\ead[label=e3]{anatoli.juditsky@imag.fr}}

\address{LJK, Universit\'e Grenoble Alpes\\
700 Avenue Centrale 38041 Domaine Universitaire
de Saint-Martin-d'H\`{e}res, France\\
\printead{e3}}

\author{\fnms{Arkadi} \snm{Nemirovski}\thanksref{t4}\ead[label=e4]{nemirovs@isye.gatech.edu}}
\and
\author{\fnms{Yao} \snm{Xie}\thanksref{t1}\ead[label=e5]{yao.xie@isye.gatech.edu}}

\thankstext{t1}{Research was partially supported by NSF grants CAREER CCF-1650913, CCF-1442635, CMMI-1538746, and a Georgia Institute of Technology IMat Faculty Fellow (IFF) Seed Grant.
}
\thankstext{t2}{Research was partially supported by an FGV grant, CNPq grants 307287/2013-0 and 401371/2014-0, and
FAPERJ grant E-26/201.599/2014.}
\thankstext{t3}{Research was supported by the LabEx PERSYVAL-Lab (ANR-11-LABX-0025) and {CNPq grant 401371/2014-0}.}
\thankstext{t4}{Research was supported by NSF grants CCF-1523768, and CCF-1415498 and {CNPq grant 401371/2014-0}.}
\runauthor{Cao et al.}


\end{aug}

\begin{abstract}
The goal of the paper is to {develop a specific application of} the convex optimization based hypothesis testing techniques  developed in
{\sl A. Juditsky, A. Nemirovski, ``Hypothesis testing via affine detectors,'' Electronic Journal of
Statistics} {\bf 10}:2204--2242, 2016. Namely, we consider the Change Detection problem as follows:
observing one by one noisy observations of outputs of a discrete-time linear dynamical system, we intend to decide, in a sequential fashion, on the null
{hypothesis that} the input to the system is a nuisance, vs. the {alternative that} the input is a ``nontrivial signal,'' with both the nuisances and the nontrivial signals modeled as inputs belonging to  finite unions of some given convex sets. Assuming the observation noises are zero mean sub-Gaussian, we develop ``computation-friendly'' sequential decision rules and demonstrate that in our context these rules are provably near-optimal.
\end{abstract}

\begin{keyword}[class=MSC]
\kwd[Primary ]{62C20}
\kwd[; secondary ]{90C22}
\end{keyword}

\begin{keyword}
change-point detection, semi-definite program
\end{keyword}
\tableofcontents

\end{frontmatter}

\section{Introduction}
Quick detection of change-points from data streams is a classic and fundamental problem in signal processing and statistics, with a wide range of applications from cybersecurity \cite{LakhinaCrovellaDiot2004} to gene mapping \cite{Siegmund07}. Classical statistical change-point detection \cite{Siegmund1985,non-parametric-change93,Chen12,Veeravalli13book,Tartakovsky2014}, where one monitors {\it i.i.d.} univariate and low-dimensional multivariate observations is a well-developed area. Outstanding  contributions include Shewhart's control chart \cite{Shewhart31}, Page's CUSUM procedure \cite{Page1954},
Shiryaev-Roberts procedure \cite{Shiryaev1963}, Gordon's non-parametric procedure \cite{gordon1994efficient}, and window-limited procedures \cite{lai1995sequential}.
Various asymptotic (see, e.g., \cite{Lorden1971,Pollak85,Pollak87,lai1995sequential,lai1998information}) and nonasymptotic \cite{Moustakides86} results have been established for these classical methods. High-dimensional change-point detection (also referred to as the multi-sensor change-point detection) is a more recent topic, and various statistical procedures are proposed including \cite{ingster02,ingster09,Ingster10,korostelev2008,Mei2008,VVV08,levy09,Fellouris14,Zaid2014,Xie2015,Mei16}.
However, there has been very little research on the computational aspect of change-point detection,  especially in the high-dimensional setting.

\subsection{Outline}\label{sectOutline}
This paper presents a computational framework to solve change-point detection problems which is completely general: it can process many high-dimensional situations achieving improved false detection control.  The main idea is to adapt the framework for hypothesis testing using convex optimization \cite{PartI} to change-point detection. Change-point detection can be viewed as a multiple-testing problem, where at each time, one has to test whether there has been no change, or there already has been a change-point. With our approach, at each time a detector is designed by convex optimization to achieve the above goal.
The convex optimization framework is computationally efficient and can control false detection uniformly according to a pre-specified level.

Since change-point detection in various settings is the subject of huge literature (see, e.g.,
\cite{Bass1993,neumann1997optimal,TartVeer2004,GJTZ2008,goldenshluger2008change,Poor2009,Xie2015,Chen12,lai1998information,Siegmund1985,non-parametric-change93,Veeravalli13book,Tartakovsky2014} and references therein), it would be too
time-consuming to position our developments w.r.t. those presented in the literature. Instead,
we illustrate our approach by its application to a  simple example and  {then comment on the ``spirit'' of} our constructions and results (which, we believe, is somehow different from majority of traditional approaches to change detection).\par
\paragraph{Illustrating problem.} We consider a simple version of the classical  problem of change detection in the input of a dynamical system (see, e.g., \cite{gholson1977,willsky1985,mazor1998} and
references therein), where we observe
noisy outputs $\omega_t\in\bR^\nu$ of a discrete time linear time invariant system on time horizon $t=1, \ldots, d$:
\begin{equation}\label{ieq1}
\begin{array}{rcl}
x_t&=&Ax_{t-1}+bu_t,\\
\omega_t&=&Cx_t+\xi_t,
\end{array}
\end{equation}
where the inputs $u_t$ are scalars, $A$, $b$, $C$ are known, and the observation noises $\xi_t\sim\N(0,I_\nu)$ are independent across time $t=1,...,d$. The input $u=[u_1;...;u_d]$ to the system can be either zero ({\sl nuisance hypothesis}), or a signal of ``some shape $\tau\in\{1,...,d\}$ and some magnitude $\geq \rho>0$,'' meaning that $u_t=0$ for $t<\tau$, and $u_\tau\geq\rho$ (so that $\tau$ represents the change-point location in time); we refer to the latter option as to the {\sl {signal} hypothesis}. We observe $\omega_t$'s one by one, and our goal is to design
decision rules $\{\T_t:1\leq t\leq d\}$ and thresholds $\{\rho_{t\tau}>0,1\leq\tau\leq t\leq d\}$
in such a way that
\par$\bullet$ rule $\T_t$ is invoked at time $t$. Depending solely on the observations $\omega_1,...,\omega_t$ available at this time, this rule
\begin{itemize}
\item either accepts the {signal} hypothesis, in which case we terminate with ``{signal}'' conclusion,
 \item or claims that so far the nuisance hypothesis {is not rejected} (``nuisance conclusion at time $t$''), in which case we pass to time instant $t+1$ (when $t<d$) or terminate (when $t=d$);
\end{itemize}
\par$\bullet$
whenever the nuisance hypothesis is true, the probability of false alarm ({signal} conclusion somewhere on time horizon $t=1,...,d$) is at most a given $\epsilon\in(0,1/2)$;
\par$\bullet$
for every $t\leq d$ and every $\tau\leq t$, whenever the input is a signal of shape $\tau$ and magnitude $\geq\rho_{t\tau}$, the probability of {signal} conclusion at time $t$ or earlier is at least $1-\epsilon$. In other words, for every input of shape $\tau$ and magnitude $\geq\rho_{t\tau}$, the probability  of the nuisance conclusions at  all time instants $1,2,...,t$ should be at most $\epsilon$.
\par\noindent
{In what follows we refer to $\epsilon$ as to {\em risk of the collection} $\{\T_t,\,1\leq t\leq d\}$.} Needless to say, we would like to meet
the outlined design specifications with as small thresholds $\rho_{t\tau}$ as possible.\par
Our related results can be summarized as follows: we develop specific decision rules $\T_t$ and thresholds $\rho_{t\tau}$ meeting the design  specifications  and such that
\par$\bullet$ $\T_t$ and $\rho_{t\tau}$ are yielded by explicit convex optimization problems and thus can be built in a computationally efficient fashion; moreover,
the decision rules $\T_t$ are easy to implement;

\par$\bullet$ the resulting inference procedure is near-optimal in some precise sense. Specifically,
for every $\tau$ and $t$, $1\leq \tau\leq t\leq d$, consider the testing problem where, given the observations $\omega_1,...,\omega_t$, we want to decide on only two hypotheses on input $u$ underlying the observations: the hypothesis $H_1$ ``$u=0$'' and the alternative $H_2(\rho)$ ``$u$ is a signal of shape $\tau$ and magnitude $\geq\rho$,'' where $\rho>0$ is a parameter. It may happen that these two hypotheses can be decided upon with risk $\leq \epsilon$, meaning that ``in the nature'' there exists a test which, depending on observations $\omega_1,...,\omega_t$, accepts exactly one of the hypotheses with  error  probabilities (i.e., probability to reject $H_1$ when $u=0$ and the probability to reject $H_2(\rho)$ when $u$ is a signal of shape $\tau$ and magnitude $\geq\rho$) at most $\epsilon$.  One can easily find the smallest $\rho=\rho^*_{t\tau}$ for which such a test exists\footnote{Note that
the observation $(\omega_1,...,\omega_t)$ is of the form ${\bar A}_tu+\xi^t$ with standard (zero mean, unit covariance matrix) Gaussian noise $\xi^t=(\xi_1,...,\xi_t)$. It is immediately seen
that $\rho^*_{t\tau}$ is the smallest $\rho$ for which the distance $\inf_{v\in V_{t\tau}}\|v\|_2$ from the origin to the convex set
$V_{t\tau}=\{{\bar A}_tu:\;u_s=0\mbox{ for }s<\tau,\,u_\tau\geq\rho\}$ is at least $2\ErfInv(\epsilon)$, where $\ErfInv$ is the inverse error function, see (\ref{inverf}).}. Clearly, by construction, $\rho^*_{t\tau}$ is a lower bound on the threshold $\rho_{t\tau}$ of any inference routine which meets the design specifications we are dealing with.
 Near-optimality of our inference routine means, essentially, that our thresholds $\bar{\rho}_{t\tau}$ are close to the ``ideal'' thresholds $\rho^*_{t\tau}$ independently of particular  values of parameters of model \rf{ieq1}:
    \[
    {\bar{\rho}_{t\tau}\over \rho^*_{t\tau}}\leq {1\over 2}\left[1+{\ErfInv(\epsilon/d^2)\over\ErfInv(\epsilon)}\right]
    \] (for details, see Proposition \ref{prop40}).
\paragraph{Paper's scope.} The developments to follow are in no sense restricted to the simplest model of nuisance and signal inputs we have considered so far. In fact, we
{allow  nuisance} inputs to vary in a  prescribed set $N\ni 0$,  and for signal inputs to have $K$  different ``shapes,'' with signals of ``shape $k\leq K$ and magnitude $\rho>0$'' varying in prescribed sets  $U_k(\rho)$ shrinking as $\rho>0$ grows. We treat
 {two cases separately}:
\par {\bf I.} {\sl ``Decision rules based on affine detectors,''} in Section  \ref{ChPADSetUp}. In this case, $N\ni 0$ is a convex compact set, and $U_k(\rho)=N+\rho W_k$, where $W_k$ are closed convex sets  not containing the origin and such that $\rho W_k\subset W_k$ whenever $\rho\geq1$, implying that $U_k(\rho)$ indeed shrinks as $\rho$ grows. As far as the observation noises are concerned, we require the vector $\xi^d=(\xi_1,...,\xi_d)$ to be zero mean sub-Gaussian, with the (perhaps, unknown)  matrix {parameter} (see \itemSG, Section \ref{sectnotconv}) belonging to a given convex compact set. This case
 {covers the } example we have started with.
\par{\bf II.} {\sl ``Decision rules based on quadratic detectors,''} in Section \ref{QDSection}. In this case, $N\ni 0$ is a bounded set given by a finite system of quadratic inequalities, and  $U_k(\rho)$, $1\leq k\leq K$, is given by a parametric system of quadratic inequalities of appropriate structure (for details, see Section \ref{newsetup}). The simplest illustration here is the case when $u_t$ in (\ref{ieq1}) are allowed to be vectors, the only nuisance input is $u=0$, and a signal input of shape $\tau\leq d$ of magnitude $\geq \rho$ is a block-vector $[u_1;...;u_d]$ with $u_1=...=u_{\tau-1}=0$ and $\|u_\tau\|_2\geq\rho$.
 {The noise $\xi^d=[\xi_1;...;\xi_d]$ is assumed} to be zero mean Gaussian, with (perhaps, unknown) covariance matrix varying in a known convex compact set.
\paragraph{Comments.} To complete the introduction, let us comment on the ``spirit'' of our constructions and results, which we refer to as {\em operational}. Following the line of research in \cite{GJN,Seq2015,PartI}, we allow for rather general {\sl structural} assumptions on the components of our setup (system \rf{ieq1} and descriptions of nuisance and signal inputs) and are looking for {\sl computation-friendly} inference routine
meaning that our easy-to-implement routines and their performance characteristics are given by {\sl efficient computation} (usually  based on Convex Optimization). This {appears to be in sharp} contrast {with} the traditional in statistics ``closed analytical form'' {\sl descriptive} procedures and performance characteristics. While closed analytical form results possess strong explanatory power, these results usually impose severe restrictions on the underlying setup and in this respect are much more restrictive than operational results. We believe that in many applications, including those considered in this paper, the relatively broad applicability of operational results more than compensates for the lack of explanatory
 {power that is typical of computation-based} constructions. It should be added that under favorable circumstances (which, in the context of this paper, do take place in case I), the operational procedures we are about to develop
 {are provably} near-optimal in a certain precise sense (see Section \ref{sect:assessing}). Therefore,  their performance, whether good or bad from the viewpoint of a particular application, is nearly the best possible under the circumstances.

\subsection{Terminology and notation}\label{sectnotconv}
In what follows,
\par 1.
 All vectors are column vectors.
 \par
 2.
 We use ``MATLAB notation:'' for matrices $A_1,...,A_k$ of common width, $[A_1;A_2;...;A_k]$ stands for the matrix obtained by {(up-to-down) vertical concatenation of $A_1, A_2,..., A_k$;} for matrices $A_1,...,A_k$ of common height, $[A_1,A_2,...,A_k]$ is the matrix obtained by {(left-to-right) horizontal concatenation of $A_1, A_2, ...,A_k$.}
\par
3.
 $\bS^n$ is the space of $n\times n$ real symmetric matrices, and $\bS^n_+$ is the cone of positive semidefinite matrices from $\bS^n$. Relation $A\succeq B$ ($A\succ B$) means that $A$, $B$ are symmetric matrices of the same size such that $A-B$ is positive semidefinite (respectively, positive definite), and $B\preceq A$ ($B\prec A$) is the same as $A\succeq B$ (respectively, $A\succ B$).

\par
\itemSG. $\SG[U,\U]$, where $U$ is a nonempty subset of $\bR^n$, and $\U$ is a nonempty subset of $\bS^n_+$, stands for the family of all Borel sub-Gaussian
probability distributions on $\bR^n$ with sub-Gaussianity parameters from $U\times\U$. In other words, $P\in \SG[U,\U]$ if and only if $P$ is a probability distribution such that for some $u\in U$ and $\Theta\in \U$ one has  $\ln(\int {\rm e}^{h^Ty}P(dy))\leq u^Th+{1\over 2}h^T\Theta h$ for all $h\in \bR^n$ (whenever this is the case, $u$ is the expectation of $P$); we refer to $\Theta$ as to {\sl sub-Gaussianity matrix} of $P$. For a random variable $\xi$ taking values in $\bR^n$, we write $\xi\sim \SG[U,\U]$ to express the fact that the distribution $P$ of $\xi$ belongs to $\SG[U,\U]$.\par
     Similarly, $\G[U,\U]$ stands for the family of all Gaussian distributions $\N(u,\Theta)$ with expectation $u\in U$ and covariance matrix $\Theta\in\U$, and $\xi\sim \G[U,\U]$ means that $\xi\sim \N(u,\Theta)$ with $u\in U$, $\Theta\in\U$.
 \par
\itemdrisk.
  Given two families $\P_1$, $\P_2$ of Borel probability distributions on $\bR^n$ and a {\sl detector} $\phi$ (a Borel real-valued function on $\bR^n$),
$\Risk(\phi|\P_1,\P_2)$ stands for the risk of the detector \cite{GJN} taken w.r.t. the families $\P_1$, $\P_2$, that is, the smallest $\epsilon$ such that
    \begin{equation}\label{riskis}
    \begin{array}{lrcl}
    (a)&\int {\rm e}^{-\phi(y)}P(dy)&\leq&\epsilon\,\;\forall P\in\P_1,\\
    (b)&\int {\rm e}^{\phi(y)}P(dy)&\leq&\epsilon\,\;\forall P\in\P_2.\\
    \end{array}
    \end{equation}
    When $\T$ is a test deciding on $\P_1$ and $\P_2$ via random observation $y\sim P\in\P_1\cup\P_2$ (that is,  $\T:\bR^n\to\{1,2\}$ is a Borel function, with $\T(y)=1$ interpreted as ``given observation $y$, the test accepts the hypothesis $H_1:P\in\P_1$ and rejects the hypothesis $H_2:P\in\P_2$,'' and $\T(y)=2$ interpreted as ``given observation $y$, $\T$  accepts $H_2$ and rejects $H_1$'')
$$\begin{array}{rcl}
\Risk_1(\T|\P_1,\P_2)&=&\sup_{P\in\P_1} \Prob_{y\sim P}\{y:\T(y)=2\},\\
\Risk_2(\T|\P_1,\P_2)&=&\sup_{P\in\P_2}\Prob_{y\sim P}\{y:\T(y)=1\}\\
\end{array}
$$
stand for the partial risks of the test,  and
\[
    \Risk(\T|\P_1,\P_2)=\max[\Risk_1(\T|\P_1,\P_2),\Risk_2(\T|\P_1,\P_2)]
    \] stands for the risk of the test.
    \par
    A detector $\phi(\cdot)$ and a real $\alpha$ specify a test $\T^{\phi,\alpha}$ which accepts $H_1$ ($\T^{\phi,\alpha}(y)=1$) when $\phi(y)\geq\alpha$, and accepts $H_2$ ($\T^{\phi,\alpha}(y)=2$) otherwise. From (\ref{riskis}) it is immediately seen that
   \begin{equation}\label{mainbasis}
    \begin{array}{l}
    \Risk_1(\T^{\phi,\alpha}|\P_1,\P_2)\leq{\rm e}^{\alpha}\Risk(\phi|\P_1,\P_2),\\
    \Risk_2(\T^{\phi,\alpha}|\P_1,\P_2) \leq  {\rm e}^{-\alpha}\Risk(\phi|\P_1,\P_2).
    \end{array}
    \end{equation}
\par

All proofs are transferred to the appendix.

\section{Dynamic change detection: preliminaries}\label{sect:ChPD}
In the sequel, we address the situation which can be described informally as follows.
We observe
noisy outputs of a linear system at times $t=1,...,d$, the input to the system  being an unknown vector $x\in\bR^n$. Our ``full observation'' is
\begin{equation}\label{eqOS}
y^d=\bar{A}_dx+\xi^d,
\end{equation}
where $\bar{A}_d$ is a given $\nu_d\times n$ {\sl sensing matrix}, and $\xi^d\sim\SG[\{0\},\U]$ (see item \itemSG\ in Section \ref{sectnotconv}), where $\U$ is a given nonempty convex compact subset of $\inter\bS_+^{\nu_d}$.
\par
Observation $y^d$ is obtained in $d$ steps; at a step (time instant) $t=1,...,d$, the observation is
\begin{equation}\label{eqOSt}
y^t=\bar{A}_tx+\xi^t\equiv S_{t}[\bar{A}_dx+\xi^d]\in\bR^{\nu_t},
\end{equation}
where $1\leq \nu_1\leq \nu_2\leq...\leq\nu_d$, $S_t$ is $\nu_t\times \nu_d$ matrix of rank $\nu_t$  and $y^t$ ``remembers'' $y^{t-1}$, meaning that  $S_{t-1}=R_tS_t$ for some matrix $R_t$.
Clearly, $\xi^t$ is sub-Gaussian with parameters $(0,\Theta_t)$,
with
 \begin{equation}\label{thetat}
\Theta_t={S_t}\Theta S_t^T\subset \U_t:=\{{S_t}\Theta S_t^T:\,\Theta\in\U\};
\end{equation} note that $\U_t$, $1\leq t\leq d$, are convex compact sets comprised of positive definite $\nu_t\times \nu_t$ matrices.
\par
Our goal is to build a {\sl dynamic test} for deciding on the {\sl null}, or {\sl nuisance}, hypothesis, stating that the input to the system underlying our observations is a nuisance, vs. the alternative of a {\sl signal} input.
Specifically, at every {time $t=1,...,d$,} given observation $y^t$, we can either decide that the input is a signal and terminate (``termination at step $t$ with a signal conclusion,'' or, equivalently, ``detection of a signal input at time $t$''), or to decide (``nuisance conclusion at step $t$'') that so far, the nuisance hypothesis holds true, and to pass to the next time instant $t+1$ (when $t<d$) or to terminate (when $t=d$).
\par
Given an upper bound $\epsilon$ on the probability of {a false alarm} (detecting a signal input somewhere on the time horizon $1,...,d$ in the situation when the true input is a nuisance), our informal goal is to build a dynamic test which respects the false alarm bound and under this restriction, detects signal inputs ``as fast as possible.''  \par
We consider two different types of detection procedures, those based on {\sl affine} and on {\sl quadratic detectors}, each type dealing with its own structure of nuisance and signal inputs.
\section{Change detection via affine detectors}\label{ChPADSetUp}
We start with
describing the structure of nuisance and signal inputs that we intend to deal with.
\subsection{Setup}\label{ADSetUp}
Consider the setup as follows.
\par1.
Inputs to the system belong to a given convex compact set $X\subset\bR^n$, and {\sl nuisance} inputs form a given
closed and convex subset $N$ of $X$, with $0\in N$.
\par2. {\sl Informal} description of a signal input $x$ is as follows: $x\in X$ is obtained from some nuisance input $v$ by adding an ``activation'' $w$ of some {\sl shape} and some {\sl magnitude}. There are $K$ possible shapes, $k$-th of them represented by a closed convex set $W_k\subset\bR^n$ such that
    \begin{enumerate}
    \item[2.1.] $0\not\in W_k$;
    \item[2.2.] $W_k$ is {\sl semi-conic}, meaning that when $w\in W_k$ and $\rho\geq1$, it holds $\rho w\in W_k$.
    \end{enumerate}
    The magnitude of an activation is just a positive real, and an activation of shape $k$ and magnitude {\sl at least} $\rho>0$ {is an element} of the set
    $$
    W_k^\rho=\{w=\rho y:y\in W_k\}.
    $$
    \begin{quote}
    {\small
{\bf Example:} Let $K=n$ and let $W_k$ be the set of all inputs $w\in\bR^n$ with the first $k-1$ entries in $w$ equal to zero, and $k$-th entry $\geq1$. In this case, the shape of an activation $w\in\bR^n$  is its ``location'' -- the index of the first nonzero entry in $w$, and activations of shape $k$ and magnitude $\geq\rho$ are vectors $w$  from $\bR^n$ with the first nonzero entry in position $k$ and the value of this entry at least $\rho$.\\
We have presented the simplest formalization of what informally could be called ``activation up.''  To get equally simple formalization of an ``activation down,'' one should take $K=2n$ and define $W_{2i-1}$  and $W_{2i}$, $i\leq n$, as the sets of all vectors from $\bR^n$ for  which the first nonzero entry is in position $i$, and the value of this entry is at least 1 for $W_{2i-1}$ (``activation up'' of magnitude $\geq1$ at time $i$)   or is at most $-1$ for $W_{2i}$ (``activation down'' of magnitude $\geq1$ at time $i$).}
\end{quote}
\par3. {\sl The formal} description of ``signal'' inputs is as follows: these are vectors $x$ from $X$ which for some $k\leq K$ can be represented as $x=v+w$ with $v\in V_k$ and $w\in W_k^\rho$ for some $\rho >0$, where $W_k$ are as described above, and $V_k$, $0\in V_k$, are nonempty compact convex subsets of $X$.\footnote{In the informal description of signals, $V_k$ were identified with the set $N$ of nuisances; now we lift this restriction in order to add more flexibility.} Thus, when speaking about signals (or signal inputs), we assume that we are given  $K$ nonempty closed convex sets $W_k$, $k\leq K$, each of them semi-conic and not containing the origin, and $K$ nonempty compact convex sets $V_k\subset X$. These sets give rise to single-parametric families of compact convex sets
\[
\begin{array}{rclcrcl}
W^\rho_k&=&\{\rho y:y\in W_k\},&
X^\rho_k&=&[V_k+W^\rho_k]\bigcap X,\\
\end{array}
\]
indexed by ``activation  shape'' $k$ and
parameterized by ``activation  magnitude'' $\rho>0$. Signals are exactly the {elements of} the set $\widehat{X}=\bigcup_{{\rho>0,k\leq K}}X^\rho_k$. In the sequel, we refer to inputs from $N$ as to {\sl feasible nuisances}, to inputs from $X^\rho_k$ as to {\sl feasible signals with activation of shape $k$ and magnitude $\geq\rho$}, and to inputs from $\widehat{X}$ as to {\sl feasible signals}.
To save words, in what follows `` {a} signal of shape $k$ and magnitude $\geq\rho$'' means exactly the same as `` {a} signal with activation of shape $k$ and magnitude $\geq\rho$.''
\par
From now on, we make the following assumption:
\begin{assump}\label{ass:0}
For every $k\leq K$, there exists $R_k>0$ such that the set $X^{R_k}_k$ is nonempty.
\end{assump}
Since $X^\rho_k$ shrinks as $\rho$ grows due to semi-conicity of $W_k$, it follows that {\sl for every $k$, the sets $X^\rho_k$ are nonempty for all small enough positive $\rho$.}
\subsection{Construction}\label{sgcase}
\subsubsection{Outline}\label{sec:outline}
Given an upper bound $\epsilon\in(0,1/2)$ on the probability of false alarm, our course of actions {is as follows}.
 \par1. We select $d$ positive reals $\epsilon_t$, $1\leq t\leq d$, such that $\sum_{t=1}^d\epsilon_t=\epsilon$; $\epsilon_t$ will be an upper bound on the probability {of a false alarm} at time $t$.
 \par2.
We select thresholds $\rho_{tk}>0$, $1\leq k\leq K$
in such a way that a properly designed test $\T_t$ utilizing the techniques of \cite[Section 3]{PartI} is able to distinguish reliably, given an observation $y^t$, between the hypotheses $H_{1,t}:x\in N$ and $H_{2,t}:x\in\bigcup\limits_{k=1}^KX^{\rho_{tk}}_k$ on the input $x$ underlying observation $y^t$. After $y^t$ is observed, we
 apply test $\T_t$ to this observation, and, according to what the test says,
 \begin{itemize}
 \item either claim that the input is a signal, and terminate,
 \item or claim that {\sl so far}, the hypothesis of nuisance input  seems to be valid, and either pass to the next observation (when $t<d$), or terminate (when $t=d$).
 \end{itemize}
The generic construction we intend to use when building the test $\T_t$ stems from \cite{GJN,PartI}.
 \subsubsection{Implementation: preliminaries}
\paragraph{Building block: affine detectors for sub-Gaussian families.} Our principal building block originates from \cite{PartI} and is as follows. Let $\U$ be a convex compact set comprised of positive definite $\nu\times \nu$ matrices, and  $U_1$, $U_2$ be two closed  nonempty convex subsets in $\bR^\nu$, with $U_1$ bounded. The following result was proved in \cite{PartI}:
\begin{proposition}\label{subG}{\rm \cite[Propositions 3.3 and 3.4]{PartI}} With $\U$, $U_1$, $U_2$ as above, consider the convex-concave saddle point problem
\[
{\cal SV}=\min_{h\in\bR^\nu}\max_{\theta_1\in U_1,\theta_2\in U_2,\Theta\in\U}\left[\Phi(h;\theta_1,\theta_2,\Theta):=\half h^T[\theta_2-\theta_1] + \half h^T\Theta h\right].
\]
This saddle point problem is solvable, and a saddle point $(h_*=\half\Theta_*^{-1}[\theta_1^*-\theta_2^*];\theta_1^*,\theta_2^*,\Theta_*)$ induces affine detector
\[
\phi_*(\omega)=h_*^T(\omega-w_*),\,\,w_*=\half[\theta_1^*+\theta_2^*]
\]
 for the families of distributions $\P_1=\SG[U_1,\U]$ and $\P_2=\SG[U_2,\U]$ (for notation, see item \itemSG\ in Section \ref{sectnotconv}), and the risk
  $\Risk(\phi_*|\P_1,\P_2)$ of this detector (see item \itemdrisk\ in Section \ref{sectnotconv}) is upper-bounded by
\[
\epsilon_\star={\rm e}^{{\cal SV}}=\exp\left({-\half  {h_*^T}\Theta_*h_*}\right)=\exp\left({-\eight[\theta_1^*-\theta_2^*]^T\Theta_*^{-1}[\theta_1^*-\theta_2^*]}\right).
\]
Moreover,
let
$$
\delta=\sqrt{h_*^T\Theta_*h_*},
$$
 and let $\alpha\leq\delta^2$, $\beta\leq\delta^2$.
Then
\begin{equation}\label{eq1100Aff}
\begin{array}{ll}
(a)&\forall (\theta\in U_1,\Theta\in \U):\; \Prob_{\omega\sim \N(\theta,\Theta)}\{\phi_*(\omega)\leq\alpha\}\leq \Erf(\delta-\alpha/\delta),\\
(b)&\forall (\theta\in U_2,\Theta\in \U): \;\Prob_{\omega\sim \N(\theta,\Theta)}\{\phi_*(\omega)\geq-\beta\}\leq \Erf(\delta-\beta/\delta),
\end{array}
\end{equation}
where
$$
\Erf(s)={1\over\sqrt{2\pi}}\int_s^\infty\exp\{-r^2/2\}dr
$$
is the normal error function. In particular, when deciding, via a single observation $\omega$, on Gaussian hypotheses $H^G_\chi$, $\chi=1,2$, with $H_\chi^G$ stating that $\omega\sim\N(\theta,\Theta)$ with $(\theta,\Theta)\in U_\chi\times\U$,
the risk  of the test which accepts $H_1^G$ when $\phi_*(\omega)\geq0$ and accepts $H_2^G$ otherwise is at most $\Erf(\delta)$.
\end{proposition}
Given $k\in\{1,...,K\}$, observe that the set $X^\rho_k$ is nonempty when $\rho>0$ is small enough (this was already assumed) and is empty for all large enough values of $\rho$ (since $X$ is compact and $W_k$ is a nonempty closed convex set not containing the origin). From these observations and compactness of $X$ it follows that {\sl there exists the largest $\rho=R_k>0$ for which $X^\rho_k$ is nonempty.}
\par
Let us fix $t\in\{1,...,d\}$, and let
\begin{equation}\label{letUt}
\U_t=\{S_t\Theta S_t^T:\Theta\in\U\}
 \end{equation}
 be the set of allowed covariance matrices of the observation noise $\xi^t$ in observation $y^t$, so that $\U_t$ is a convex compact subset of the interior of $\bS^{\nu_t}_+$.
According to our assumptions, for any nuisance input the distribution of the associated observation $y^t$, see (\ref{eqOSt}), belongs to the family
$\SG[N^t,\U_t]$, with
\begin{equation}\label{letNt}
N^t=\bar{A}_tN,
\end{equation}
where $N\subset X$ is the convex compact set of nuisance inputs. Given, along with $t$, an integer $k\leq K$ and a real $\rho\in(0,R_k]$, we can define the set
\begin{equation}\label{letUtkrho}
U^t_{k\rho}=\{\bar{A}_t x:x\in X_k^\rho\};
\end{equation}
whatever be a signal input from $X_k^\rho$, the distribution of
 {observation $y^t$ associated with $x$}  belongs to the family $\SG[U^t_{k\rho},\U_t]$. Applying
Proposition \ref{subG} to data $U_1=N^t$, $U_2=U^t_{k\rho}$, and $\U=\U_t$,  we arrive at the convex-concave saddle point problem
\begin{equation}\label{ccproblem}
{\cal SV}_{tk}(\rho)=\min\limits_{h\in\bR^{\nu_t}}\max\limits_{\theta_1\in N^t,\theta_2\in U^t_{k\rho},\Theta\in \U_t}
\left[\half h^T[\theta_2-\theta_1]+\half h^T\Theta h\right].
\end{equation}
The corresponding saddle point
$$
(h_{tk\rho};\theta_{tk\rho}^1,\theta_{tk\rho}^2,\Theta_{tk\rho})
$$
does exist and gives rise to the affine detector
\begin{equation}\label{eq600}
\phi_{tk\rho}(y^t)=h_{tk\rho}^T[y^t-w_{tk\rho}],\,\,w_{tk\rho}=\half[\theta_{tk\rho}^1+\theta_{tk\rho}^2],
\end{equation}
and risk
\begin{equation}\label{eq601}
\begin{array}{l}
\Risk(\phi_{tk\rho}|\SG[N^t,\U_t],\SG[U^t_{k\rho},\U_t])\leq \epsilon_{tk\rho}:=\exp\{{\cal SV}_{tk}(\rho)\}\\
\multicolumn{1}{r}{=
\exp\{-{\eight}[\theta_{tk\rho}^1-\theta_{tk\rho}^2]^T[\Theta_{tk\rho}]^{-1}[\theta_{tk\rho}^1-\theta_{tk\rho}^2]\}.}\\
\end{array}
\end{equation}
Therefore, in view of (\ref{mainbasis}),
\begin{equation}\label{eq602}
\begin{array}{rcll}
\int\limits_{\bR^{\nu_t}}\exp\{-\phi_{tk\rho}(y^t)\}P(dy^t)&\leq&\epsilon_{tk\rho}&\forall P\in\SG[N^t,\U_t],\\
\int\limits_{\bR^{\nu_t}}\exp\{\phi_{tk\rho}(y^t)\}P(dy^t)&\leq&\epsilon_{tk\rho}&\forall P\in\SG[U^t_{k\rho},\U_t].
\end{array}
\end{equation}
To proceed, we need the following simple observation:
\begin{lemma}\label{lem1}
For every $t\in\{1,...,d\}$ and $k\in\{1,...,K\}$, the function ${\cal SV}_{tk}(\rho)$ is concave, nonpositive and nonincreasing continuous function of $\rho\in(0,R_k]$, and
$\lim_{\rho\to+0}{\cal SV}_{tk}(\rho)=0$.\par
Moreover, if $\U$ contains a $\succeq$-largest element $\overline{\Theta}$, that is, $\overline{\Theta}\succeq \Theta$ for some $\overline{\Theta}\in \U$ and all $\Theta\in \U$, then $\Gamma_{tk}(\rho)=\sqrt{-{\cal SV}_{tk}(\rho)}$ is a nondecreasing continuous convex nonnegative function on $\Delta_k=(0,R_k]$.
\end{lemma}

\subsubsection{Implementation: construction}\label{impl:constr}
Recall that we have split the required false alarm probability $\epsilon$ between decision steps $t=1,...,d$:
$$
\epsilon=\sum_{t=1}^d\epsilon_t.\eqno{[\epsilon_t>0\,\forall t]}
$$
At time instant $t\in\{1,...,d\}$ we act as follows:
\par1. For $\varkappa\in(0,1]$, let
\[
\begin{array}{rcl}
\K_t(\varkappa)&=&\{k\leq K: {\cal SV}_{tk}(R_k)< \ln(\varkappa)\},\\
K_t(\varkappa)&=&\Card\K_t(\varkappa),
\end{array}
\]
so that $K_t(\varkappa)$ is nondecreasing and continuous from the left,
and let\footnote{Specific choices of parameters $\varkappa_t$, $K_t(\varkappa_t)$, etc., allow to control false alarm and signal miss probabilities; the rationale behind these choices becomes clear from the proof of Proposition \ref{prop16}.}
\begin{equation}\label{varkappat}
\varkappa_t=\sup\left\{\varkappa\in(0,1]:K_t(\varkappa)\leq{\epsilon\epsilon_t\over\varkappa^2}\right\}.
\end{equation}
Clearly, $\varkappa_t$ is well defined, takes values in $(0,1]$, and since $K_t(\varkappa)$ is continuous from the left, we have
\begin{equation}\label{eq700}
K_t(\varkappa_t)\leq {\epsilon\epsilon_t\over\varkappa_t^2}.
\end{equation}
For $k\in\K_t(\varkappa_t)$, we have $0=\lim_{\rho\to+0}{\cal SV}_{tk}(\rho)>\ln(\varkappa_t)$ and ${\cal SV}_{tk}(R_k)<\ln(\varkappa_t)$. Invoking Lemma \ref{lem1}, there exists (and can be rapidly approximated to high accuracy by  bisection) $\rho_{tk}\in(0,R_k)$ such that
\begin{equation}\label{balance}
{\cal SV}_{tk}(\rho_{tk})=\ln(\varkappa_t).
\end{equation}
After $\rho_{tk}$ is specified, we build the associated detector $\phi_{tk}(\cdot)\equiv \phi_{tk\rho_{tk}}(\cdot)$ according to (\ref{eq600}). Note that the risk (\ref{eq601}) of this detector is $\epsilon_{tk\rho_{tk}}=\varkappa_t$.
\par
For $k\not\in\K_t(\varkappa_t)$, we set $\rho_{tk}=+\infty$.
\par3.
Finally, we set
$
\alpha_t=\ln(\varkappa_t/\epsilon).
$
and process observation $y^t$ at step $t$ as follows:
\begin{itemize}
\item if there exists $k$ such that $\rho_{tk}<\infty$ and  $\phi_{tk}( {y^t})<\alpha_t$, we claim that the input underlying observation is a signal and terminate;
\item otherwise, we claim that so far, the nuisance hypothesis is not rejected, and pass to the next time instant $t+1$ (when $t<d$) or terminate (when $t=d$).
\end{itemize}

\subsubsection{Characterizing performance}
The performance of the above inference procedure can be described as follows:
\begin{proposition}\label{prop16}  {For any} zero mean sub-Gaussian, with parameter $\Theta\in\U$,  distribution of observation noise $\xi^d$ in {\rm (\ref{eqOS})}, one has:
\par {\rm(i)} when the input is a feasible nuisance, the probability of terminating with the signal conclusion at time $t\in\{1,...,d\}$ does not exceed $\epsilon_t$, and thus the probability
 {of a false alarm} is at most $\epsilon=\sum_{t=1}^d\epsilon_t$;
\par{\rm (ii)} when $t\in\{1,...,d\}$ and $k\in\{1,...,K\}$ are such that $\rho_{tk}<\infty$, and the input belongs to a set $X^\rho_k$ with $\rho\geq \rho_{tk}$, then the probability to terminate at step $t$ with the  signal conclusion is at least $1-\epsilon$.
\end{proposition}

\subsection{Refinement in the Gaussian case}\label{sect:refinement}
In the case when observation noise $\xi$ in (\ref{eqOS}) is $\N(0,\Theta)$ with $\Theta\in \U$, the outlined construction can be refined. Specifically, at a  time instant $t\leq d$ we now act as follows.
\subsubsection{Construction}
\par1.  {Let $\ErfInv$ stand for} the inverse error function:
\begin{equation}\label{inverf}
\ErfInv(\kappa)=\left\{\begin{array}{ll}\min\left\{r:\Erf(r)\leq \kappa\right\},&0\leq\kappa\leq1/2,\\
0,&1/2\leq\kappa<\infty.\\
\end{array}\right.
\end{equation}
Assuming $\epsilon<1/2$ and given $t\in\{1,...,d\}$, we set for $\delta\geq0$
\[
\begin{array}{rcl}
\L_t(\delta)&=&\{k\leq K: {\cal SV}_{tk}(R_k)<-\half\delta^2\},\\
L_t(\delta)&=&\Card \L_t(\delta),\\
\end{array}
\]
so that $L_k(\delta)$ is  {a} continuous from the right non-increasing function of $\delta\geq0$.
We put
\begin{equation}\label{deltat}
\begin{array}{c}
\delta_t=\inf\left\{\delta:\delta\geq \half\left[\ErfInv(\epsilon_t/L_t(\delta))+\ErfInv(\epsilon)\right]\right\}\;\;
\left[\mbox{with }\ErfInv(\epsilon_t/0):=0\right].\end{array}
\end{equation}
Clearly,  $\delta_t$ is well defined, is positive, and
\begin{equation}\label{eq801}
\delta_t\geq\half\left[\ErfInv(\epsilon_t/L_t(\delta_t))+\ErfInv
(\epsilon)\right],
\end{equation}
since $L_k(\delta)$ defined above is continuous from the right.
\par2. For $k\in\L_t(\delta_t)$, we have ${\cal SV}_{tk}(R_k)<-\half\delta_t^2$, and ${\cal SV}_{tk}(\rho)>-\half\delta_t^2$ for all small enough $\rho>0$. Invoking Lemma \ref{lem1}, there exists (and can be rapidly approximated to high accuracy by bisection) $\rho=\rho_{tk}\in\Delta_k$ such that
\begin{equation}\label{balance1}
{\cal SV}_{tk}(\rho_{tk})=-\half\delta_t^2.
\end{equation}
After $\rho_{tk}$ is specified, we define the associated detector $\phi_{tk}(\cdot)\equiv \phi_{tk\rho_{tk}}(\cdot)$ by
 {applying the construction from Proposition \ref{subG}
to the data $U_1=N^t$, $U_2=U^t_{k\rho_{tk}}$, $\U=\U_t$}(see (\ref{letUt}), (\ref{letNt}), (\ref{letUtkrho})),  that is, find a
saddle point
$(h_*;\theta_1^*,\theta_2^*,\Theta_*)$ of the convex-concave function
$$
\half h^T[\theta_2-\theta_1]+\half h^T\Theta h:\bR^{\nu_t}\times(N^t\times U^t_{k\rho_{tk}}\times\U_t)\to\bR
$$
(such a saddle point does exist).
By Proposition \ref{subG}, the affine detector
$$
\phi_{tk}(y^t)= h_*[y^t-w_*],\,\,w_*=\half[\theta_1^*+\theta_2^*]
$$
has the risk  {bounded by}
\begin{equation}\label{eq367}
\exp\{-\half\delta^2\}=\epsilon_\star=\exp\{\half h_*^T[\theta_2^*-\theta_1^*]+\half  {h_*^T}\Theta_*h_*\}.
\end{equation}
Moreover (see (\ref{eq1100Aff})), for all $\alpha\leq\delta^2$ and $\beta\leq\delta^2$ it holds
\begin{equation}\label{eq1111}
\begin{array}{lll}
(a)&\forall (\theta\in N^t,\Theta\in \U_t):&\Prob_{y^t\sim \N(\theta,\Theta)}\{\phi_{tk}(y^t)\leq\alpha\}\leq \Erf(\delta-\alpha/\delta),\\
(b)&\forall (\theta\in U^t_{k\rho_{tk}},\Theta\in \U_t):&\Prob_{y^t\sim \N(\theta,\Theta)}\{\phi_{tk}(y^t)\geq-\beta\}\leq \Erf(\delta-\beta/\delta).
\end{array}
\end{equation}
Comparing the second equality in (\ref{eq367}) with the description of ${\cal SV}_{tk}(\rho_{tk})$, we see that $\epsilon_\star=\exp\{{\cal SV}_{tk}(\rho_{tk})\}$, which
combines with the first equality in (\ref{eq367}) and with (\ref{balance1}) to imply that $\delta$ in (\ref{eq367}) is nothing but $\delta_t$ as given by (\ref{deltat}). The bottom line is that
\begin{quotation}\noindent
$(\#)$ {\sl For $k\in\L_t(\delta_t)$, we have defined reals $\rho_{tk}\in\Delta_k$ and affine detectors $\phi_{tk}(y^t)$ such that relations {\rm (\ref{eq1111})} are satisfied with $\delta=\delta_t$ given by {\rm (\ref{deltat})} and  every $\alpha\leq\delta_t^2,\beta\leq\delta_t^2$.}
\end{quotation}
For $k\not\in\L_t(\delta_t)$, we set $\rho_{tk}=\infty$.
\par\itemffinal.\
Finally, we process observation $y^t$ at step $t$ as follows. We set
\begin{equation}\label{eqalphabeta}
\alpha=-\beta={\delta_t\over 2}\left[\ErfInv(\epsilon)-\ErfInv(\epsilon_t/L_t(\delta_t))\right],
\end{equation}
thus ensuring, in view of (\ref{eq801}), that $\alpha\leq\delta_t^2,\beta\leq\delta_t^2$. Next, given observation $y^t$, we look at  {the} $k$'s with finite $\rho_{tk}$ (that is, at $k$'s from $\L_t(\delta_t)$) and check whether for at least one of these $k$'s the relation $\phi_{tk}(y^t)<\alpha$ is satisfied. If this is the case, we terminate and claim that the input is a signal, otherwise we claim that so far, the nuisance hypothesis seems to be true, and pass to time $t+1$ (if $t<d$) or terminate (when $t=d$).
\subsubsection{Characterizing performance}
The performance of the above  inference procedure can be described as follows (cf. Proposition \ref{prop16}):
\begin{proposition}\label{prop17} Let the observation noise $\xi\sim\N(0,\Theta)$ with $\Theta\in\U$. Then
\par{\rm (i)} when the input is a feasible  {nuisance}, the probability to terminate with the signal conclusion at time $t\in\{1,...,d\}$ does not exceed $\epsilon_t$, and thus the probability
of  {a} false alarm is at most $\epsilon=\sum_{t=1}^d\epsilon_t$ {\rm ({we know this already} from Proposition \ref{prop16})}
\par{\rm (ii)} when $t\in\{1,...,d\}$ and $k\in\{1,...,K\}$ are such that $\rho_{tk}<\infty$, and the input belongs to a set $X^\rho_k$ with $\rho\geq \rho_{tk}$, then the probability to terminate at step $t$ with the signal conclusion is at least $1-\epsilon$.
\end{proposition}

\subsection{Near-optimality}\label{sect:assessing}
 Our goal now is to understand how good are the inference procedures we have developed.
 For the sake of definiteness, assume that
 $$
 \epsilon_t=\epsilon/d,\,1\leq t\leq d.
 $$
 We consider two assumptions about the observation noise $\xi$ (\ref{eqOS}) along with two respective change inference procedures:
 \begin{itemize}
 \item {\sl Sub-Gaussian case}, where $\xi$ is known to be sub-Gaussian with parameters $(0,\Theta)$ and $\Theta$ known to belong to $\U$; the corresponding inference procedure is built in Section \ref{sgcase};
 \item {\sl Gaussian case}, where $\xi\sim\N(0,\Theta)$ with $\Theta\in\U$; the corresponding inference procedure is described in Section \ref{sect:refinement}.
\end{itemize}
Let us fix time instant $t\leq d$ and signal shape $k\leq K$.
\par
Given $\epsilon\in(0,1/2)$, it may happen that ${\cal SV}_{tk}(R_k)>-\half {\ErfInv^2(\epsilon)}$. In this case, informally speaking, even the  feasible signal of shape $k$ and the largest possible magnitude $R_k$ does not
allow
to claim at time $t$ that the input is signal ``$(1-\epsilon)$-reliably.''
\begin{quote}
Indeed, denoting by $(h_*;\theta_1^*,\theta_2^*,\Theta_*)$ the saddle point of the convex-concave function (\ref{ccproblem}) with $\rho=R_k$, we have
$
\theta_1^*=\bar{A}_tz_*$ with some $z_*\in N$, $\theta_2^*=\bar{A}_t[v_*+R_kw_*]$ with $v_*\in V_k$, $w_*\in W_k$ and $v_*+R_kw_*\in X$, and
$$\half\|[\Theta^*]^{-1/2}[\theta_1^*-\theta_2^*]\|_2=\sqrt{-2{\cal SV}_{tk}(R_k)}<\ErfInv(\epsilon).$$
The latter implies that when $\xi^t\sim\N(0,\Theta_*)$ (which is possible),  there is no test which allows distinguishing via observation $y^t$ with risk $\leq\epsilon$ between the feasible nuisance input $z_*$ and the feasible signal $v_*+R_k w_*$ of shape $k$ and magnitude $\geq R_k$. In other words, even after the nuisance hypothesis is reduced to a single nuisance input $z_*$, and the alternative to this hypothesis is reduced to a single signal $v_*+R_k w_*$ of shape $k$ and magnitude $R_k$, we are still unable to distinguish $(1-\epsilon)$-reliably between these two hypotheses via observation $y^t$ available at time $t$.
\end{quote}
Now consider the situation where
\begin{equation}\label{assum1}
{\cal SV}_{tk}(R_k)\leq-\half {\ErfInv^2(\epsilon)},
\end{equation}
so that there exists $\rho_{tk}^*\in(0,R_k)$ such that
\begin{equation}\label{rhostar}
{\cal SV}_{tk}(\rho_{tk}^*)=-\half {\ErfInv^2(\epsilon)}.
\end{equation}
Similarly to the above, $\rho_{tk}^*$ is  {just} the smallest magnitude of signal of shape $k$ which is distinguishable from nuisance at time $t$, meaning that for every $\rho'<\rho_{tk}^*$ there
exist a feasible nuisance input $u$ and feasible signal input of shape $k$ and magnitude $\geq\rho'$ such that these two inputs cannot be distinguished via $y^t$ with risk $\leq \epsilon$. A natural way to quantify the quality of an inference procedure is to look at the smallest magnitude $\rho$ of  a feasible signal  of shape $k$  which, with probability $1-\epsilon$, ensures the {signal} conclusion and termination at time $t$. We can quantify the performance of a procedure by the ratios $\rho/\rho_{tk}^*$ stemming from various $t$ and $k$, the closer these ratios are to 1, the better. The result of this quantification of the inference procedures we have developed is as follows:
\begin{proposition}\label{prop40} Let $\epsilon\in(0,1/2)$, $t\leq d$ and $k\leq K$ be such that {\rm (\ref{assum1})}  {is satisfied. Let} $\epsilon_t=\epsilon/d$, $1\leq t\leq d$, and let $\rho_{tk}^*\in(0,R_k)$ be given by {\rm (\ref{rhostar})}. Let, further, a real $\chi$ satisfy
$\chi>\underline{\chi}$ where
\begin{equation}\label{chi}
\underline{\chi}=\left\{
\begin{array}{ll}{1\over 2}\left[1+{\ErfInv(\epsilon/(Kd))\over\ErfInv(\epsilon)}\right],&\hbox{Gaussian case}\\
& \hbox{$\U$ contains $\succeq$-largest element,}\\
\left({1\over 2}\left[1+{\ErfInv(\epsilon/(Kd))\over\ErfInv(\epsilon)}\right]\right)^2,&\hbox{Gaussian case}\\&\hbox{$\U$ does not contain $\succeq$-largest element,}\\
{\ln(Kd/\epsilon^2)\over\ErfInv^2(\epsilon)},&\hbox{sub-Gaussian case.}
\end{array}\right.
\end{equation}
Then, whenever the input is a feasible signal of shape $k$ and magnitude at least $\chi\rho_{tk}^*$, the probability for the inference procedure from Section \ref{sgcase} in the sub-Gaussian case, and the procedure from  Section \ref{sect:refinement} in the Gaussian case, to terminate at time $t$ with the signal inference is at least $1-\epsilon$.
\end{proposition}

\paragraph{Discussion.} Proposition \ref{prop40} states that when (\ref{assum1}) holds (which, as was explained, just says that feasible signals of shape $k$ of the largest possible magnitude $R_k$  can be $(1-\epsilon)$-reliably detected at time $t$), the ratio $\chi$ of
 {the} magnitude of {a signal} of shape $k$ which is detected $(1-\epsilon)$-reliably by the inference procedure we have developed to the lower bound $\rho_{tk}^*$ on the magnitude of activation of shape $k$ detectable $(1-\epsilon)$-reliably at time $t$ by {\sl any other} inference procedure can be made arbitrarily close to the right hand side quantities in (\ref{chi}). It is immediately seen that the latter quantities are upper-bounded by $\bar{\chi}=O(1)\ln(Kd/\epsilon)/\ln(1/\epsilon)$, provided $\epsilon\leq 0.5$. We see that {\sl unless $K$ and/or $d$ are  extremely large, $\bar{\chi}$ is a moderate constant}. Moreover, when $K$, $d$ remain fixed and $\epsilon\to+0$, we have $\bar{\chi}\to 1$, which, informally speaking, means that {\sl with $K,d$ fixed, the performance of the inference routines in this section
approaches the optimal performance as $\epsilon\to+0$.}
\subsection{Numerical illustration}\label{numill}
The setup of the numerical experiment we are about to report upon is as follows. We observe on time horizon $\{t=1,2,...,d=16\}$ the output $z_1,z_2,...$ of the dynamical system
\begin{equation}\label{finiteD}
(I-\Delta)^3 z=\kappa (I-\Delta)^2(u+\zeta),\,\,\kappa=(0.1d)^{-3}\approx 0.244,
\end{equation}
where $\Delta$ is the shift in the space of two-sided sequences: $(\Delta z)_t=z_{t-1}$, $\{u_t:-\infty< t<\infty\}$ is the input, and $\zeta=\{\zeta_t<-\infty<t<\infty\}$ is the random input noise with zero mean independent Gaussian components $\zeta_t$ with variances
varying in $[\sigma^2,1]$, with some given $\sigma\in(0,1]$. Our goal is to
 {dynamically test}
the nuisance hypothesis about system's input vs. a signal alternative. We start with  specifying the model of the system input.
Note that, aside from noise and the system input $u^d=[u_1;...;u_d]$ on the time horizon we are interested in, the observed output $[z_1;...;z_d]$ depends on the past -- prior to time instant $t=1$ -- outputs and inputs. The influence of this past on the observed  behavior of the system can be summarized by the initial conditions $v$ (in the case of the dynamics described by \rf{finiteD}, $v\in\bR^3$). We could augment the  input $u^d$  by these initial conditions to consider as the input the pair $x=[v;u^d]$, and express our hypotheses on input in terms of $x$, thus bringing the situation back to that considered in Section \ref{ADSetUp}. It turns out, however, that
when no restrictions  {are imposed on the initial conditions}, our inferential procedure may become numerically unstable. On the other
hand, note that by varying  {the} initial conditions we  shift the trajectory $z^t=[z_1;...;z_t]$ along the low-dimensional linear subspace $E_t\in\bR^t$ (in  {the} case of \rf{finiteD} $E_t$ is the space of collections $(z_1,...,z_s)$ with entries $z_s$ quadratically depending on $s$). Given $t$, we can
project  {the} observed $z^t$ onto the orthogonal complement of $E_t$ in $\bR^t$ and treat this projection, $y^t$, as the observation we have at time $t$. It is immediately seen that the resulting observation scheme is of the form
(\ref{eqOSt}):
\begin{equation}\label{eqOStnew}
y^t=\bar{A}_tu^t+\xi^t,\,\,4\leq t\leq d,
\end{equation}
with matrix $\bar{A}_t$ readily given by $t$, and zero mean Gaussian noise $\xi^t$ with covariance matrix $\Theta$ belonging to the ``matrix interval'' $\U_t=\{\sigma^2\Theta_t\preceq \Theta\preceq \Theta_t$\}, with $\Theta_t=\bar{A}_t\bar{A}_t^T$. Note that the restriction $t\geq4$ reflects the fact that for $t\leq 3$, $E_t=\bR^t$, and thus our observations $z^t$, $t\leq 3$, bear no information on the input $u^d$.
\par
Now we have reduced the problem to the framework of Section \ref{ADSetUp}, with inputs to the system being the actual
external inputs $u^d=[u_1;...;u_d]$ on the observation horizon. In our experiments, the
 nuisance and signal inputs were as follows:
\par$\bullet$ The set $X$ of all allowed inputs was
$$
X=\{u\in\bR^d:\|u\|_\infty\leq R=10^4\};
$$
\par$\bullet$ The set $N$ of nuisances was just the origin: $N=\{0\}\subset\bR^d$;
\par$\bullet$ The sets $V_k$ and $W_k$, $1\leq k\leq K$, responsible for signal inputs, were as follows: the number $K$ of these sets was set to $d=16$,
and we used $V_k=\{0\}$, $k\leq K$. We have considered three scenarios for the sets $W_k$ of ``activations of shape $k$ and magnitude at least $1$:''
\begin{enumerate}
\item {[pulse]} $W_k=\{u\in\bR^d: u_t=0,t\neq k,u_k\geq 1\}$, $1\leq k\leq K=d$;
\item {[jump up]} $W_k=\{u\in\bR^d: u_t=0,t<k,u_t\geq 1, k\leq t\leq d\}$, $1\leq k\leq K=d$,
\item {[step]} $W_k=\{u\in\bR^d: u_t=0,t<k,u_k=u_{k+1}=...=u_d\geq 1\}$, $1\leq k\leq K=d$.
\end{enumerate}
In other words, in our experiments, signals of shape $k$ are exactly the same as ``pure activations'' of shape $k$ --
these are the sequences $u_1,...,u_d$ which ``start'' at time $k$ (i.e., $u_t=0$ for $t<k$), of magnitude which is the value of $u_k$. In addition, there are some restrictions, depending on the scenario in question, on $u_t$'s for $t>k$.
\par
 {In this situation}, the detection problem becomes a version of the standard problem of detecting sequentially a pulse of a given form in the (third) derivative of a time series observed in Gaussian noise.
The goal of our experiment was to evaluate  {the performance of the inference procedure from Section \ref{sect:refinement} for this example}. The procedure was tuned to the probability of false alarm $\epsilon=0.01$, equally distributed between the $d=16$ time instants, that is, we used $\epsilon_t=0.01/16$, $t=1,...,d=16$.
\par
\begin{figure}
$$
\begin{array}{cc}
\epsfxsize=170pt\epsfysize=150pt\epsffile{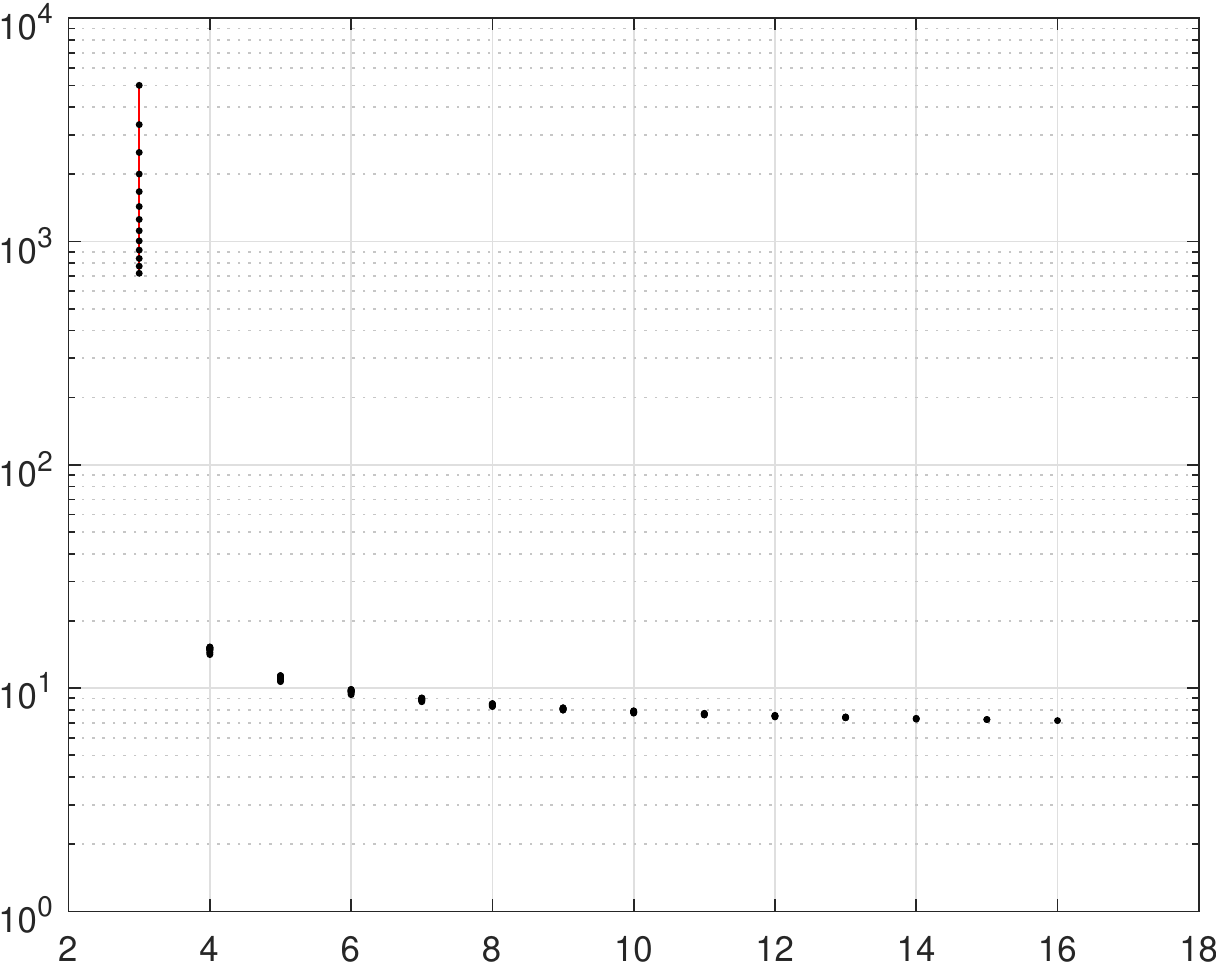}&
\epsfxsize=170pt\epsfysize=150pt\epsffile{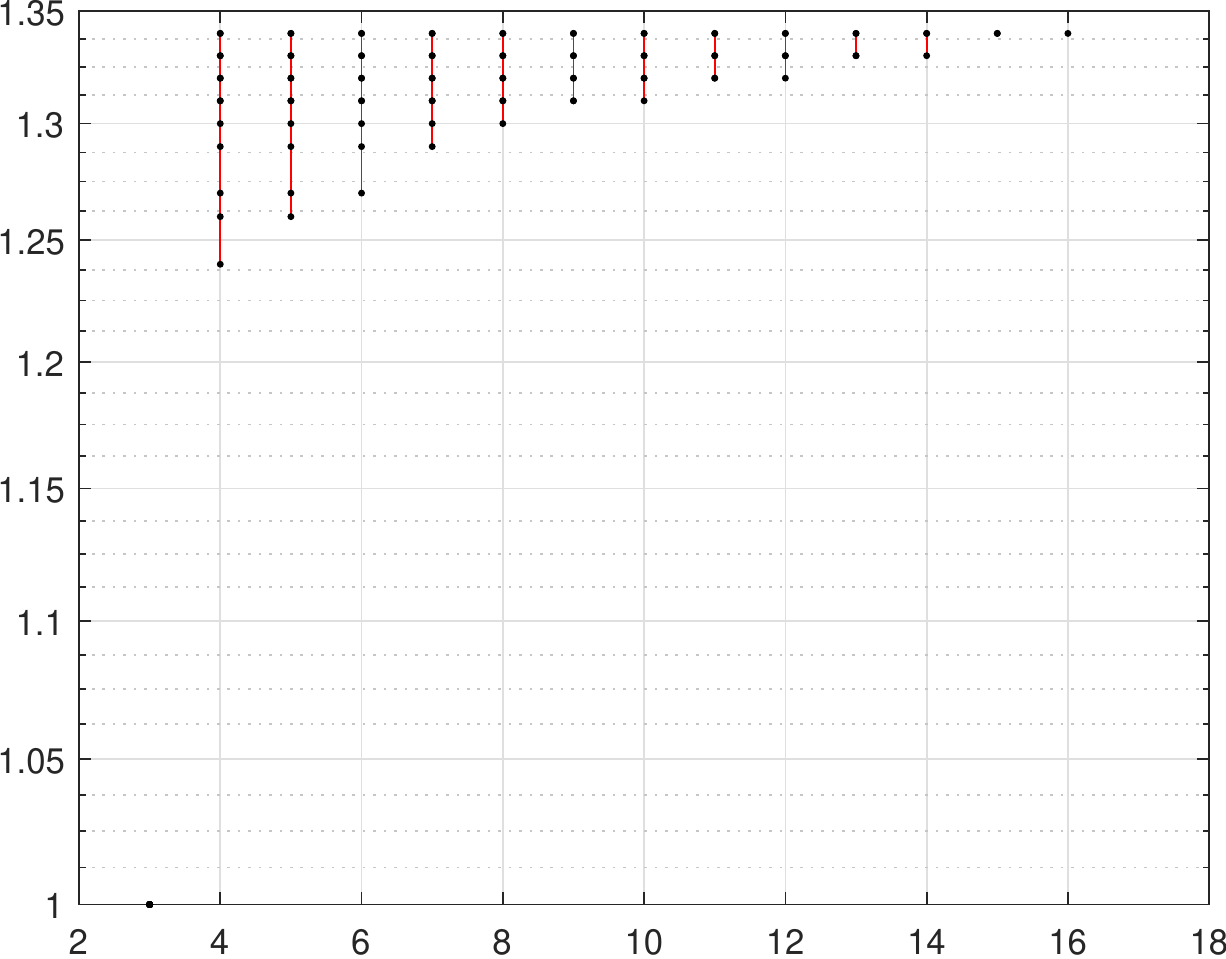}\\
\epsfxsize=170pt\epsfysize=150pt\epsffile{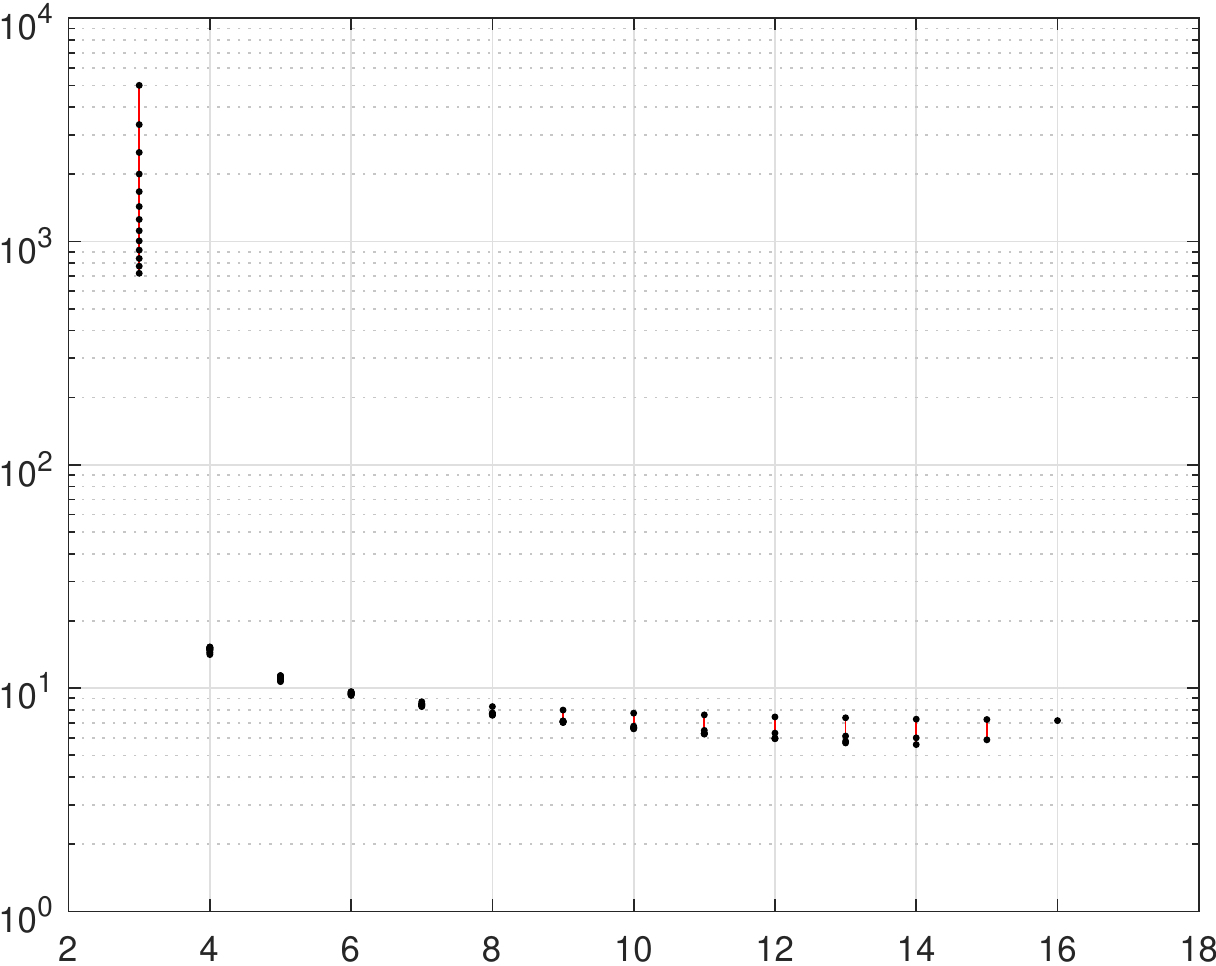}&
\epsfxsize=170pt\epsfysize=150pt\epsffile{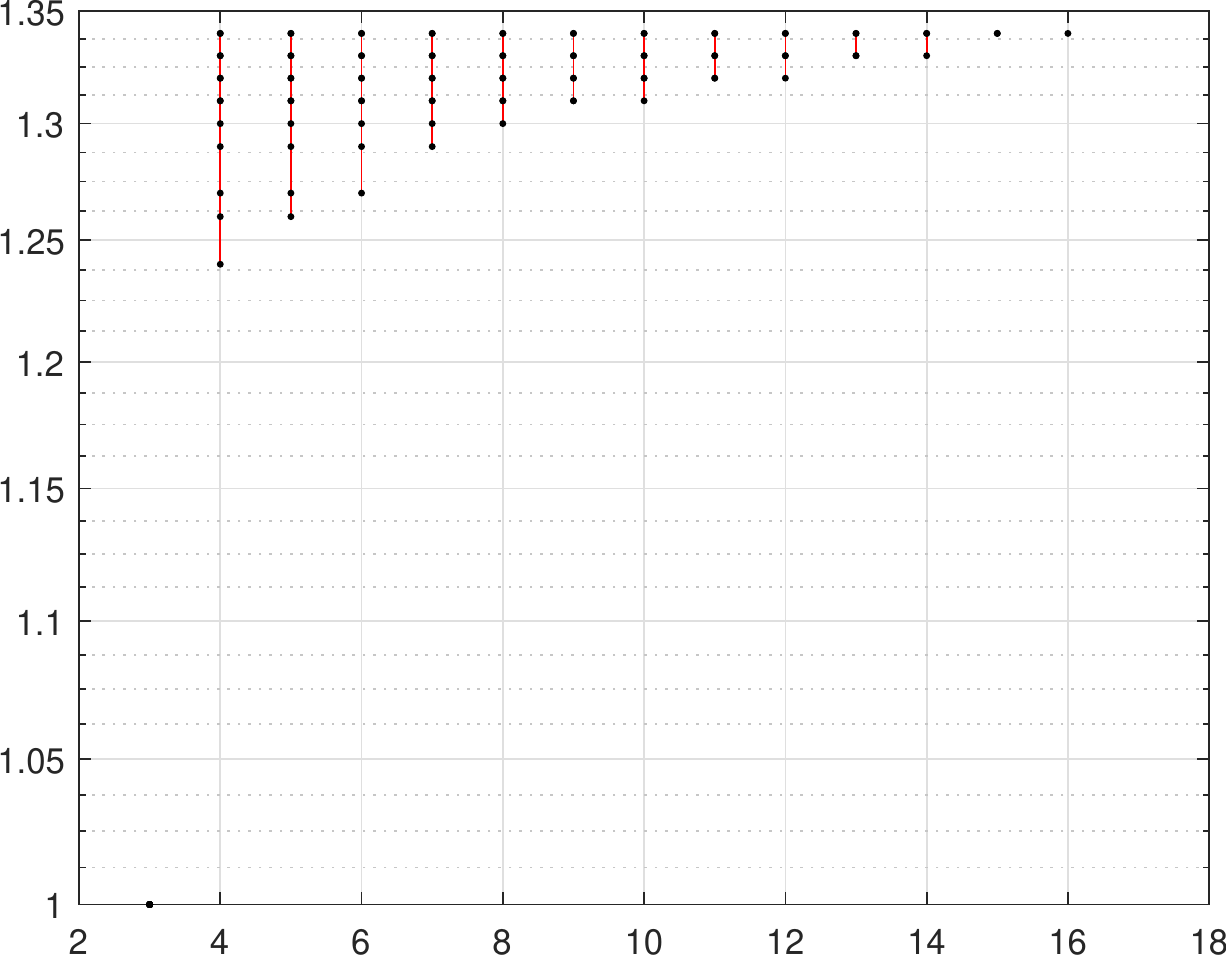}\\
\epsfxsize=170pt\epsfysize=150pt\epsffile{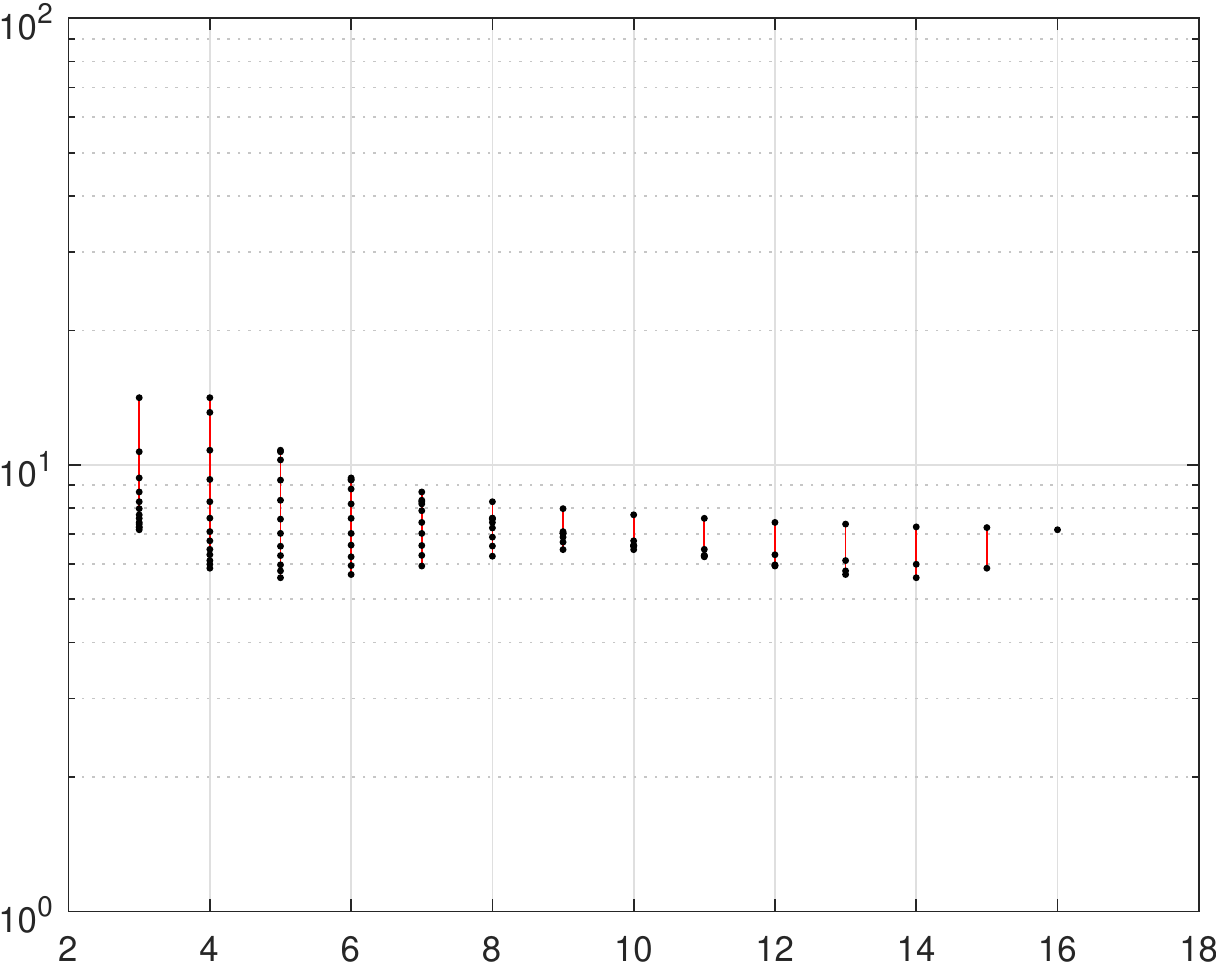}&
\epsfxsize=170pt\epsfysize=150pt\epsffile{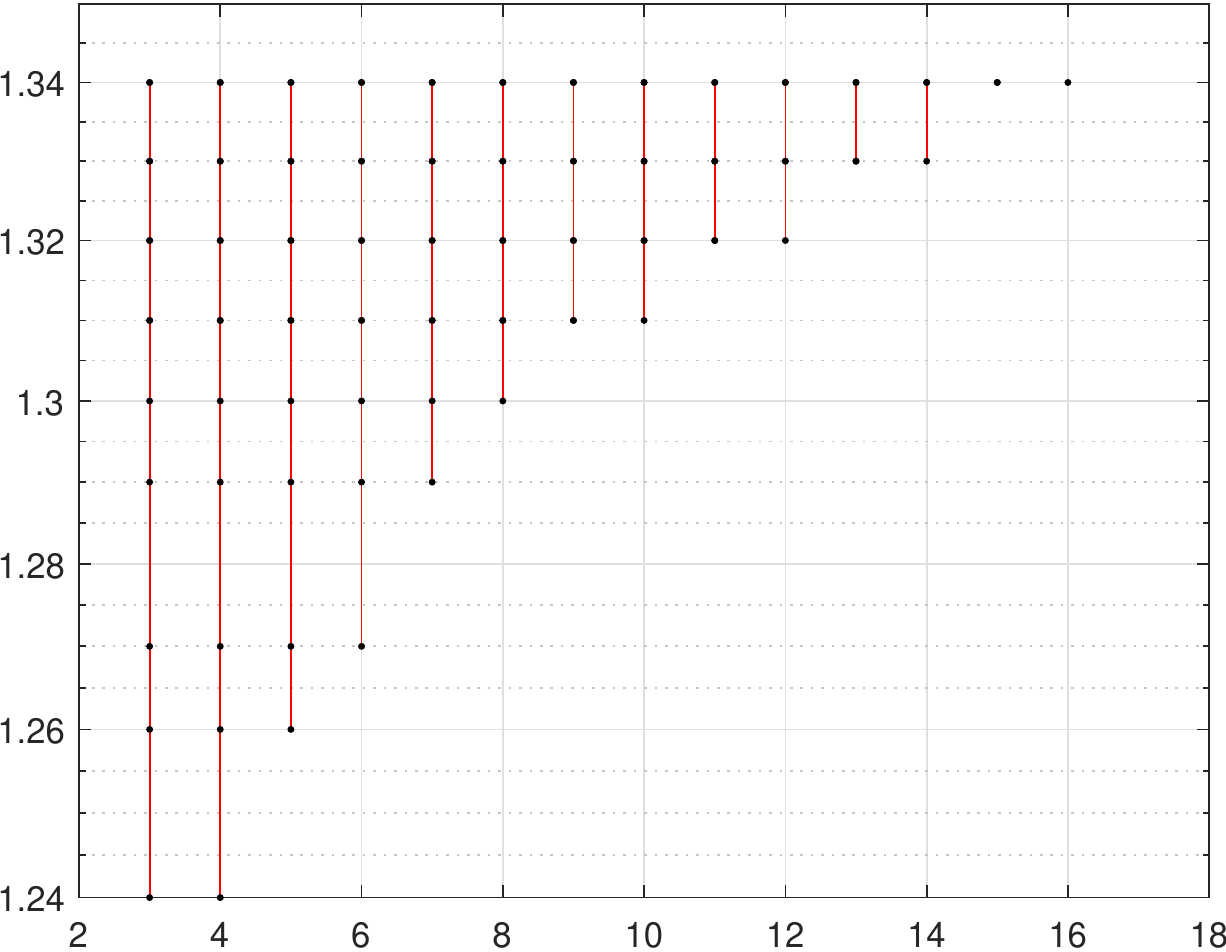}
\end{array}
$$

\caption{\label{fig111}. Performance of detector from Section \ref{sect:refinement}, dynamics (\ref{finiteD}). Left: $\rho_{tk}$, $\max[k,4]\leq t\leq 16$ (ranges and values) vs. $k$, $3\leq k\leq 16$.
Right:  $\rho_{tk}/\rho_{tk}^*$, $\max[k,4]\leq t\leq 16$ (ranges and values) vs $k$, $3\leq k\leq 16$.
Activation geometry: pulses  for top plots, jumps up for middle plots, and steps for bottom plots.
}
\end{figure}

We present the numerical results in Figure \ref{fig111}. We denote by $\rho_{tk}$  the magnitude of an activation of shape $k$ which is provably detected at time $k$ with confidence level $0.99$; we also denote
by $\rho^*_{tk}$ the ``oracle'' lower bound on this quantity.\footnote{$\rho_{tk}^*$ (defined in Section \ref{sect:assessing}) is the minimal magnitude of activation  of shape $k$ such that the ``ideal'' inference which knows $k$ in advance, tuned for
reliability $0.99$, terminates with a signal conclusion at time $t$. When $\rho_{tk}^*>R$, the maximal allowed activation magnitude, we set
$\rho^*_{tk}=+\infty$.
Recall that in the reported experiments $R=10^4$ is used.} Figure \ref{fig111} displays the dependence of $\rho_{tk}$ (left plots) and the ratio $\rho_{tk}/\rho^*_{tk}$ (right plots) on $k$ (horizontal axis) for different activation geometries (pulses, jumps up, and steps). We display  these data only for the pairs $t,k$ with finite $\rho^*_{tk}$; recall that $\rho^*_{tk}=\infty$ means that with the upper bound $R=10^4$ on the uniform norm of a feasible input, even the ideal inference does not
allow us to detect 0.99-reliably  an activation of shape $k$ at time $t$.\par Our experiment shows that $\rho^*_{tk}$ is finite in the domain $\{(t,k):\, 4\leq t\leq 16,\,3\leq k\leq t\}$. The restriction  $k\leq t$ is quite natural: we cannot detect a signal of shape $k$ before the corresponding activation starts. Note that signals of shapes $k=1,2$ are ``undetectable,'' and that no signal inputs can be detected at time $t=3$ seemingly due to the fact that activation can be completely masked by the initial conditions in the case of ``early'' activation and/or short observation horizon. Our experiment shows that this phenomenon affects equally the inference routines from Sections \ref{sgcase} and \ref{sect:refinement}, and the ideal detection, and disappears when the initial conditions for (\ref{eqOStnew}) are set to 0 and our inferences are adjusted to this a priori information.
\par
The data in Figure \ref{fig111} show that the ``non-optimality ratios'' $\rho_{tk}/\rho^*_{tk}$ of the proposed inferences as compared to the ideal detectors are quite moderate -- they never exceed 1.34;  not that bad, especially when taking into account that the ideal detection assumes {\em a priori} knowledge of the activation shape (position).

\subsection{Extension: union-type nuisance}\label{sect:exten}
So far, we have considered the case of a single nuisance hypothesis and multiple signal alternatives. The proposed approach can be easily extended to the case of multiple nuisance hypotheses, namely, to the situation differing from the one described in Section \ref{ADSetUp} in exactly one point -- instead of assuming that nuisances belong to a closed convex set $N\subset X$, we can assume that nuisance inputs run through the union $\bigcup_{m=1}^MN_m$ of given closed convex sets $N_m\subset X$, with $0\in N_m$ for all $m$. The implied modifications of our  constructions and results are as follows.
\paragraph{Sub-Gaussian case.} In this case, the construction of Section \ref{impl:constr}  in \cite{PartI}, as applied to $N_m$ in the role of $N$, gives rise to  $M$ functions
\begin{equation}\label{ccproblemm}
{\cal SV}_{tk}^m(\rho)=\min\limits_{h\in\bR^{\nu_t}}\max\limits_{\theta_1\in N^{tm},\theta_2\in U^t_{k\rho},\Theta\in \U_t}
\left[\half h^T[\theta_2-\theta_1]+\half h^T\Theta h\right],\,\,\,N^{tm}=\bar{A}_tN_m,
\end{equation}
$1\leq m\leq M$, and thus - to the parametric families
\begin{equation}\label{Ktdeltam}
\begin{array}{rcl}
\K_t(\varkappa)&=&\{k\leq K: {\cal SV}_{tk}^m(R_k)< \ln(\varkappa),\,\forall 1\leq m\leq M\},\\
K_t(\varkappa)&=&\Card\K_t(\varkappa),
\end{array}
\end{equation}
so that $K_t(\varkappa)$ is nondecreasing and continuous from the left.
At time instant $t$ we act as follows:
\begin{enumerate}
\item We define the quantity
$$
\varkappa_t=\sup\left\{\varkappa\in(0,1]:MK_t(\varkappa)\leq{\epsilon\epsilon_t\over\varkappa^2}\right\}.
$$
Clearly,
$\varkappa_t$ is well defined, takes values in $(0,1)$, and since $K_t(\varkappa)$ is continuous from the left, we have
\begin{equation}\label{eq700m}
MK_t(\varkappa_t)\leq {\epsilon\epsilon_t\over\varkappa_t^2}.
\end{equation}
We set
$
\alpha_t=\ln(\varkappa_tM/\epsilon).
$
\item
For $k\in\K_t(\varkappa_t)$, we have $0=\lim_{\rho\to+0}{\cal SV}_{tk}^m(\rho)>\ln(\varkappa_t)$ and ${\cal SV}_{tk}^m(R_k)<\ln(\varkappa_t)$. Invoking Lemma \ref{lem1}, there exists (and can be rapidly approximated to high accuracy by bisection) $\rho_{tk}\in(0,R_k)$ such that
\begin{equation}\label{balancem}
\max_{m\leq M}{\cal SV}_{tk}^m(\rho_{tk})=\ln(\varkappa_t).
\end{equation}
Given $\rho_{tk}$, we define the affine detectors
$$
\phi_{tk}^m(y^t)=h_{tkm}^T[y^t-w_{tkm}],\,\,w_{tkm}=\half[\theta_{1,tkm}+\theta_{2,tkm}],
$$
where $(h_{tkm};\theta_{1,tkm},\theta_{2,tkm},\Theta_{tkm})$ is a solution to the saddle point problem
(\ref{ccproblemm}) with $\rho=\rho_{tk}$.
\par
For $k\not\in\K_t(\varkappa_t)$, we set $\rho_{tk}=+\infty$.
\item\label{fffinal}
Finally, we process the observation $y^t$ at step $t$ as follows:
\begin{itemize}
\item if there exists $k$ such that $\rho_{tk}<\infty$ and  $\phi_{tk}^m(y^t)<\alpha_t$ for all $m\leq M$, we claim that the observed input is a signal, and terminate;
\item otherwise, we claim that so far, the nuisance hypothesis is not rejected, and pass to the next time instant $t+1$ (when $t<d$) or terminate (when $t=d$).
\end{itemize}
\end{enumerate}
 {The performance} of the inference policy we have described is given by the following analogue of Proposition \ref{prop16}:
\begin{proposition}\label{prop16m}  {For any} zero mean sub-Gaussian, with parameter $\Theta\in\U$,  distribution of observation noise on time horizon $1,...,d$, 
\begin{itemize}
\item when the input is a nuisance (i.e., belongs to $\cup_{m\leq M}N_m$), the probability to terminate with the signal conclusion at time $t\in\{1,...,d\}$ does not exceed $\epsilon_t$, and thus the probability of  {a} false alarm is at most $\epsilon=\sum_{t=1}^d\epsilon_t$;
\item when $t\in\{1,...,d\}$ and $k\in\{1,...,K\}$ are such that $\rho_{tk}<\infty$, and the input belongs to a set $X^\rho_k$ with $\rho\geq \rho_{tk}$, then the probability to terminate at step $t$ with the  signal conclusion is at least $1-\epsilon$.
\end{itemize}
Furthermore, let us assume that $\epsilon_t=\epsilon/d$ for all $t$ and that for some $t\leq d$ and $k\leq K$ we have
\begin{equation}\label{eq393}
{\cal SV}_{tk}^m(R_k)<-\half\ErfInv^2(\epsilon),\,\,\forall 1\leq m\leq M,
\end{equation}
so that the quantities $\rho_{tk}^{m*}\in(0,R_k)$ such that
$$
{\cal SV}_{tk}^m(\rho_{tk}^{m*})=-\half\ErfInv^2(\epsilon)
$$
are well defined (for ``lower bound interpretation'' of these quantities, see comments after (\ref{rhostar})). Then for every $\chi$ satisfying
$$
\chi>{\ln(dKM/\epsilon^2)\over \ErfInv^2(\epsilon)},
$$
and every feasible signal input of shape $k$ and magnitude $\geq \chi\max\limits_{m\leq M} \rho_{tk}^{m*}$, the probability of termination with signal conclusion at time $t$ is $\geq 1-\epsilon$.
\end{proposition}
Proof of Proposition \ref{prop16m} is given by a straightforward modification of the proofs of Propositions \ref{prop16} and \ref{prop40}.

\section{Change detection via quadratic detectors}\label{QDSection}
\subsection{Outline}
In Section \ref{ChPADSetUp}, we were interested in deciding as early as possible upon the hypotheses about the input $x$ underlying observations \rf{eqOSt} in the situation where both signals and nuisances
formed finite unions of convex sets. Solving this problem was reduced to decisions on pairs of {\sl convex} hypotheses -- those stating
that the expectation of a  (sub-)Gaussian random vector with partly known covariance matrix belongs to the union of convex sets associated with the hypotheses, and we could make decisions looking at the (signs of) properly
 {built} affine detectors -- affine functions of observations. Now we intend to address
the case when the signals (or nuisances) are specified by {\sl non-convex} restrictions, such as ``$u$ belongs to a given linear subspace and has Euclidean norm at least $\rho>0$.''
This natural setting is difficult to capture via convex hypotheses: in such an attempt, we are supposed to ``approximate'' the restriction ``the $\|\cdot\|_2$-norm of vector $x$ is $\geq \rho$'' by the union of convex hypotheses like ``$i$-th entry in $x$ is $\geq\rho'$''/``$i$-th entry in $x$ is $\leq-\rho'$''; the number of these hypotheses grows with the input's dimension, and the ``quality of approximation,''
whatever be its definition, deteriorates as the dimension grows.
\par
In this situation, a natural way to proceed is to look at ``quadratic liftings'' of inputs and observations. Specifically, given a vector $w$ of dimension $m$, let us associate with it its ``quadratic lifting'' -- the symmetric $(m+1)\times (m+1)$ matrix $Z(w)=[w;1][w;1]^T$. Observe that the restrictions on $w$ expressed by linear and quadratic constraints induce {\sl linear} restrictions
on $Z(w)$.
Secondly, given noisy observation $y^t=\bar{A}_tx+\xi^t$ of signal $x$, the quadratic lifting $Z(y^t)$ can be thought of as noisy observation of
an {\sl affine image} $\widehat{A}Z(x)\widehat{A}^T$ of $Z(x)$, where \[\widehat{A}=\left[\begin{array}{c|c}\bar{A}_t&\cr\hline&1\cr\end{array}\right],\]
(here and in what follows the empty block refers to the null matrix). As a result, roughly speaking, linear and quadratic constraints on the input translate into linear constraints on the expectation of ``lifted observation'' $Z(y^t)$, and different hypotheses on input, expressed by linear and quadratic constraints, give rise to convex hypotheses on $Z(y^t)$. Then, in order to decide on the resulting
convex hypotheses, we can use affine in $Z(y^t)$, that is, {\sl quadratic} in $y^t$, detectors, and this is what we intend to do.
\subsection{Preliminaries}
\subsubsection{Gaussian case}
In the sequel, the following result (which is a slightly modified concatenation of Propositions 3.1 and 5.1 of \cite{PartI}) is used:

\begin{proposition}\label{concatenation}
 \item[$\qquad${\rm (i)}] Let $\U$ be a convex compact set contained in the interior of the cone $\bS^\nu_+$ of positive semidefinite $\nu\times\nu$ matrices in the space $\bS^\nu$ of symmetric $\nu\times\nu$ matrices.
 Let $\Theta_*\in \bS^\nu_+$ be such that $\Theta_*\succeq \Theta$ for all $\Theta\in\U$, and let $\delta\in[0,2]$ be such that
\begin{equation}\label{56delta}
\|\Theta^{1/2}\Theta_*^{-1/2}-I_\nu\|\leq\delta,\;\;\forall \Theta\in\U,
\end{equation}
where $\|\cdot\|$ is the spectral norm.\footnote{with $\delta=2$, (\ref{56delta}) is satisfied for all $\Theta$ such that $0\preceq\Theta\preceq \Theta_*$, so that the restriction $\delta\leq2$ is w.l.o.g.} Finally,
let $\gamma\in(0,1)$,  $A$ be a $\nu\times(n+1)$ matrix,  $\Z$ be a nonempty convex compact subset of the set $\Z^+=\{Z\in\bS^{n+1}_+:Z_{n+1,n+1}=1\}$, and let
\begin{equation}\label{phiZ}
\phi_\Z(Y):=\max_{Z\in\Z} \Tr(ZY)
\end{equation}
be the support function of $\Z$. These data specify the closed convex set
\begin{equation}\label{PcH}
\H=\H^\gamma:=\{(h,H)\in \bR^\nu\times \bS^\nu:-\gamma\Theta_*^{-1}\preceq H\preceq \gamma \Theta_*^{-1}\},
\end{equation}
the matrix
\begin{equation}\label{matrixB}
B=\left[\begin{array}{c}A\cr [0,...,0,1]\cr\end{array}\right]
\end{equation}
and the function $\Phi_{A,\Z}:\,\H\times\U\to\bR$,
\begin{equation}\label{phi}
\begin{array}{rcl}
\Phi_{A,\Z}(h,H;\Theta)&=&-\half\ln\Det(I-\Theta_*^{1/2}H\Theta_*^{1/2})+\half \Tr([\Theta-\Theta_*]H)\\
&&+{\delta(2+\delta)\over 2(1-\|\Theta_*^{1/2}H\Theta_*^{1/2}\|)}\|\Theta_*^{1/2}H\Theta_*^{1/2}\|_F^2\\
&&+{1\over 2}\phi_\Z\left(B^T\left[\hbox{\small$\left[\begin{array}{c|c}H&h\cr\hline h^T&\end{array}\right]+
\left[H,h\right]^T[\Theta_*^{-1}-H]^{-1}\left[H,h\right]$}\right]B\right)
,
\end{array}
\end{equation}
where $\|\cdot\|_F$ is the Frobenius norm of a matrix.
\par
Function $\Phi_{A,\Z}
$ is continuous on its domain, convex in $(h,H)\in\H$ and concave in $\Theta\in\U$ and possesses the following property:
\begin{quotation}
\noindent
Whenever $u\in \bR^n$ is such that $[u;1][u;1]^T\in\Z$ and $\Theta\in\U$, the Gaussian random vector $\zeta\sim\N(A[u;1],\Theta)$ satisfies the relation
\begin{equation}\label{moments}
\forall (h,H)\in\H: \;\;\ln\left(\bE_{\zeta\sim \N(A[u;1],\Theta)}\left\{{\rm e}^{\half\zeta^TH\zeta+h^T\zeta}\right\}\right)\leq \Phi_{A,\Z}(h,H;\Theta).
\end{equation}
\end{quotation}
Besides this, $\Phi_{A,\Z}$ is coercive in $(h,H)$: $\Phi_{A,\Z}(h_i,H_i;\Theta)\to+\infty$ as $i\to\infty$ whenever $\Theta\in \U$,  $(h_i,H_i)\in\H$ and $\|(h_i,H_i)\|\to\infty$, $i\to\infty$.
\item[$\qquad${\rm (ii)}] Let  two collections of data from {\rm (i):} $(\U_\chi,\Theta_*^{(\chi)},\delta_\chi,\gamma_\chi,A_\chi,\Z_\chi)$, $\chi=1,2$,   with common $\nu$ be given, giving rise to the sets $\H_\chi $, matrices $B_\chi $, and functions $\Phi_{A_\chi ,\Z_\chi }(h,H;\Theta)$, $\chi=1,2$. These collections
specify the families of normal distributions
\[
\G_\chi =\{\N(v,\Theta): \Theta\in\U_\chi \ \&\ \exists u: v=A_\chi [u;1], [u;1][u;1]^T\in\Z_\chi \},\,\chi=1,2.
\]
Consider the convex-concave saddle point problem
\begin{equation}\label{SPPLift}
{\cal SV}=\min\limits_{(h,H)\in\H_1\cap\H_2}\max\limits_{\Theta_1\in\U_1,\Theta_2\in\U_2} \underbrace{\half\left[\Phi_{A_1,\Z_1}(-h,-H;\Theta_1)+
\Phi_{A_2,\Z_2}(h,H;\Theta_2)\right]}_{\Phi(h,H;\Theta_1,\Theta_2)}.
\end{equation}
A saddle point  {$(h_*, H_* ;\Theta_1^*,\Theta_2^*)$}  does exist in this problem, and the induced quadratic detector
\begin{equation}\label{iquaddet}
\phi_*(\omega)=\half\omega^TH_*\omega+h_*^T\omega+\underbrace{\half\left[\Phi_{A_1,\Z_1}(-h_*,-H_*;\Theta^*_1)
-\Phi_{A_2,\Z_2}(h_*,H_*;\Theta^*_2)\right]}_{a},
\end{equation}
satisfies
\begin{equation}\label{riskagain}\begin{array}{llll}
(a)&~~\int_{\bR^\nu}{\rm e}^{-\phi_*(\omega)}P(d\omega)&\leq \epsilon_\star:={\rm e}^{{\cal SV}}\;\;&\forall P\in \G_1,\\
(b)&~~\int_{\bR^\nu}{\rm e}^{\phi_*(\omega)}P(d\omega)&\leq \epsilon_\star\;\;&\forall P\in \G_2.
\end{array}
\end{equation}
That is, the risk, as defined in item \itemdrisk\ of Section \ref{sectnotconv}, of the detector $\phi_*$ on the families $\G_1,\G_2$ satisfies
$$
\Risk(\phi_*|\G_1,\G_2)\leq\epsilon_\star.
$$
\end{proposition}
For the proof, see \cite{PartI}; for the reader's convenience, we reproduce the proof in Section \ref{AppConcProof}. The justification for the remark below can be found in appendix \ref{just_remark41}.

\begin{remark}\label{remark41} {\rm Note that the computational effort of solving {\rm (\ref{SPPLift})} reduces dramatically in the ``{\em easy case}'' of the situation described in  item (ii) of Proposition \ref{concatenation}, specifically,
in the case where
\begin{itemize}
\item the observations are {\sl direct}, meaning that $n=\nu$ and $A_\chi [u;1] \equiv u$, $u\in\bR^\nu$, $\chi=1,2$;
\item the sets $\U_\chi$ are comprised of positive definite {\sl diagonal} matrices, and the matrices $\Theta^{(\chi)}_*$ are diagonal as well, $\chi=1,2$;
\item the sets $\Z_\chi$, $\chi=1,2$, are convex compact sets of the form
$$
\Z_\chi=\{Z\in\bS^{\nu+1}_+: Z\succeq0,\,\Tr(ZQ^\chi_j)\leq q^\chi_j,\,1\leq j\leq J_\chi\}
$$
with {\sl diagonal} matrices $Q^\chi_j$, \footnote{In terms of the sets $U_\chi$, this assumption means that the latter sets
are given by linear inequalities on the {\sl squares} of entries in $u$.} and these sets intersect the interior of the positive semidefinite cone $\bS^{\nu+1}_+$.
\end{itemize}
In this case, the convex-concave saddle point problem {\rm (\ref{SPPLift})} admits a saddle point $(h_*,H_*;\Theta_1^*,\Theta_2^*)$
where $h_*=0$ and $H_*$ is diagonal, and restricting $h$ to be zero and $H$ to be diagonal reduces drastically the design dimension of the saddle point problem.
}
\end{remark}

\subsubsection{Sub-Gaussian case}
Sub-Gaussian version of Proposition \ref{concatenation} is as follows:
\begin{proposition}\label{concatenationSG}
 \item[$\qquad${\rm (i)}] Let $\U$ be a convex compact set contained in the interior of the cone $\bS^\nu_+$ of positive semidefinite $\nu\times\nu$ matrices in the space $\bS^\nu$ of symmetric $\nu\times\nu$ matrices, let $\Theta_*\in \bS^\nu_+$ be such that $\Theta_*\succeq \Theta$ for all $\Theta\in\U$, and let $\delta\in[0,2]$ be such that {\rm (\ref{56delta})} holds true.
Finally, let $\gamma, \gamma^+$ be such that $0<\gamma<\gamma^+<1$, $A$ be $\nu\times(n+1)$ matrix,  $\Z$ be a nonempty convex compact subset of the set $\Z^+=\{Z\in\bS^{n+1}_+:Z_{n+1,n+1}=1\}$, and let $\phi_\Z(Y)$ be the support function of $\Z$, see {\rm (\ref{phiZ})}.
These data specify the closed convex sets
\[
\begin{array}{rcl}
\H &=& \H^\gamma:=\{(h,H)\in \bR^\nu\times \bS^\nu:\;-\gamma\Theta_*^{-1}\preceq H\preceq \gamma \Theta_*^{-1} \},\\
\widehat{\H} & = & \widehat{\H}^{\gamma,\gamma^+} =\left\{(h,H,G)\in \H^\gamma  \times \bS^\nu:
0\preceq G\preceq \gamma^+ \Theta_*^{-1},\,H\preceq G\right\},
\end{array}
\]
the matrix $B$ given by {\rm (\ref{matrixB})},
and the functions
\begin{equation}\label{phiSG}
\hbox{\footnotesize$\begin{array}{l}
\Psi_{A,\Z}(h,H,G)\\
\quad=-\half\ln\Det(I-\Theta_*^{1/2}G\Theta_*^{1/2})\\
~~~~~ +\half\phi_\Z\left(B^T\left[\hbox{\scriptsize$\left[\begin{array}{r|r}H&h\cr\hline h^T&\cr\end{array}\right]$}+ [H,h]^T[\Theta_*^{-1}-G]^{-1}[H,h]\right]B\right): {\widehat{\H}}\times\Z\to\bR,\\
\Psi^\delta_{A,\Z}(h,H,G;\Theta)\\
\quad=-\half\ln\Det(I-\Theta_*^{1/2}G\Theta_*^{1/2})+\half\Tr([\Theta-\Theta_*]G)+{\delta(2+\delta)\over2(1-\|\Theta_*^{1/2}G\Theta_*^{1/2}\|)}\|
\Theta_*^{1/2}G\Theta_*^{1/2}\|_F^2\\
\quad+\half\phi_\Z\left(B^T\left[\hbox{\scriptsize$\left[\begin{array}{r|r}H&h\cr\hline h^T&\cr\end{array}\right]$}+ [H,h]^T[\Theta_*^{-1}-G]^{-1}[H,h]\right]B\right):
{ {\widehat{\H}}\times\{0\preceq\Theta\preceq \Theta_*\}\to\bR,}\\
\Phi_{A,\Z}(h,H)=\min\limits_G\left\{\Psi_{A,\Z}(h,H,G):(h,H,G)\in  {\widehat{\H}}\right\}:\H\to\bR,\\
\Phi^\delta_{A,\Z}(h,H;\Theta)=\min\limits_G\left\{\Psi^\delta_{A,\Z}(h,H,G;\Theta):
(h,H,G)\in\widehat{\H}\right\}:\H\times\{0\preceq\Theta\preceq \Theta_*\}\to\bR,\\
\end{array}$}
\end{equation}
where, same as in {\rm (\ref{phi})}, $\|\cdot\|$ is the spectral, and $\|\cdot\|_F$ is the Frobenius norm of a matrix.
\par
Function $\Phi_{A,\Z}(h,H)$ is convex and continuous on its domain, while function $\Phi^\delta_{A,\Z}(h,H;\Theta)$ is continuous on its domain, convex in $(h,H)\in\H$ and concave in
$\Theta\in\{0\preceq\Theta\preceq\Theta_*  {\}}$.
Besides this,
\begin{quotation}\noindent
Whenever $u\in\bR^n$ is such that $[u;1][u;1]^T\in\Z$ and $\Theta\in\U$,
the sub-Gaussian random vector $\zeta$ with parameters $(A[u;1], \Theta)$ satisfies the relation
\begin{equation}\label{momentsSG}
\begin{array}{ll}
\multicolumn{2}{l}{\forall (h,H)\in\H:}\\
(a)&\ln\left(\bE_{\zeta}\left\{{\rm e}^{\half\zeta^TH\zeta+h^T\zeta}\right\}\right)\leq \Phi_{A,\Z}(h,H),\\
(b)&\ln\left(\bE_{\zeta}\left\{{\rm e}^{\half\zeta^TH\zeta+h^T\zeta}\right\}\right)\leq \Phi^\delta_{A,\Z}(h,H;\Theta).\\
\end{array}
\end{equation}
\end{quotation}
Besides this, $\Phi_{A,\Z}$ and $\Phi^\delta_{A,\Z}$ are coercive in $(h,H)$: $\Phi_{A,\Z}(h_i,H_i)\to+\infty$ and $\Phi^\delta_{A,\Z}(h_i,H_i;\Theta)\to+\infty$ as $i\to\infty$ whenever $\Theta\in \U$,  $(h_i,H_i)\in\H$ and $\|(h_i,H_i)\|$ $\to\infty$, $i\to\infty$.
\item[$\qquad${\rm (ii)}] Let  two collections of data from {\rm (i):} $(\U_\chi,\Theta_*^{(\chi)},\delta_\chi,\gamma_\chi, {\gamma^+_\chi},A_\chi,\Z_\chi)$, $\chi=1,2$,   with common $\nu$ be given, giving rise to the sets $\H_\chi $, matrices $B_\chi $, and functions $\Phi_{A_\chi ,\Z_\chi }(h,H)$, $\Phi^{\delta_\chi}_{A_\chi,\Z_\chi }(h,H;\Theta)$, $\chi=1,2$. These collections specify the families of distributions
$
\SG_\chi,
$ $\chi=1,2$, where $\SG_\chi$ is comprised of all sub-Gaussian distributions  with parameters $v,\Theta$, such that $v$ can be represented as $A_\chi[u;1]$ for  some $u$ with $[u;1][u;1]^T\in\Z_\chi$, and $\Theta\in \U_\chi$.
Consider the convex-concave saddle point problem
\begin{equation}\label{SPPLiftSG}
{\cal SV}=\min\limits_{(h,H)\in\H_1\cap\H_2}\max\limits_{\Theta_1\in\U_1,\Theta_2\in\U_2} \underbrace{\half\left[\Phi^{\delta_1}_{A_1,\Z_1}(-h,-H;\Theta_1)+
\Phi^{\delta_2}_{A_2,\Z_2}(h,H;\Theta_2)\right]}_{\Phi^{\delta_1,\delta_2}(h,H;\Theta_1,\Theta_2)}.
\end{equation}
A saddle point $(h_*,H_*;\Theta_1^*,\Theta_2^*)$  does exist in this problem, and the induced quadratic detector
\[
\phi_*(\omega)=\half\omega^TH_*\omega+h_*^T\omega+\underbrace{\half\left[\Phi^{\delta_1}_{A_1,\Z_1}(-h_*,-H_*;\Theta^*_1)
-\Phi^{\delta_2}_{A_2,\Z_2}(h_*,H_*;\Theta^*_2)\right]}_{a},
\]
when applied to the families of sub-Gaussian distributions $\SG_\chi $, $\chi=1,2$, has the risk $$
\Risk(\phi_*|\SG_1,\SG_2)\leq \epsilon_\star:={\rm e}^{{\cal SV}},$$ that is
\[
\begin{array}{lrl}
(a)&~~\int_{\bR^\nu}{\rm e}^{-\phi_*(\omega)}P(d\omega)\leq \epsilon_\star\;\;&\forall P\in \SG_1,\\
(b)&~~\int_{\bR^\nu}{\rm e}^{\phi_*(\omega)}P(d\omega)\leq \epsilon_\star\;\;&\forall P\in \SG_2.
\end{array}
\]
Similarly, the convex minimization problem
\begin{equation}\label{SPPLiftConvSG}
\Opt=\min\limits_{(h,H)\in\H_1\cap\H_2}\underbrace{\half\left[\Phi_{A_1,\Z_1}(-h,-H)+
\Phi_{A_2,\Z_2}(h,H)\right]}_{\Phi(h,H)}
\end{equation}
is solvable, and the induced by its optimal solution $(h_*,H_*)$ quadratic detector
\[
\phi_*(\omega)=\half\omega^TH_*\omega+h_*^T\omega+\underbrace{\half\left[\Phi_{A_1,\Z_1}(-h_*,-H_*)
-\Phi_{A_2,\Z_2}(h_*,H_*)\right]}_{a},
\]
when applied to the families of sub-Gaussian distributions $\SG_\chi $, $\chi=1,2$, has the risk $$
\Risk(\phi_*|\SG_1,\SG_2)\leq\epsilon_\star:={\rm e}^{\hbox{\scriptsize\rm Opt}},$$ so that for just defined $\phi_*$ and $\epsilon_\star$ relation {\rm (\ref{SPPLiftConvSG})} takes place.
\end{proposition}
\begin{remark}\label{remrem} {\rm
Proposition \ref{concatenationSG} offers two options for building quadratic detectors for the families $\SG_1$, $\SG_2$, those based on the saddle point of {\rm (\ref{SPPLiftSG})} and on the optimal solution to {\rm (\ref{SPPLiftConvSG})}. Inspecting the proof, the number of options can be increased to 4: we can replace any  of the functions $\Phi^{\delta_\chi}_{A_\chi,\Z_\chi}$, $\chi=1,2$  (or both these functions simultaneously) with $\Phi_{A_\chi,\Z_\chi}$. The second of the original two options is exactly what we get
when replacing both $\Phi^{\delta_\chi}_{A_\chi,\Z_\chi}$, $\chi=1,2$, with  {$\Phi_{A_\chi,\Z_\chi}$}. It is easily seen that depending on the data, {each} of these 4 options can result in the smallest risk bound. Thus, it makes sense to keep all these options in mind and to use the one which, under the circumstances, results in the best risk bound. Note that the risk bounds are efficiently computable, so that identifying the best option is easy.}
\end{remark}
\subsection{Setup}\label{newsetup}
We continue to consider the situation described in Section \ref{sect:ChPD}, but with different specifications of noise and of nuisance and signal inputs, as compared to Section \ref{ADSetUp}.
\par
We define nuisance and signal inputs as follows.
\par1. Admissible inputs, nuisance and signal alike, belong to a bounded  set $X\subset\bR^n$ containing the origin cut off $\bR^n$  by a system of quadratic inequalities:
\begin{equation}\label{seq100}
X=\{x\in\bR^n: \Tr(Q_iZ(x))\leq q_i,\,1\leq i\leq I\},\end{equation}
where $Q_i$ are $(n+1)\times (n+1)$ symmetric matrices.
We assume w.l.o.g. that the first constraint defining $X$ is $\|x\|_2^2\leq R^2$, that is, $Q_1$ is the diagonal matrix with the diagonal $1,...,1,0$, and $q_1=R^2$. We set
\begin{equation}\label{seq101}
\X=\{W\in\bS^{n+1}_+: W_{n+1,n+1}=1, \Tr(WQ_i)\leq q_i,1\leq i\leq I\},
\end{equation}
so that $\X$ is a convex compact set in $\bS^{n+1}_+$, and $Z(x)\in\X$ for all $x\in X$.
\par2. The set $N$ of nuisance inputs contains the origin and is cut off $X$  by a system of quadratic inequalities, so that
\begin{equation}\label{seq102}
N=\{x\in\bR^n: \Tr(Q_iZ(x))\leq q_i,\,1\leq i\leq I_+\}, \;I_+>I.\end{equation}
We set
\begin{equation}\label{seq101A}
\N=\{W\in\bS^{n+1}_+: W_{n+1,n+1}=1, \Tr(WQ_i)\leq q_i,\,1\leq i\leq I_+\},
\end{equation}
so that $\N\subset\X$ is a convex compact set in $\bS^{n+1}_+$, and $Z(x)\in\N$ for all $x\in N$.
\par3. Signals belonging to $X$ are
of different shapes and magnitudes, with signal of shape $k$, $1\leq k\leq K$, and magnitude $\geq 1$ defined as a vector from the set
    $$
    W_k=\{x\in\bR^n: \Tr(Q_{ik}Z(x))\leq b_{ik},\,1\leq i\leq I_k\}
    $$
    with two types of quadratic constraints:
    \begin{itemize}
    \item constraints of type A: $b_{ik}\leq0$, the symmetric matrices $Q_{ik}$ have zero North-West (NW) block of size  $n\times n$, and zero South-East (SE) diagonal entry; these constraints are just linear constraints on $x$;
    \item constraints of type B: $b_{ik}\leq0$, the only nonzeros in $Q_{ik}$  are in the NW block of size  $n\times n$.
    \end{itemize}
    We denote the sets of indices $t$ of constraints of these two types by $\Ii_k^A$ and $\Ii_k^B$ and assume that at least one of the right hand sides $b_{ik}$ is strictly negative, implying that $W_k$ is at a positive distance from the origin.
\par
We define a {signal} of shape $k$ and magnitude $\geq\rho>0$ as a vector from the set $W^\rho_k=\rho W_k$; note that
{\small $$
W^\rho_k=\{x\in\bR^n: \Tr(Q_{ik}Z(x))\leq \rho b_{ik},i\in \Ii_k^A,\Tr(Q_{ik}Z(x))\leq \rho^2 b_{ik},i\in \Ii_k^B\}.
$$}\noindent
We set
\[\begin{array}{l}
\W_k^\rho=\{W\in\bS^{n+1}_+: \,W_{n+1,n+1}=1,\, \Tr(Q_{ik}W)\leq \rho b_{ik},\,i\in \Ii_k^A,\\
\multicolumn{1}{r}{\Tr(Q_{ik}W)\leq \rho^2 b_{ik},\,i\in \Ii_k^B\},}\\
\end{array}
\]
ensuring that $Z(x)\in \W^\rho_k$ whenever $x\in W^\rho_k$. Note that  sets $\W_k^\rho$ shrink as $\rho>0$ grows due to $b_{ik}\leq0$. We assume that for small $\rho>0$, the sets $\W_k^\rho\cap \X$ are nonempty (this is definitely the case when some signals of shape $k$ and  positive magnitude are admissible inputs -- otherwise signals of shape $k$ are of no interest in our context, and we can ignore them). Since $\X$ is compact and some of $b_{ik}$ are negative, the sets $\W^\rho_k$ are empty for large enough values of $\rho$.
As a byproduct of the compactness of $\X$, it is immediately seen that there exists $R_k\in(0,\infty)$ such that $W_k^\rho$ is nonempty when $\rho\leq R_k$ and is empty when $\rho>R_k$.
\subsection{Change detection via quadratic detectors, Gaussian case}\label{sectquaddet}
In this section, we consider the situation of Section \ref{sect:ChPD}, {\sl assuming the noise $\xi^d$ in {\rm (\ref{eqOS})} to be zero mean Gaussian: $\xi^d\sim \N(0,\Theta)$.}
\subsubsection{Preliminaries}\label{ssprelim}
Given $t\leq d$, let us set
$$
A_t=[\bar{A}_t,0], B_t=\left[\begin{array}{c}A_t\cr [0,...,0,1]\cr\end{array}\right],
$$
so that the observation $y^t\in\bR^{\nu_t}$ at time $t$ is Gaussian with the expectation  $A_t[x;1]$ and covariance matrix $\Theta$ belonging to the convex compact subset $\U_t$ of the interior of the positive semidefinite cone $\bS^{\nu_t}_+$, see (\ref{eqOSt}), (\ref{thetat}).
\par
We fix $\gamma\in(0,1)$, and  $\Theta_{*,d}\in\bS^{\nu_d}_+$ such that $\Theta_{*,d}\succeq\Theta$ for all $\Theta\in\U_d$. For $1\leq t\leq d$, we
set  {$\Theta_{*,t}=S_t\Theta_{*,d}S_t^T$}, so that $\Theta_{*,t}\succ0$ is such that $\Theta_{*,t}\succeq\Theta$ for all $\Theta\in \U_t$. Further, we specify reals $\delta_t\in[0,2]$ and $\gamma\in (0,1)$  such that
$$
\|\Theta^{1/2}[\Theta_{*,t}]^{-1/2}-I_{\nu_t}\|\leq\delta_t\,\,\forall \Theta\in \U_t,
$$
and set\aic{}{\footnote{Note that parameter $\gamma\in (0,1)$ is introduced to prevent $\Phi_t$ to become infinite. Therefore, the larger $\gamma$ is, the better the computed quadratic detector would be. In practice, $\gamma=0.999$ would fit most applications.}}
$$
\H_t=\{(h,H)\in\bR^{\nu_t}\times \bS^{\nu_t}: -\gamma\Theta_{*,t}^{-1}\preceq H\preceq \gamma\Theta_{*,t}^{-1}\}.
$$
Finally, given $t$, we put
$$
\begin{array}{l}
\Phi_{t}(h,H;\Theta)=-\half\ln\Det(I-\Theta_{*,t}^{1/2}H\Theta_{*,t}^{1/2})+\half \Tr([\Theta-\Theta_{*,t}]H)\\
\quad+{\delta_t(2+\delta_t)\over 2(1-\|\Theta_{*,t}^{1/2}H\Theta_{*,t}^{1/2}\|)}\|\Theta_{*,t}^{1/2}H\Theta_{*,t}^{1/2}\|_F^2\\
\quad+{1\over 2}\phi_{\N}\left(B_t^T\left[\hbox{\small$\left[\begin{array}{c|c}H&h\cr\hline h^T&\end{array}\right]+
\left[H,h\right]^T[\Theta_{*,t}^{-1}-H]^{-1}\left[H,h\right]$}\right]B_t\right):\H_t\times\U_t\to\bR,\\
\end{array}
$$
and given $t$, $k$ and $\rho\in(0,R_k]$, we set
$$
\begin{array}{l}
\Phi_{tk\rho}(h,H;\Theta)=-\half\ln\Det(I-\Theta_{*,t}^{1/2}H\Theta_{*,t}^{1/2})+\half \Tr([\Theta-\Theta_{*,t}]H)\\
\quad+{\delta_t(2+\delta_t)\over 2(1-\|\Theta_{*,t}^{1/2}H\Theta_{*,t}^{1/2}\|)}\|\Theta_{*,t}^{1/2}H\Theta_{*,t}^{1/2}\|_F^2\\
\quad+{1\over 2}\phi_{\Z^\rho_k}\left(B_t^T\left[\hbox{\small$\left[\begin{array}{c|c}H&h\cr\hline h^T&\end{array}\right]+
\left[H,h\right]^T[\Theta_{*,t}^{-1}-H]^{-1}\left[H,h\right]$}\right]B_t\right):\H_t\times\U_t\to\bR,\\
\end{array}
$$
and $\Z^\rho_k=\W^\rho_k\bigcap \X.$
\par
Invoking Proposition \ref{concatenation}, we obtain the following
\begin{corollary}\label{corLift}
Given $t\leq d$, $k\leq K$ and $\rho\in(0,R_k]$, consider the convex-concave saddle point problem
\[
{\cal SV}_{tk}(\rho)=\min_{(h,H)\in \H_t}\max_{\Theta_1,\Theta_2\in \U_t}\half\left[\Phi_t(-h,-H;\Theta_1)+\Phi_{tk\rho}(h,H;\Theta_2)\right].
\]
This saddle point problem is solvable, and a saddle point $(h_*,H_*;\Theta_1^*,\Theta_2^*)$ induces quadratic detector
\[
\begin{array}{rcl}
\phi_{tk\rho}(\omega^t)&=&\half[\omega^t]^TH_*\omega^t+h_*^T\omega^t+a:\bR^{\nu_t}\to\bR,\\
a&=&\half\left[\Phi_t(-h_*,-H_*;\Theta_1^*)-\Phi_{tk\rho}(h_*,H_*;\Theta_2^*)\right],\\
\end{array}
\]
such that,
when applied to observation $y^t=\bar{A}_tx+\xi^t$, see {\rm (\ref{eqOSt})}, we have:
\par{\rm (i)} whenever $x\in X$ is a nuisance input,
\begin{equation}\label{case1}
\bE_{y^t}\left\{{\rm e}^{-\phi_{tk\rho}(y^t)}\right\}\leq \epsilon_{tk\rho}:=\exp\{{\cal SV}_{tk}(\rho)\};
\end{equation}
\par{\rm (ii)} whenever  $x\in X$ is a signal  of shape $k$ and magnitude $\geq\rho$,
\begin{equation}\label{case2}
\bE_{y^t}\left\{{\rm e}^{\phi_{tk\rho}(y^t)}\right\}\leq \epsilon_{tk\rho}.
\end{equation}
\end{corollary}

\subsubsection{Construction and performance characterization} The construction to follow is similar to that from Section \ref{impl:constr}. Given $t\leq d$ and $k\leq K$,
it is easily seen that the function ${\cal SV}_{tk}(\rho)$ possesses the following properties:
\begin{itemize}
\item it is nonpositive on $\Delta_k=(0,R_k]$ and nonincreasing in $\rho$ (indeed, $\Phi_t(0,0;\cdot)\equiv\Phi_{tk\rho}(0;0;\cdot)=0$ and $\Phi_{tk\rho}(\cdot,\cdot)$ decreases as $\rho$ grows since $\Z^\rho_k$ shrinks as $\rho$ grows, implying that $\phi_{\Z^\rho_k}(\cdot)$ decreases as $\rho$ grows);
\item the function tends to 0 as $\rho\to +0$;
\item the function is continuous on $\Delta_k$.
\end{itemize}
Given an upper bound $\epsilon\in(0,1/2)$ on the probability of a false alarm, let us set
$$
\epsilon_t={\epsilon\over d},\,\,\varkappa={\epsilon\over\sqrt{dK}},\,\,\alpha=-\ln(dK)/2.
$$
Given $t$, $k$, we define $\rho_{tk}$ as follows: if ${\cal SV}_{tk}(R_k)>\ln(\varkappa)$, we set $\rho_{tk}=+\infty$, otherwise we use bisection to find
$\rho\in(0,R_k]$ such that
$$
{\cal SV}_{tk}(\rho_{tk})=\ln(\varkappa).
$$
Our change detection procedure is as follows: at a step $t=1,2,...,d$, given the observation $y^t$, we look at all values $k\leq K$ for which $\rho_{tk}<\infty$. If $k$ is such that $\rho_{tk}<\infty$, we check whether $\phi_{tk\rho_{tk}}(y^t)<\alpha$. If it is the case, we terminate with a signal conclusion. If $\phi_{tk\rho_{tk}}(y^t)\geq \alpha$ for all $k$ corresponding to $\rho_{tk}<\infty$, we claim that so far, the nuisance hypothesis seems to be valid, and pass
to time $t+1$ (if $t<d$) or terminate (if $t=d$).
\begin{proposition}\label{propquad} Let the input $x\in X$ be observed according to {\rm (\ref{eqOSt})}, and let the observation noise $\xi^d$ be Gaussian with zero mean and covariance matrix $\Theta\in \U_d$. Then
\begin{itemize}
\item if $x$ is a nuisance, the probability for the above detection procedure to terminate with a {signal} conclusion is at most $\epsilon$;
\item if $x$ is a {signal} of shape $k$ and magnitude $\geq\rho>0$, and $t\leq d$ is such that $\rho_{tk}\leq\rho$, then the probability for the detection procedure to terminate with a signal
conclusion at time $t$ or earlier is at least $1-\epsilon$.
\end{itemize}
\end{proposition}

\subsubsection{Numerical illustration}
Here we report on a preliminary numerical experiment with the proposed detection procedure via quadratic detectors.
\paragraph{Observation scheme} we deal with is given by
\begin{equation}\label{ldsm}
\begin{array}{rcl}
z_t&=&Az_{t-1}+Bx_t,\\
w_t&=&Cz_t+\xi_t;\\
\end{array}
\end{equation}
here $z_t$, $x_t$, $w_t$ are, respectively, the states, the inputs and the outputs of a linear dynamical system, of dimensions $n_z$, $n_x$, $n_w$, respectively, and $\xi_t$ are independent across $t$ standard  Gaussian noises.
We assume that the observation at time $t$, $1\leq t\leq d$, is the collection  $w^t=[w_1;w_2;...;w_t]$. In order to account for the initial state $x_0$ and to make the expectations of observations known linear functions of the inputs, we, same as in Section \ref{numill}, define $E_t$ as the linear subspace in $\bR^{n_wt}$ comprised by all collections of accumulated outputs $[w_1;...;w_t]$ of the zero-input system
$$
z_s=Az_{s-1},\,w_s=Cz_s,
$$
and define our (accumulated) observation $y^t$ at time $t$ as the projection of the  observation $w^t$ onto the orthogonal complement $E_t^\perp$ of $E_t$. We represent this projection by the vector $y^t$ of its coordinates in an orthonormal basis of $E_t^\perp$ and set $\nu_t=\dim E_t^\perp$. Note that in this case the corresponding noises $\xi^t$, $1\leq t\leq d$, see (\ref{eqOSt}), are standard Gaussian of dimensions $\nu_t$ (as projections of standard Gaussian vectors), so
that we are in the situation of $\U_t=\{I_{\nu_t}\}$, see (\ref{thetat}). Therefore we can set
$\Theta_{*,t}=I_{\nu_t}$, and $\delta_t=0$, see Section \ref{ssprelim}. \par
We define the admissible nuisance and signals inputs as follows:
\begin{itemize}
\item the admissible inputs $x=[x_1;...;x_d]$, $x_t\in\bR^{n_x}$,  are those with $\|x\|_2\leq R$ (we set $R=10^4$);
\item the only nuisance input is $x=0\in\bR^n$, $n=n_xd$;
\item there are $K=d$ signal shapes, signal of shape $k$ and magnitude $\geq 1$ being a vector of the form $x=[0;...;0;x_k;x_{k+1};...;x_d]$ with $\|x_k\|_2\geq 1$ (``signal of shape $k$ and magnitude $\geq1$ starts at time $k$  with block $x_k$ of energy $\geq1$''). We consider three different types of the signal behavior after time $k$:
    \begin{itemize}
    \item {\sl pulse}: $x_{k+1}=...=x_d=0$,
    \item {\sl step}: $x_k=x_{k+1}=...=x_d$,
    \item  {\sl free jump}:  $x_{k+1},...,x_d$ may be arbitrary.
\end{itemize}
\end{itemize}
The description of the matrix $\bar{A}_t$ arising in (\ref{eqOSt}) is self-evident. The description, required in Section \ref{newsetup}, of
the nuisance set $N$ by quadratic constraints imposed on the quadratic lifting of an input is equally self-evident. The corresponding descriptions of signals of shape $k$ and magnitude $\geq1$ are as follows:
\begin{itemize}
\item {\sl pulse:} $Q_{1k}$ is the diagonal $(n+1)\times (n+1)$ matrix with the only nonzero diagonal entries,  {equal to -1}, in positions $(i,i)$,
$i\in J_k:=\{i:(k-1)n_x+1\leq i\leq kn_x\}$\footnote{VG: $I_k$ was an integer and is now defined as a set! }, and $b_{1k}=-1$. The constraint
$\Tr(Q_{1k}Z(x))\leq b_{1k}$ says exactly that $\|x_k\|_2^2\geq1$. The remaining constraints are homogeneous and express the facts that
    \begin{itemize}
    \item the entries in $Z(x)$ with indices $(i,n+1)$ and $i\leq n$, except for those with $i\in J_k$, are zeros, which can be easily expressed by homogeneous constraints of type A, and
    \item the entries in $Z(x)$ with indices $(i,j)$, $i\leq j\leq n$, except for those with $i,j\in J_k$, are zeros, which can be easily expressed by homogeneous constraints of type B;
    \end{itemize}
\item {\sl step:} $Q_{1k}$ and $b_{1k}$ are exactly as above. The remaining constraints are homogeneous and express the facts that
\begin{itemize}
\item the entries in $Z(x)$ with indices $(i,n+1)$ and  {$i\leq i_k:=n_x(k-1)$} are zero (homogeneous constraints of type A);
\item the entries in $Z(x)$ with indices $i\leq j\leq n$,  {$i \leq i_k$}, are zero (homogeneous constraints of type B);
\item the entries in $Z(x)$ with indices $(i,n+1)$ and $(i',n+1)$  such that  {$i_k < i,i'\leq n$} and $i-i'$ is an integer multiple of $n_x$, are equal to each other (homogeneous constraints of type A);
\item the entries with indices $(i,j)$ and $(i',j')$ such that  {$i_k< i,i',j,j'\leq n$} and both $i-i'$, $j-j'$ are integer multiples of $n_x$, are equal to each other (homogeneous constraints of type B);
    \end{itemize}
\item {\sl free jump:} $Q_{1k}$ and $b_{1k}$ are exactly as above, the remaining constraints are homogeneous and express the facts that
\begin{itemize}
\item the entries in $Z(x)$ with indices $(i,n+1)$,   {$i \leq i_k$}, are zeros (homogeneous constraints of type A);
\item the entries in $Z(x)$ with indices $(i,j)$ such that   {$i \leq i_k$} and $i\leq j$ are zeros (homogeneous constraints of type B).
\end{itemize}
\end{itemize}
\paragraph{Numerical results.} The discrete time dynamical system (\ref{ldsm}) we consider  is obtained by the discretization of the continuous-time model
$$
{d\over ds}\left[\begin{array}{c}u(s)\cr v(s)\cr\end{array}\right]=\underbrace{\left[\begin{array}{c|c}&I_2\cr\hline
&\cr\end{array}\right]}_{A_c}\left[\begin{array}{c}u(s)\cr v(s)\cr\end{array}\right]+\underbrace{\left[\begin{array}{c}\cr I_2\cr\end{array}\right]}_{B_c}x(s),
$$
 with unit time step, assuming the input $x(s)=x_t$ constant on consecutive segments $[t-1,t]$.  We obtain the discrete-time system
$$
\underbrace{\left[\begin{array}{c}u_t\cr v_t\cr\end{array}\right]}_{z_t}=Az_{t-1}+Bx_t,\,\,A=\exp\{A_c\},\,B=\int_0^1\exp\{(1-s)A_c\}B_cds,
$$
or, which is the same, the system
$$
\begin{array}{rcrrl}
u_t&=&u_{t-1}&+v_{t-1}&+\half x_t,\\
v_t&=&&v_{t-1}&+x_t.\\
\end{array}
$$
The system output $u_t$ is observed with the standard Gaussian noise at times $t=1,2,...,d$. Our time horizon was $d=8$, and required probability of false alarm was $\epsilon=0.01$ \par
The results of experiments are presented in Table \ref{table16}; the cells $t,k$ with $k>t$ are blank, because signals of shape $k>t$ start after time $t$ and are therefore ``completely invisible'' at this time. Along with the quantity $\rho_{tk}$ -- the magnitude of the signal of shape $k$ which makes it detectable, with probability $1-\epsilon=0.99$ at time $t$ (the first number in a cell) we present the ``non-optimality index'' (second number in a cell) defined as follows. Given $t$ and $k$, we
compute the largest  {$\rho=\rho_{t k}^*$} such that for a signal  $x^{tk}$ of shape $k$ and magnitude $\geq\rho$, the $\|\cdot\|_2$-norm of $\bar{A}_tx$, see (\ref{eqOSt}), is $\leq 2\ErfInv(\epsilon)$. The latter implies that if all we need to decide at time $t$  is whether the input is the  signal $\theta x^{tk}$ with $\theta<1$, or is identically zero, a $(1-\epsilon)$-reliable decision would be impossible.\footnote{According to our convention, meaningful inputs should be of Euclidean norm at most $10^4$. Consequently, in the case $\rho_{tk}^*>10^4$, we put $\rho_{tk}^*=\infty$.} Since $\theta$ can be made arbitrarily close to 1, $\rho_{tk}^*$ is a lower bound on the magnitude of a signal of shape $k$ which can be detected $(1-\epsilon)$-reliably, by a procedure utilizing observation $y^t$ (cf. Section \ref{sect:assessing}). The non-optimality index reported in the table is the ratio $\rho_{tk}/\rho_{tk}^*$. Note that the computed values of this ratio are neither close to one (which is a bad news for us), nor ``disastrously large'' (which is a good news). In this respect it should be mentioned  that $\rho_{tk}^*$
are overly optimistic estimates of the performance of an ``ideal'' change detection routine.

\begin{center}
\begin{table}
\centering
$$
\begin{array}{|c|}
\hline\hline
\hbox{\small$\begin{array}{||c|c|c|c|c|c|c|c|c||}
\hline
\hbox{\scriptsize$\begin{array}{cc}&k\cr
t&\cr\end{array}$}&1&2&3&4&5&6&7&8\\
\hline
1& \infty/1.00&&&&&&&\\ \hline
2& \infty/1.00& \infty/1.00&&&&&&\\ \hline
3& \infty/1.00&37.8/1.66&37.8/1.66&&&&&\\ \hline
4& \infty/1.00&28.5/1.68&15.6/1.68&28.5/1.67&&&&\\ \hline
5& \infty/1.00&24.8/1.69&11.4/1.69&11.4/1.69&24.8/1.69&&&\\ \hline
6& \infty/1.00&23.0/1.70& 9.6/1.70& 7.9/1.70& 9.6/1.70&23.0/1.70&&\\ \hline
7& \infty/1.00&21.7/1.71& 8.6/1.71& 6.4/1.71& 6.4/1.71& 8.6/1.71&21.7/1.71&\\ \hline
8& \infty/1.00&20.9/1.72& 8.0/1.71& 5.6/1.72& 5.1/1.72& 5.6/1.72& 8.0/1.71&20.9/1.72\\ \hline
\end{array}$}\\
\hbox{Signal geometry: pulse}\\
\hline\hline
\hbox{\small$\begin{array}{||c|c|c|c|c|c|c|c|c||}
\hline
\hbox{\small$\begin{array}{cc}&k\cr
t&\cr\end{array}$}&1&2&3&4&5&6&7&8\\
\hline
1& \infty/1.00&&&&&&&\\ \hline
2& \infty/1.00& \infty/1.00&&&&&&\\ \hline
3&19.0/1.67&19.0/1.67&37.8/1.66&&&&&\\ \hline
4& 7.8/1.68& 7.8/1.68&10.3/1.68&28.5/1.67&&&&\\ \hline
5& 4.2/1.70& 4.2/1.70& 4.9/1.69& 7.9/1.69&24.8/1.69&&&\\ \hline
6& 2.6/1.70& 2.6/1.70& 2.8/1.70& 3.8/1.71& 6.9/1.70&23.0/1.70&&\\ \hline
7& 1.7/1.71& 1.7/1.71& 1.9/1.72& 2.2/1.71& 3.3/1.71& 6.3/1.71&21.7/1.71&\\ \hline
8& 1.2/1.72& 1.2/1.72& 1.3/1.72& 1.5/1.73& 1.9/1.72& 2.9/1.72& 5.9/1.72&20.9/1.72\\ \hline
\end{array}$}\\
\hbox{Signal geometry: step}\\
\hline\hline
\hbox{\small$\begin{array}{||c|c|c|c|c|c|c|c|c||}
\hline
\hbox{\small$\begin{array}{cc}&k\cr
t&\cr\end{array}$}&1&2&3&4&5&6&7&8\\
\hline
1& \infty/1.00&&&&&&&\\ \hline
2& \infty/1.00& \infty/1.00&&&&&&\\ \hline
3& \infty/1.00& \infty/1.00&37.8/1.66&&&&&\\ \hline
4& \infty/1.00& \infty/1.00&38.3/1.68&28.5/1.67&&&&\\ \hline
5& \infty/1.00& \infty/1.00&38.5/1.69&28.7/1.69&24.8/1.69&&&\\ \hline
6& \infty/1.00& \infty/1.00&38.8/1.70&28.9/1.70&25.0/1.70&23.0/1.70&&\\ \hline
7& \infty/1.00& \infty/1.00&39.0/1.71&29.1/1.72&25.3/1.72&23.2/1.72&21.7/1.71&\\ \hline
8& \infty/1.00& \infty/1.00&39.2/1.72&29.1/1.72&25.3/1.72&23.2/1.72&21.8/1.72&20.9/1.72\\ \hline
\end{array}$}\\
\hbox{Signal geometry: free jump}\\
\hline
\end{array}
$$
\caption{\label{table16} Change detection via quadratic detectors. The first number in a cell is $\rho_{tk}$, the second is the nonoptimality index $\rho_{tk}/\rho_{tk}^*$.}
\end{table}
\end{center}

\subsection{Change detection via quadratic detectors, sub-Gaussian case}\label{sectquaddetSG}
Using Proposition \ref{concatenationSG} in the role of Proposition \ref{concatenation}, the constructions and the results of Section \ref{sectquaddet} can be easily adjusted to the situation when the noise $\xi^d$ in {\rm (\ref{eqOS})} is zero mean sub-Gaussian, $\xi^T\sim\SG(0,\Theta)$, rather than Gaussian. In fact, there  are two options
for  such an adjustment, based on quadratic detectors yielded by saddle point problem (\ref{SPPLiftSG}) and convex minimization problem (\ref{SPPLiftConvSG}), respectively. To save space, we restrict ourselves with the first option; utilizing the second option is completely similar.
\par
The only modification of the contents of Section \ref{sectquaddet} needed to pass from Gaussian to sub-Gaussian observation noise is the redefinition of the functions $\Phi_{t}(h,H;\Theta)$  and   $\Phi_{tk\rho}(h,H;\Theta)$
 introduced in Section \ref{ssprelim}. In our  {present} situation,
 \begin{itemize}
 \item  $\Phi_{t}(h,H;\Theta)$  should be redefined as the function $\Phi^{\delta_t}_{A_t,\N}(h,H;\Theta)$ given by relation
(\ref{phiSG}) as applied to $\delta_t$ in the role of $\delta$, $A_t=[\bar{A}_t,0]$ in the role of $A$, the set $\N$, see (\ref{seq101A}), in the role of $\Z$, and the matrix $\Theta_{*,t}$ in the role of $\Theta_*$.
\item  $\Phi_{tk\rho}(h,H;\Theta)$ should be redefined as the function $\Phi^{\delta_t}_{A_t,\Z^\rho_k}(h,H;\Theta)$ given by (\ref{phiSG}) with $\Z^\rho_k$ in the role of $\Z$ and
the just specified  {$\delta_t,A_t,\Theta_*$}.
\end{itemize}
With this redefinition of $\Phi_{t}(h,H;\Theta)$  and   $\Phi_{tk\rho}(h,H;\Theta)$, Corollary \ref{corLift} and Proposition \ref{propquad} (with the words  ``let the observation noise $\xi^d$ be Gaussian with zero mean and covariance matrix $\Theta\in \U_d$'' replaced with ``let the observation noise $\xi^d$ be sub-Gaussian with zero mean and matrix parameter $\Theta\in \U_d$'') remain intact.

\section{{Rust signal detection}}

\subsection{Situation}\label{ysituation}
In this Section, we present an example motivated by material science applications, in which one aims to detect the onset of a rust signal in a piece of metal from a sequence of noisy images. In general, this setup can be used to detect degradation in systems
of a similar nature.

The rust signal occurs at some time, and its energy grows in the subsequent images.
This can be modeled as follows. At times $t=0,1,...,d$,
we observe vectors
\begin{equation}
\label{yeq1}
y_t=y+x_t+\xi_t\in\bR^\nu,
\end{equation}
where
\begin{itemize}
\item $y$ is a fixed deterministic ``background,''
\item $x_t$ is a deterministic {\sl spot},  which may  {correspond} to a rust signal at time $t$, and
\item $\xi_t$ are independent across all $t$ zero mean Gaussian observation noises with covariance matrices $\Sigma_t$.
\end{itemize}
We assume that $x_0=0$, and our ideal goal is to decide on the nuisance hypothesis $x_t=0$, $1
\leq t\leq d$, versus the alternative that the input (``signal'') $x=[x_1;...;x_d]$ is of some shape and some positive magnitude. We specify the shape and the magnitude below.

\subsubsection{Assumptions on observation noise}\label{yobsnoise}
Assume that the observation noise covariance matrices $\Sigma_t$, for all $t$,  are known to belong to a given convex compact subset $\Xi$ of the interior of the
positive semidefinite cone $\bS^\nu_+$. We allow the following two scenarios:
\begin{itemize}
\item[{\bf C.1}]: $\Sigma_t=\Sigma\in\Xi$ for all $t$;
\item[{\bf C.2}]: $\Sigma_t$ can vary with $t$, but stay all the time in $\Xi$.
\end{itemize}
\subsubsection{Assumptions on spots}\label{sec:spots}
We specify signals  $x=[x_1;...;x_d]$ by {\sl shape} $k\in\{1,...,K\}$, $K = d$, and {\sl magnitude} $\rho>0$. Namely, signal $x=[x_1;...;x_d]$ of shape $k$ and magnitude $\geq\rho>0$ ``starts'' at time $k$, meaning that $x_t=0$ when $t<k$. {After the ``change'' happens, the signal satisfies}
\begin{equation}\label{yeq2}
\|x_t\|_2^2\geq \sum_{s=1}^{t-k}\alpha_{t, s}\|x_{t-s}\|_2^2+\aic{\rho}{\airho} p(t-k+1),\,k\leq t\leq K=d,
\end{equation}
where $\alpha_{t, s}$ are given nonnegative coefficients responsible for dynamics of the energies $\|x_t\|_2^2$, and $p(s)$, $1\leq s\leq d$, are given nonnegative coefficients with $p(1)=1$.
For example,
\begin{enumerate}
\item Setting $p(1)=1$ and  $p(s)=0$ for $s>1$,
\begin{itemize}
\item with $\alpha_{t,s}\equiv 0$, we get an ``occasional spot'' of magnitude $\geq {\rho}$ and shape $k$: $x_t=0$ for $t<k$, the energy of $x_k$ is at least $\aic{\rho}{\airho}$, and there are no restrictions on the energy of $x_t$ for $t>k$;
\item with $\alpha_{t,1}=\lambda_t\geq0$ and $\alpha_{t,s}=0$ when $s>1$, we get $x_t=0$ for $t<k$, $\|x_k\|_2^2\geq\aic{\rho}{\airho}$,
    and $\|x_t\|_2^2\geq\lambda_t \|x_{t-1}\|_2^2$ for $t>k$. In other words, the energy of the signal of shape $k$ increases or decreases in a prescribed way after the instant $k$.
    \end{itemize}

\item Setting $\alpha_{t,s}\equiv0$, we get signals of shape $k$ with $x_t=0$ for $t<k$, and energies satisfying $\|x_t\|_2^2 \geq \aic{\rho}{\airho} p(t-k+1)$ for $t\geq k$.
\end{enumerate}
On top of (\ref{yeq2}), we impose on signal $x$ of shape $k$ and magnitude $\geq\rho$  a system (perhaps, empty) of linear constraints
\begin{equation}\label{yeq3}
C_k x \leq {\rho} c_k
\end{equation}
{\sl with $c_k\leq0$}.
\subsection{Processing the situation: formulation}
Let us treat as our observation at time $t$, $t=1,...,d$,  the vector $y^t$ with blocks  $y_i-y_0$, $1\leq i\leq t$, arriving at the observation scheme
\begin{equation}\label{yeq4}
y^t:=[y_1-y_0;...;y_t-y_0]=\bar{A}_tx+\xi^t=S_t[\bar{A}_dx+\xi^d],
\end{equation}
where
\begin{itemize}
\item $\bar{A}_d$ is the unit matrix of size $n=\nu d$; and $S_t$ is the natural projection of $\bR^n$ onto the space of the first $\nu_t=\nu t$ coordinates;
\item $\xi^d\sim\N(0,\Theta)$, where $\Theta$ is a positive semidefinite $d\times d$ block matrix with $\nu\times \nu$ blocks $\Theta^{t\tau}$, $1\leq t,\tau\leq d$, given by\footnote{Indeed, the $t$-th
block $y_t-y_0$ of $y^d$ satisfies
$y_t-y_0=x_t+\xi_t-\xi_0$  {(recall that $x_0=0$)}. Thus, 
for the $t$-th block $\xi_t-\xi_0$ of $\xi^d$ we have
 $\bE\{(\xi_t-\xi_0)(\xi_t-\xi_0)^T\}=\Sigma_t+\Sigma_0$, while
 $\bE\{(\xi_t-\xi_0)(\xi_\tau-\xi_0)^T\}=\Sigma_0$
when $t\neq \tau$, giving rise to \rf{yeq5}.}
\begin{equation}\label{yeq5}
\Theta^{t\tau}=\left\{\begin{array}{ll}\Sigma_0,&t\neq\tau;\\
\Sigma_0+\Sigma_\tau,&t=\tau.\\
\end{array}\right.
\end{equation}
\end{itemize}
We can easily  translate a priori information on $\Sigma_s$, $0\leq s\leq d$, described in Section \ref{yobsnoise}, into a convex compact subset $\U$ of the interior of $\bS^{\nu d}_+$ such that $\Theta$ always belongs to $\U$.
We now cast the above ``spot detection'' problem into the setup from Section \ref{newsetup} as follows. We set $Z(x)=[x;1][x;1]^T$.
\par1.
We assume that the magnitudes of all entries in a meaningful input are bounded by $R$, for a given $R>0$, and put
\[
X=\{x\in \bR^n:\,\Tr(Z(x)Q_i)\leq R^2,\,1\leq i\leq n\}, \;\;Q_i=\Diag\{e_i\},
\]
where $e_i$ is $i$th canonical basis vector in $\bR^{n+1}$.
 We further set $I=n$ and (cf. (\ref{seq101}))
$$
\X=\{W\in\bS^{n+1}_+: W_{n+1,n+1}=1, \,\Tr(WQ_i)\leq R^2,\,1\leq i\leq n\}.
 $$
\par2.
In our current situation, the nuisance set $N$ is the origin. To represent this set in the form (\ref{seq102}), it suffices to set $I_+=I+1={n}+1$, $q_{{n}+1}=0$, and to take, as $Q_{{n}+1}$, the $(n+1)\times (n+1)$ diagonal matrix with the diagonal entries $1,...,1,0$.
We put (cf. \rf{seq101A})
$$
\N=\{W\in\bS^{n+1}_+:\, W_{n+1,n+1}=1,\, \Tr(WQ_i)\leq q_i,\,1\leq i\leq n+1\}.
$$
\par3. Sets $W_k$ of signals of shape $k$ and magnitude $\geq1$, as described in Section \ref{sec:spots}, are given by quadratic constraints on $x=[x_1;...;x_d]$:
\\
$\bullet$ linear constraints on the traces of diagonal blocks $Z_t(x)=x_tx_t^T$  in $Z(x)=[x_1;...;x_d;1][x_1;...;x_d;1]^T$, $1\leq t\leq d$, namely,
\begin{equation}\label{yeq30}
\begin{array}{c}
\Tr(Z_t(x))\leq 0,\;1\leq t\leq  {k-1};\;\;\;-\Tr(Z_k(x))\leq -p(1)=-1;\\
-\Tr(Z_t(x))+\sum_{s=1}^{t-k}  {\alpha_{t,s}\Tr(Z_{t-s}(x))}\leq -p(t-k+1),\;k<t\leq d.\\
\end{array}
\end{equation}
In the terminology of Section \ref{newsetup}, these are type B constraints;\\
\par\noindent
$\bullet$ in addition, linear constraints
$
C_kx\leq c_k
$ defined in (\ref{yeq3})  map to linear constraints
on the first $n$ entries in the last column of $Z(x)$. All these constraints are of type A (recall that $c_k\leq 0$).\par
Observe that among the right hand sides of the constraints (\ref{yeq30}) there is a $(-1)$, implying that all $W_k$ are at a positive distance from the origin.
\par Finally, we put $W^\rho_k=\rho W_k$ and convert these sets, as described in Section \ref{newsetup}, into sets $\W^\rho_k$ such that $Z(x)\in \W^\rho_k$ whenever $x\in W^\rho_k$.
\par
{Note that with our $\X$, all sets $\W_k^\rho$ with small positive $\rho$ do intersect with $\X$.}
\par\noindent
We have covered the problem posed in Section \ref{ysituation} by the setup of Section \ref{newsetup}, and, consequently, can apply the machinery from Section \ref{sectquaddet} to process the problem.
\subsection{Processing the situation: computation}
A computational issue related to this approach stems from the fact that in our intended application $y$ and $x_t$ are images, implying that $\nu=\dim y=\dim x_t$ can be in the range of tens of thousands.
This would make our approach completely unrealistic computationally, unless we can ``kill'' the huge dimensions of the arising convex programs. We are about to demonstrate that under some meaningful structural assumptions this indeed can be done. These assumptions, in their simplest version, are as follows:
\par1. Matrices $\Sigma_t$, $0\leq t\leq d$, are equal to each other and are of the form $\theta\sigma^2I_\nu$, with known $\sigma>0$ and known range $[\vartheta,1]$ of the factor $\theta$,
with $\vartheta\in(0,1]$.
\par2. The only restrictions on the activation signal, apart from the component-wise boundedness, are energy constraints in \rf{yeq2} (e.g., linear constraints as in \rf{yeq3} are not allowed).
\def\SP{{\cal SP}}
\par
Now, computational problems we should solve in the framework of the approach developed in Section \ref{sectquaddet} reduce to building and solving, for given $t\in\{1,...,d\}$, $k\in\{1,...,t\}$, and $\rho>0$, saddle point problems associated with $t,k,\rho$.  Let us fix $t\in\{1,...,d\}$, $k\in\{1,...,t\}$, and $\rho>0$, and let $\SP(t,k,\rho)$ denote the corresponding saddle point problem. This problem is built as follows.
\par1) We deal with observations
\[
y^t=x^t+\xi^t,\,\,\xi^t\sim\N(0,\Theta)
\]
where
\begin{enumerate}
\item[(a)]
 $y^t$, $x^t$ are block vectors with $t$ blocks, $y_i$ and $x_i$, respectively; dimension of every block is $\nu$;
\item[(b)] $\Theta\in \U_t$, where $\U_t$ is comprised of matrices $\Theta=\Theta_\theta$ with $t\times t$ blocks $\Theta^{ij}_\theta$ of size $\nu\times\nu$ such that
\begin{equation}\label{givingrise}
\Theta_\theta^{ij}=\left\{\begin{array}{ll}\theta\sigma^2I_\nu,&i\neq j\\
2\theta\sigma^2I_\nu,&i=j\\
\end{array}\right.,
\end{equation}
with parameter $\theta$ running through $[\vartheta,1]$ (cf. \rf{yeq5}). In other words, denoting by $J_t$ the $t\times t$ matrix with diagonal entries equal to 2 and off-diagonal entries equal to 1, we have
$$
\U_t=\{ {J_t \otimes \theta\sigma^2 I_\nu}:\,\vartheta\leq\theta\leq 1\},
$$
where $A\otimes B$ is the Kronecker product of matrices $A$, $B$: $A\otimes B$ is block matrix obtained by replacing entries {$A_{ij}$ in $A$} with blocks {$A_{ij}B$}.
\end{enumerate}
\par
It is immediately seen that $\U_t$ has the $\succeq$-largest element, specifically, the matrix
\[
\Theta_{*,t}=\sigma^2{J_t\otimes I_\nu}.
\]
Note that
\[
\Theta_{*,t}^{1/2}=\sigma{J_t^{1/2}\otimes I_\nu} \hbox{\ and\ } \Theta\in\U_t\Rightarrow  \|\Theta^{1/2}\Theta_{*,t}^{-1/2}-I_{\nu t}\|\leq \delta:=1-\sqrt{\vartheta}.
\]
\par2) We specify the set $\Z_{tk\rho}\subset \bS^{\nu t+1}_+$ as follows:
\[\begin{array}{rcl}
\Z_{tk\rho}&=&\left\{Z\in\bS^{\nu t+1}_+: Z_{\nu t+1,\nu t+1}=1,
\Tr\left(Z\Diag\{\D_{tks},0\}\right)\leq \aic{\rho}{\airho} d_{tks},\;1\leq s\leq S_{tk}\right\},\\
\D_{tks}&=&{D_{tks}\otimes I_\nu}\\
\end{array}
\]
with {\sl diagonal} $t\times t$ matrices $D_{tks}$ readily given by the coefficients in \rf{yeq30}.
\par
Now, we are in the situation where functions $\Phi_t$ and $\Phi_{tk\rho}$ from Section \ref{ssprelim} are as follows:
{\small\[
\begin{array}{rcl}
\Phi_t(h,H;\Theta)&=&-\half\ln\Det\left(I_{\nu t}-\Theta_{*,t}^{1/2}H\Theta_{*,t}^{1/2}\right)+\half\Tr\left([\Theta-\Theta_{*,t}]H\right)\\
&&+
{\delta(2+\delta)\over 2(1-\|\Theta_{*,t}^{1/2}H\Theta_{*,t}^{1/2}\|)}\|\Theta_{*,t}^{1/2}H\Theta_{*,t}^{1/2}\|_F^2+\half h^T[\Theta_{*,t}^{-1}-H]^{-1} h\\
\Phi_{tk\rho}(h,H;\Theta)
&=&-\half\ln\Det\left(I_{\nu t}-\Theta_{*,t}^{1/2}H\Theta_{*,t}^{1/2}\right)+\half\Tr\left([\Theta-\Theta_{*,t}]H\right)\\
&&+
{\delta(2+\delta)\over 2(1-\|\Theta_{*,t}^{1/2}H\Theta_{*,t}^{1/2}\|)}\|\Theta_{*,t}^{1/2}H\Theta_{*,t}^{1/2}\|_F^2\\
&&+\half\max\limits_{Z\in \Z_{tk\rho}}\Tr\left(Z\hbox{\scriptsize$\left[\begin{array}{c|c}H+H[\Theta_{*,t}^{-1}-H]^{-1}H&h+H[\Theta_{*,t}^{-1}-H]^{-1}h\cr
\hline
h^T+h^T[\Theta_{*,t}^{-1}-H]^{-1}H&h^T[\Theta_{*,t}^{-1}-H]^{-1}h\cr\end{array}\right]$}\right).\\
\end{array}
\]}
The saddle point problem $\SP(t,k,\rho)$ reads
\begin{equation}\label{yeq444}
\begin{array}{c}
\min\limits_{(h,H)\in\H}\left[\Psi(h,H):=\max\limits_{\Theta_1,\Theta_2\in\U_t}\half\left[\Phi_t(-h,-H;\Theta_1)+\Phi_{tk\rho}(h,H;\Theta_2)\right]\right],\\
\H=\{(h,H):-\gamma \Theta_{*,t}^{-1}\preceq H\preceq \gamma\Theta_{*,t}^{-1}\}.
\end{array}
\end{equation}
Observe that
when $(h,H)\in\H$ and $\Theta\in\U_t$, we clearly have $\Phi_t(h,H;\Theta)=\Phi_t(-h,H;\Theta)$ and $\Phi_{tk\rho}(h,H;\Theta)=\Phi_{tk\rho}(-h,H;\Theta)$, where the concluding relation is due to the fact that whenever $Z\in\Z_{tk\rho}$, we
also have $EZE\in\Z_{tk\rho}$, where $E$ is the diagonal matrix with diagonal $1,1,...,1,-1$. As a result, (\ref{yeq444}) has a saddle point with $h=0$, and
building such a saddle point reduces to solving the problem
\begin{equation}\label{yeq4444}
\min\limits_{H\in\widehat{\H}}\left[\widehat{\Psi}(H):=\max\limits_{\Theta_1,\Theta_2\in\U_t}\half\left[\widehat{\Phi}_t(-H;\Theta_1)+\widehat{\Phi}_{tk\rho}(H;\Theta_2)\right]\right],\\
\end{equation}
where
\[\begin{array}{rcl}
\widehat{\Phi}_t(H;\Theta)&=&-\half\ln\Det\left(I_{\nu t}-\Theta_{*,t}^{1/2}H\Theta_{*,t}^{1/2}\right)\\
&&+\half\Tr\left([\Theta-\Theta_{*,t}]H\right)+
{\delta(2+\delta)\over 2(1-\|\Theta_{*,t}^{1/2}H\Theta_{*,t}^{1/2}\|)}\|\Theta_{*,t}^{1/2}H\Theta_{*,t}^{1/2}\|_F^2,\\
\widehat{\Phi}_{tk\rho}(H;\Theta)&=&-\half\ln\Det\left(I_{\nu t}-\Theta_{*,t}^{1/2}H\Theta_{*,t}^{1/2}\right)+\half\Tr\left([\Theta-\Theta_{*,t}]H\right)\\
&&+
{\delta(2+\delta)\over 2(1-\|\Theta_{*,t}^{1/2}H\Theta_{*,t}^{1/2}\|)}\|\Theta_{*,t}^{1/2}H\Theta_{*,t}^{1/2}\|_F^2\\
&&+\half\max\limits_{Z\in \Z_{tk\rho}}\Tr\left(\NW_{\nu t}(Z)\left[H+H[\Theta_{*,t}^{-1}-H]^{-1}H\right]\right),\\
\widehat{\H}&=&\{H:-\gamma \Theta_{*,t}^{-1}\preceq H\preceq \gamma  {\Theta_{*,t}^{-1}}\},\\
\end{array}\]
and $\NW_\ell(Q)$ is the North-Western $\ell\times\ell$ block of $(\ell+1)\times(\ell+1)$ matrix $Q$.
\par
Note that the saddle point problem in (\ref{yeq4444}) has symmetry; specifically, if $\D=I_t\otimes P$ with matrix $P$ which is obtained from $\nu\times\nu$  permutation matrix by replacing some entries equal to 1 with minus these entries, then
\begin{itemize}
\item $\D^T\Theta \D=\Theta$ for every $\Theta\in\U_t$,
\item $\Diag\{\D,1\}^TZ\Diag\{\D,1\}\in \Z_{tk\rho}$ whenever $Z\in\Z_{tk\rho}$,
\item $\D^TH\D\in\widehat{\H}$ whenever $H\in\widehat{\H}$.
\end{itemize}
Hence, as is immediately seen from (\ref{yeq4444}), it holds $\widehat{\Psi}(\D^TH\D)=\widehat{\Psi}(H)$. As a result, (\ref{yeq4444}) has a saddle point with $H=\D^TH\D$ for all indicated $\D$'s, or, which is the same, with $H={G\otimes I_\nu}$, for some $t\times t$ symmetric matrix $G$. Specifying $G$ reduces to solving saddle point problem of sizes
{\sl not affected by $\nu$}, specifically, the problem
\def\J{{\cal J}}
\begin{equation}\label{yeq44444}
\min\limits_{G\in\widehat{\G}}\left[\widetilde{\Psi}(G):=\max\limits_{\I_1,\I_2\in \J_t}\half\left[\widetilde{\Phi}_t(-G;\I_1)+\widetilde{\Phi}_{tk\rho}(G;\I_2)\right]\right],
\end{equation}
\[
\begin{array}{rcl}
\widetilde{\Phi}_t(G;\I)&=&-{\nu\over 2}\ln\Det\left(I_{t}-J_t^{1/2}GJ_t^{1/2}\right)+{\nu\over 2}\Tr\left([\I-J_t]G\right)\\
&&+
{\delta(2+\delta)\nu\over 2(1-\|J_t^{1/2}G J_t^{1/2}\|)}\|J_t^{1/2}GJ_t^{1/2}\|_F^2,\\
\widetilde{\Phi}_{tk\rho}(G;\I)&=&-{\nu\over 2}\ln\Det\left(I_{t}-J_t^{1/2}GJ_t^{1/2}\right)+{\nu\over 2}\Tr\left([\I-J_t]G\right)\\
&&+
{\delta(2+\delta)\nu\over 2(1-\|J_t^{1/2}GJ_t^{1/2}\|)}\|J_t^{1/2}GJ_t^{1/2}\|_F^2\\
&&+{\nu\over 2}\max\limits_{W\in \W_{tk\rho}}\Tr\left(W\left[G+G[J_t^{-1}-G]^{-1}G\right]\right),
\end{array}
\]
where
\[\begin{array}{rcl}
J_t&=&\sigma^2[I_t+[1;...;1][1;...;1]^T],\\
\J_t&=&\{\theta J_t:\vartheta\leq \theta\leq 1\},\\
\widehat{\G}&=&\{G\in\bS^{t}_+:-\gamma J_t^{-1}\preceq G\preceq \gamma J_t^{-1}\},\\
\W_{tk\rho}&=&\left\{W:W\in\bS^{t}_+,\Tr(WD_{tks})\leq \rho \nu^{-1}d_{tks},\,1\leq s\leq S_{tk}\right\}.\\
\end{array}\]
\begin{remark}\label{yrem1} {\rm Our approach is aimed at processing the situation where the magnitude of a spot is quantified by its energy. When $y_t$ represents an image with $\nu$ pixels, this model makes sense if changes in image are more or less spatially uniform, so that a
``typical spot of the magnitude 1'' means small (eventually, $\ll \sigma$)  change in brightness of a significant fraction of the pixels (i.e., we are in the case of  {\em dense alternatives},  in the terminology of \cite{Ingster10}). We can also easily process the model where ``typical spot of magnitude 1'' means large (of order of 1) changes in brightnesses of just few pixels (in the terminology of \cite{Ingster10}, this is the case of  {\em sparse alternatives}). In the latter situation, we
do not need quadratic lift: we can model the set of ``spots of shape $k$ and magnitude $\geq\rho>0$'' as the union of two convex sets, one where the $k$-th entry in the spot is $\geq \rho$, and the other one -- where this entry is $\leq-\rho$. In this model, all we need are affine detectors.}
\end{remark}
\subsection{Real-data example}

\begin{figure}
\begin{center}
\includegraphics[width = 0.9\textwidth]{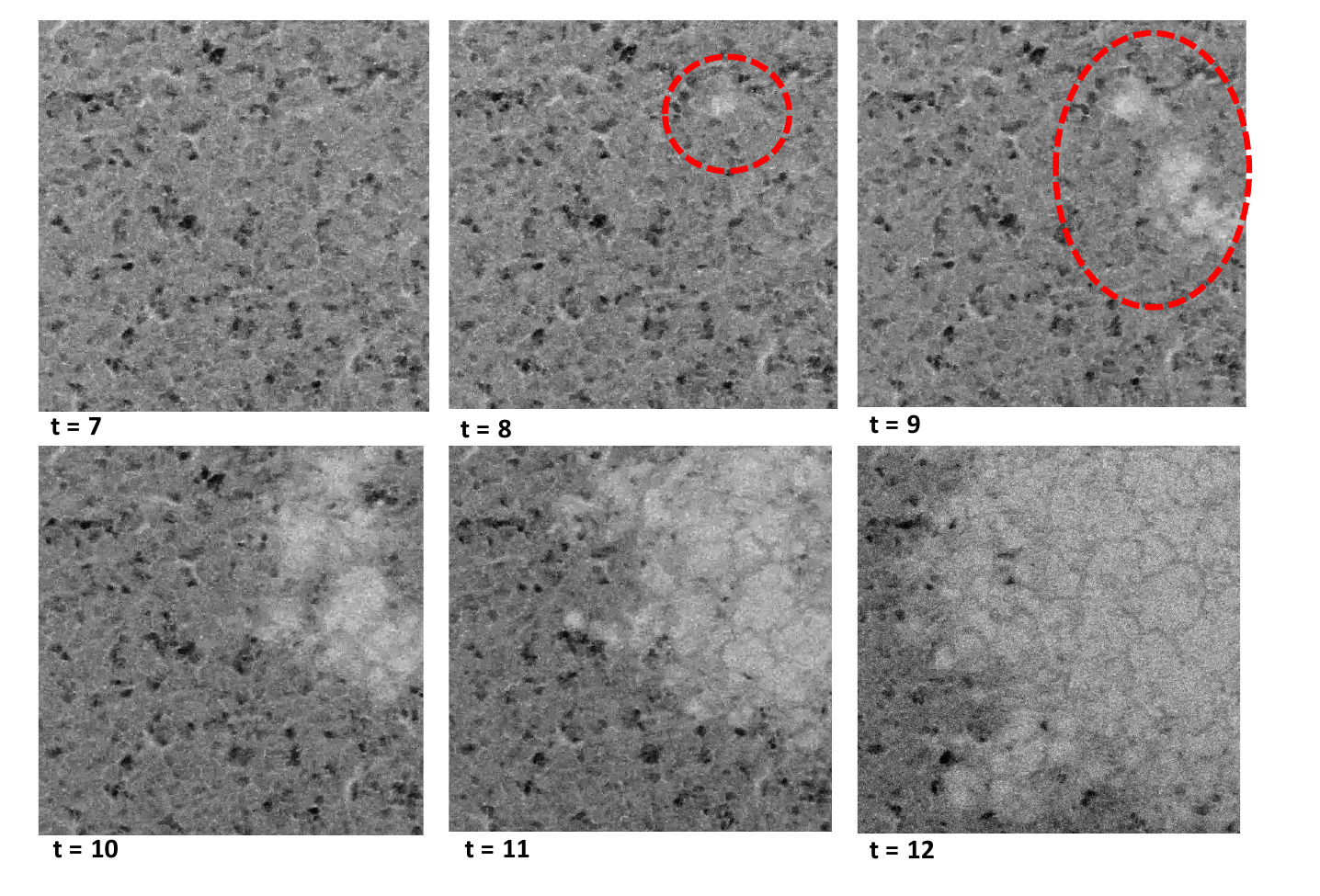}
\end{center}
\caption{A sequence of metal corrosion images. The time (index for the image in the sequence) is labeled; the corrosion initiates at time $t = 8$ (marked by red circle) and develops over time.}
\label{fig:data}
\end{figure}

In this Section, we consider a sequence of metal corrosion images captured using bright-field transmission electron microscopy.\footnote{Data courtesy of
Dr. Josh Kacher at the School of Materials Science and Engineering, Georgia Institute of Technology. More details can be found in Section 3.1 of \cite{CaoZhu17}.
} We downsize each image to 308-by-308 pixels. There are 23 gray images (frames) in the sequence and 2 frames per second. Hence, this corresponds to 11.5 seconds from the original video. At some point, a corrosion spot initiates in the image sequence. Sample images from the sequence are
illustrated in Fig. \ref{fig:data}.

The dynamics of the signal model, in terms of the definition in \rf{yeq2}, has the following parameters: $\alpha_{t,1} = 1$, and $\alpha_{t,s}=0$ for $s>1$; $p(1)=1$ and $p(s)=0$ for $s>1$;  $\rho$ is about $1.2\times 10^2$ and it is estimated from the real-data.

In the example, we set the risk tolerance $\epsilon = 0.1$,  and let $\vartheta = 0.5$. To evaluate detection performance, we run 3000 Monte Carlo trials and add zero-mean Gaussian noise (with variance 25) to the images.
To estimate the noise variance $\sigma^2$, we
use the empirical estimation obtained
taking the first 5 noisy images in the sequence (hence we  assume they do not contain a rust spot). The corresponding estimation is 25.

Since the rust signal is local, i.e., when it occurs, a cluster of pixels captures the rust, we will apply our detector in the following scheme. Break each image into (rectangular or square) patches of equal size. Design a quadratic detector as described above for a patch. Then at each time, whenever one patch detects a change, we claim there has been a change - this corresponds to a ``multi-sensor'' scheme and the local detection statistic by taking
their maximum.\par
We compare our quadratic detector  to the ``sliding window'' ({\tt{Sl-W}}) detector developed in \cite{korostelev2008,guiguesjnps12} and defined as follows. Given ``window width'' $h\in\{1,2,...\}$ and
denoting by
$y_{t j}$ the vector of observations at time $t$ in patch $j$, we build the left and the right estimates, ${\bar y}_{\ell j}^t ( h)$ and ${\bar y}_{r j}^t ( h)$, of $y_{tj}$:
$$
{\bar y}_{\ell j}^t ( h) = \frac{1}{h} \sum_{i=t-h+1}^t y_{i j}  \mbox{ and }{\bar y}_{r j}^t ( h) =\frac{1}{h} \sum_{i=t}^{t+h-1} y_{i j}.
$$
At time $t$, {\tt{Sl-W}} always accepts the nuisance hypothesis when $t \leq 2h-2$; when  $t \geq 2h-1$, the nuisance hypothesis is accepted
if for every patch $j=1,\ldots,N$, it holds
$$
\max_{h \leq \tau \leq t-h+1}\;\|{\bar y}_{\ell j}^\tau ( h) - {\bar y}_{r j}^\tau ( h)\|_{\infty} \leq \kappa,
$$
and is rejected otherwise. In our experiments, $h=2$ and $h=3$ were used. The corresponding thresholds $\kappa$ are computed using Monte-Carle simulation,  see \cite{guiguesjnps12} for details.
\par
Simulation results are presented in Table \ref{singletable}. While the performance of {\tt{Sl-W}} with properly selected $h$ and the number of patches $N$ is quite good, the quadratic detector is a clear winner in terms of reliability (zero empirical probabilities of a false  alarm and a miss), and with $N=49$, there is no delay in detecting the change.
\begin{table}
\centering
\begin{tabular}{|c|}
\hline
\\
\begin{tabular}{|c||c|c|c||}
\hline
 & \multicolumn{3}{c||}{Number $K$ of patches}\\
\hline
Detector& $K=1$  &  $K=4$  &   $K=8$  \\
\hline
{{\tt{Sl-W detector}}}, $h=2$ &[10.0,11.0,14.0] & [10.0,10.6,14] &[10.0,10.2,14]\\
\hline
{{\tt{Sl-W detector}}}, $h=3$ &[10.0,10.9,11.0] &[10.0,10.7,11.0]&[5.0,10.3,11.0]\\
\hline
{{\tt{Quadratic detector}}} & [13.0,13.0,13.0]& [11.0,11.0,11.0] & [10.0,10.0,10.0] \\
\hline
Detector& $K=16$  &  $K=28$  &   $K=49$  \\
\hline
{{\tt{Sl-W detector}}}, $h=2$ &[10.0,10.1,11.0]&[9.0,10.0,11.0]&[10.0,10.1,11]\\
\hline
{{\tt{Sl-W detector}}}, $h=3$ &[5.0,9.8,11.0]&[5.0,8.7,10.0]&[5.0,6.8,10.0]\\
\hline
{{\tt{Quadratic detector}}} &  [10.0,10.0,10.0]  &   [10.0,10.0,10.0] & [8.0,8.0,8.0]    \\
\hline
\end{tabular}\\
Stopping time. Data in a cell $[t_{\min}, \bar{t}, t_{\max}]$:  $\bar{t}$ is the mean, and $[t_{\min}, t_{\max}]$ is the range of \\ instant
where the signal conclusion has been made. The actual change occurs at time $8$.\\
\hline
\\
\begin{tabular}{|c||c|c|c|c|c|c||c|c|c|c|c|c||}
\hline
&  \multicolumn{6}{c||}{ {\tt{False alarm probability}}} &   \multicolumn{6}{c||}{ {\tt{Miss detection rate}}}  \\
\hline
Nb. of patches&  $1$  &  $4$  &   $8$  & $16$  &  $28$  &   $49$  &  $1$  &  $4$  &   $8$  & $16$  &  $28$  &   $49$  \\
\hline
{{\tt{Sl-W detector}}}, $h=2$ & 0 & 0 &  0& 0& 0& 0    &0.34 &0.01&0&0&0&0   \\
\hline
{{\tt{Sl-W detector}}}, $h=3$ & 0 & 0 &  0.006& 0.05& 0.26& 0.64    & 0&0&0&0&0&0   \\
\hline
{{\tt{Quadratic detector}}} & 0&0&0&0&0&0  &0 &0&0&0&0&0 \\
\hline
\end{tabular}\\
Probabilities of false alarm and miss rates.\\
\hline
\end{tabular}
\caption{\label{singletable} Numerical results for rust detection
}
\end{table}

\if{
Fig. \ref{fig:EDD}  illustrates the \aic{expected}{} detection delay (EDD) versus the number of patches. In all instances, the false-detection-rate and miss-detection-rate are both zero.
Note that when the number of patches is greater than 25, we can detect the change immediately after it occurs (the delay is only 1 sample).
}\fi

\if{
\begin{figure}[h!]
\begin{center}
\includegraphics[width = 0.5\textwidth]{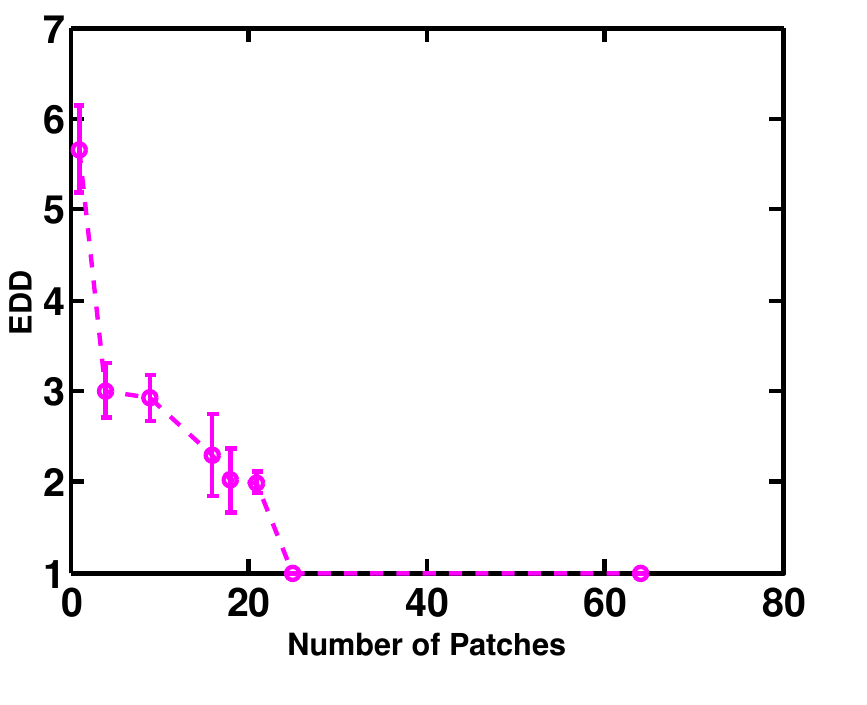}
\end{center}
\caption{The expected detection delay (EDD) versus the number of patches for the metal corrosion data. The miss-detection-rate and false-detection-rate for all instances are zero in this example.}
\label{fig:EDD}
\end{figure}
}\fi

\appendix
\section{Proofs}
\subsection{Proof of Lemma \ref{lem1}}\label{AppLem1}
Recalling what $N^t$ and $U^t_{k\rho}$ are, we have
$$
{\cal SV}_{tk}(\rho)=\min\limits_{h\in\bR^{\nu_t}}\max\limits_{\begin{array}{c}z,v,w,\Theta:\\
z\in N,v\in V_k,w\in W_k,\\
v+\rho w\in X,
\Theta\in \U_t\\
\end{array}}\left[\half h^T\bar{A}_t[v+\rho w-z]+\half h^T\Theta h\right].
$$
From compactness of $X$ and $N$ it follows that the domain of the right hand side saddle point problem is nonempty for all $\rho\in\Delta_k$, and from the fact that $\U_t$ is a compact set contained in the interior of the positive semidefinite cone  it follows that the saddle point of the right hand side exists for all $\rho\in\Delta_k$. We also clearly have
\begin{equation}\label{eq456new}
\begin{array}{rcl}
{\cal SV}_{tk}(\rho)&=&\max\limits_{\begin{array}{c}z,v,w,\Theta:\\
z\in N,v\in V_k,w\in W_k,\\
v+\rho w\in X,
\Theta\in \U_t\\
\end{array}}\underline{\Psi}_{t}(z,v+\rho w,\Theta),\\
\underline{\Psi}_t(p,q,\Theta):&=&-{1\over 8}[p-q]^T[\bar{A}_t^T\Theta^{-1}\bar{A}_t][p-q],\\
\end{array}
\end{equation}
which combines with compactness of $N$, $X$, and $\U_t$ and with the inclusion $\U_t\subset \inter \bS^{\nu_t}_+$ to imply that ${\cal SV}_{tk}(\rho)$ is
nonpositive and continuous on $\Delta_k$. From the same representation, due to semi-conicity of $W_k$, it follows that ${\cal SV}_{tk}(\rho)$ is non-increasing in
$\rho\in\Delta_k$, and that $\liminf_{\rho\to+0}{\cal SV}_{tk}(\rho)\geq0$ due to $0\in N$ and $0\in V_k$ (this was assumed in Section \ref{ADSetUp}), which combines with the fact that ${\cal SV}_{tk}(\rho)$ is nonpositive
to imply that $\lim_{\rho\to +0}{\cal SV}_{tk}(\rho)=0$.
It remains to prove that ${\cal SV}_{tk}(\rho)$ is concave. To this end note that when $\rho\in\Delta_k$, the maximum in (\ref{eq456new}) is achieved, and that $\underline{\Psi}_t(\cdot)$ is concave by the Schur Complement Lemma.
Now let $\rho',\rho''\in \Delta_k$, $\alpha\in[0,1]$, and $\beta=1-\alpha$. We can find $z',z''\in N$, $v',v''\in V_k$, $w',w''\in W_k$ and $\Theta',\Theta''\in\U_t$ such that
$$
\begin{array}{l}
 \underline{\Psi}_t(z',v'+\rho'w',\Theta')={\cal SV}_{tk}(\rho'),
 \underline{\Psi}_t(z'',v''+\rho''w'',\Theta'')={\cal SV}_{tk}(\rho''),\\
  v'+\rho' w'\in X, \,v''+\rho'' w''\in X.\\
  \end{array}
 $$
Setting
$$
\begin{array}{l}
\rho=\alpha\rho'+\beta\rho'', z=\alpha z'+\beta z'',v=\alpha v'+\beta v'',\Theta=\alpha\Theta'+\beta\Theta'',\\
w=\rho^{-1}[\alpha\rho' w'+\beta\rho''w'']=\frac{\alpha \rho'}{\alpha \rho' + \beta \rho''} w' + \frac{\beta \rho''}{\alpha \rho' + \beta \rho''} w'',\\
\end{array}
$$
we get by convexity of $N,X,V_k,W_k,\U_t$:
$$
z\in N,v\in V_k, w\in W_k,\Theta\in \U_t,v+\rho w=\alpha[v'+\rho'w']+\beta[v''+\rho''w'']\in X
$$
and
$$
[z;v+\rho w;\Theta]=\alpha [z';v'+\rho' w';\Theta']+\beta[z'';v''+\rho'' w'';\Theta''].
$$
The latter equality combines with concavity of $\underline{\Psi}_t$ to imply that
$$
\begin{array}{l}
{\cal SV}_{tk}(\alpha\rho'+\beta\rho'')\geq  \underline{\Psi}_t(z,v+\rho w,\Theta)\\
\geq
\alpha \underline{\Psi}_t(z',v'+\rho' w',\Theta')+\beta\underline{\Psi}_t(z'',v''+\rho'' w'',\Theta'')
=\alpha{\cal SV}_{tk}(\rho')+\beta{\cal SV}_{tk}(\rho'').
\end{array}
$$
The resulting inequality holds true for all $\rho',\rho''\in\Delta_k$ and all $\alpha=1-\beta\in[0,1]$, so that ${\cal SV}_{tk}(\cdot)$ is concave.

\par
Now let $\U$ contain the $\succeq$-largest element $\overline{\Theta}$, whence the $\nu_t\times \nu_t$ matrix $\overline{\Theta}_t=S_t\overline{\Theta}S_t^T$, see (\ref{thetat}), is the $\succeq$-largest element in $\U_t$. Then by (\ref{eq456new}) we have
$$
\begin{array}{rcl}
\Gamma_{tk}(\rho)&:=&\sqrt{-{\cal SV}_{tk}(\rho)}=\min\limits_{\begin{array}{c}
z,v,w:\\
z\in N, v\in V_k,w\in W_k,\\
v+\rho w\in X\end{array}}\|H_t[v+\rho w-z]\|_2,\\
H_t&=&{1\over 2\sqrt{2}}[\bar{A}_t^T[\overline{\Theta}_t]^{-1}\bar{A}_t]^{1/2},\\
\end{array}
$$
and from the part of the lemma we have just proved we know that $\Gamma_{tk}(\rho)$ is a continuous nonnegative and nondecreasing function of $\rho\in\Delta_k$ such that $\lim_{\rho\to+0}\Gamma_{tk}(\rho)=0$.
Given $\rho',\rho''\in\Delta_k$ and taking
into account the compactness of $N$, $V_k$, and $X$, we can find $z',v',w',z'',v'',w''$ such that
$$
\begin{array}{c}
z',z''\in N, v',v''\in V_k,w',w''\in W_k,v'+\rho'w'\in X,v''+\rho''w''\in X,\\
\Gamma_{tk}(\rho')=\|H_t[v'+\rho' w'-z']\|_2,\Gamma_{tk}(\rho'')=\|H_t[v''+\rho'' w''-z'']\|_2.\\
\end{array}
$$
Now, given $\alpha\in[0,1]$ and setting $\beta=1-\alpha$,
$$
[\rho;z;v]=\alpha[
\rho';z';v'] + \beta[\rho'';z'';v''],\,w=\rho^{-1}[\alpha\rho'w'+\beta\rho''w''],
$$
we clearly have $z\in N$, $v\in V_k$, $w\in W_k$, $v+\rho w=\alpha[v'+\rho'w']+\beta[v''+\rho''w'']\in X$ and therefore
$$
\begin{array}{rcl}
\Gamma_{tk}(\rho)&\leq& \|H_t[v+\rho w-z]\|_2= \|H_t[\alpha[v'+\rho'w'-z']+\beta[v''+\rho''w''-z'']]\|_2\\
&\leq&\alpha \Gamma_{tk}(\rho')+\beta
\Gamma_{tk}(\rho''),
\end{array}$$
and convexity of $\Gamma_{tk}(\cdot)$ follows.
\hfill$\Box$
\subsection{Proof of Proposition \ref{prop16}}

{\em (i)} Let the input be a nuisance, and let $t\in\{1,...,d\}$. The distribution $P$ of observation $y^t$ in this case belongs to $\SG[N^t,\U_t]$. Now let $k$ be such that $\rho_{tk}<\infty$. Invoking the first inequality in (\ref{eq602}) with $\rho$ set to $\rho_{tk}$, we see that
$$
\int_{\bR^{\nu_t}} \exp\{-\phi_{tk}(y^t)\}P(dy^t)\leq \epsilon_{tk\rho_{tk}}=\exp\{{\cal SV}_{tk}(\rho_{tk})\}=\varkappa_t
$$
(see (\ref{eq601}) and (\ref{balance})). Consequently, $P$-probability of the event $\E_k=\{\xi^t:\,\phi_{tk}(y^t)<\alpha_t\}$ is at most
$\varkappa_t\exp\{\alpha_t\}={\varkappa_t^2\over\epsilon}$.  The signal conclusion at step $t$ is made only when one of the events $\E_k,\, k\in\K_t(\varkappa_t)$, takes place, and $P$-probability of such outcome is  at most
$K_t(\varkappa_t){\varkappa_t^2\over\epsilon}$. We remark that the latter quantity is $\leq \epsilon_t$ by (\ref{eq700}).
\par
{\em (ii)} Now assume that $t$ and $k$ are such that $\rho_{tk}<\infty$, and that the input belongs to $X^\rho_k$ with $\rho\geq \rho_{tk}$. Since $X^\rho_k$ shrinks when $\rho$ grows, the input in fact belongs to $X^{\rho_{tk}}_k$, and therefore the distribution $P$ of observation $y^t$ belongs to $\SG[U^t_{k\rho_{tk}},\U_t]$. Invoking the second inequality in (\ref{eq602}), we get
$$
\int_{\bR^{\nu_t}} \exp\{\phi_{tk}(y^t)\}P(dy^t)\leq \epsilon_{tk\rho_{tk}}=\exp\{{\cal SV}_{tk}(\rho_{tk})\}=\varkappa_t
$$
(see (\ref{eq601}) and (\ref{balance})). Hence, $P$-probability of the event $\E^c_k=\{\xi^t:\,\phi_{tk}(y^t)\geq \alpha_t\}$  is at most $\exp\{-\alpha_t\}\varkappa_t=\epsilon.$ \hfill$\Box$

\subsection{Proof of Proposition \ref{prop17}}

{\em (i)} Let the input be a nuisance, and let $t\in\{1,...,d\}$. The distribution $P$ of observation $y^t$ in this case belongs to $\G[N^t,\U_t]$. Now let $k$ be such that $\rho_{tk}<\infty$. Invoking $(\#)$ and taking into account that, as it was already explained, $\alpha$ as given by (\ref{eqalphabeta}) satisfies $\alpha\leq\delta_t^2$, we conclude that 
the inequality in (\ref{eq1111}.$a$) holds, that is,
 $$
 \begin{array}{l}
 \Prob_{y^t\sim P}\left\{\phi_{tk}(y^t)<\alpha\right\}\leq \Erf(\delta_t-\alpha/\delta_t)\\
 \quad=\Erf\left(\delta_t-\half[\ErfInv(\epsilon)-\ErfInv(\epsilon_t/L_t(\delta_t))]\right)\hbox{\ [by (\ref{eqalphabeta})]}\\
 \quad\leq \Erf\left(\half[\ErfInv(\epsilon)+\ErfInv(\epsilon_t/L_t(\delta_t))]-\half[\ErfInv(\epsilon)-\ErfInv(\epsilon_t/L_t(\delta_t))]\right)\\
 \qquad\quad\hbox{\ [by (\ref{eq801}) and since $\Erf(\cdot)$ is nonincreasing]}\\
 \quad=\min[\epsilon_t/L_t(\delta_t),1/2].
 \end{array}
 $$
This implies that the probability to come to the signal conclusion at step $t$ (this conclusion is made only when $L_t(\delta_t)>0$ and $\phi_{tk}(y^t)<\alpha$ for some $k\in\L_t(\delta_t)$) is at most $L_t(\delta_t)\cdot(\epsilon_t/L_t(\delta_t))=\epsilon_t$, as claimed.
\par
{\em (ii)} Now assume that $t$ and $k$ are such that $\rho_{tk}<\infty$, and that the input belongs to $X^\rho_k$ with $\rho\geq \rho_{tk}$. Since $X^\rho_k$ shrinks when $\rho$ grows, the input in fact belongs to $X^{\rho_{tk}}_k$, and therefore the distribution $P$ of observation $y^t$ belongs to $\G[U^t_{k\rho_{tk}},\U_t]$.
 Since, as it was already explained, $\beta$ as given by (\ref{eqalphabeta}) satisfies $\beta\leq\delta_t^2$, invoking $(\#)$, we conclude that for our $P$ the inequality in (\ref{eq1111}.$b$) holds, that is,
 $$
 \begin{array}{l}
 \Prob_{y^t\sim P}\left\{\phi_{tk}(y^t)\geq\alpha\right\}=\Prob_{y^t\sim P}\left\{\phi_{tk}(y^t)\geq-\beta\right\}\leq \Erf(\delta_t-\beta/\delta_t)\\
 \quad=\Erf\left(\delta_t+\half[\ErfInv(\epsilon)-\ErfInv(\epsilon_t/L_t(\delta_t))]\right)\hbox{\ [by (\ref{eqalphabeta})]}\\
 \quad\leq \Erf\left(\half[\ErfInv(\epsilon)+\ErfInv(\epsilon_t/L_t(\delta_t))]+\half[\ErfInv(\epsilon)-\ErfInv(\epsilon_t/L_t(\delta_t))]\right)\\
 \qquad\quad\hbox{\ [by (\ref{eq801}) and since $\Erf(\cdot)$ is nonincreasing]}\\
 \quad=\epsilon.
 \end{array}
 $$
In other words, in the situation in question, $P$-probability to terminate at time $t$ with the signal conclusion (which is made when $\phi_{tk'}(y^t)<\alpha$
 for some $k'$ with $\rho_{tk'}<\infty$) is at least $1-\epsilon$. \hfill $\Box$

\subsection{Proof of Proposition \ref{prop40}}

Let us refer to the three situations listed in (\ref{chi}) as to cases I, II and III.
Let $\bar{\rho}=\chi\rho_{tk}^*$, with $\chi$ satisfying (\ref{chi}). It may happen that $\bar{\rho}>R_k$; in this case we have nothing to prove, since there do not exist feasible signals of shape $k$ and magnitude $\geq\bar{\rho}$. Now let
\begin{equation}\label{Rkislarge}
\bar{\rho}\leq R_k.
\end{equation}
The function ${\cal SV}_{tk}(\rho)$ is concave on $(0,R_k]$, tends to 0 as $\rho\to +0$ and is equal to $-\half\ErfInv^2(\epsilon)$ when $\rho=\rho_{tk}^*$. Since $\chi\geq1$, we conclude that ${\cal SV}_{tk}(\bar{\rho})\leq \chi{\cal SV}_{tk}(\rho_{tk}^*)$, implying in case II that
\begin{equation}\label{eq757II}
\begin{array}{rcl}
{\cal SV}_{tk}(\bar{\rho})&<& -\half\ErfInv^2(\epsilon)\left(\half\left[1+{\ErfInv(\epsilon/(Kd))\over\ErfInv(\epsilon)}\right]\right)^2\\
&=&
-\half\left(\half\left[\ErfInv(\epsilon)+\ErfInv(\epsilon/(Kd))\right]\right)^2,\\
\end{array}
\end{equation}
and in case III -- that
\begin{equation}\label{eq757III}
{\cal SV}_{tk}(\bar{\rho})< \ln(\bar{\varkappa}),\,\,\bar{\varkappa}=\epsilon/\sqrt{Kd}.
\end{equation}
In case I,  by Lemma \ref{lem1}, the function $\sqrt{-{\cal SV}_{tk}(\rho)}$ is convex on $(0,R_k]$, and therefore the same argument as above shows that  $\sqrt{-{\cal SV}_{tk}(\bar{\rho})}\geq \chi\sqrt{-{\cal SV}_{tk}(\rho_{tk}^*)}$. That is,
${\cal SV}_{tk}(\bar{\rho})\leq \chi^2{\cal SV}_{tk}(\rho_{tk}^*)$, and we again arrive at (\ref{eq757II}).
\par
Let us now consider the Gaussian case. As we have seen, in this case
\begin{equation}\label{eq758}
{\cal SV}_{tk}(\bar{\rho})<-\half\bar{\delta}^2,\,\,\bar{\delta}=\half\left[\ErfInv(\epsilon_t/K)+\ErfInv(\epsilon)\right].
\end{equation}
Taking into account that $L_t(\bar{\delta})\leq K$ and that $\ErfInv(\cdot)$ is non-increasing, we conclude that
$$
\bar{\delta}\geq\half\left[\ErfInv(\epsilon)+\ErfInv(\epsilon_t/L_t(\bar{\delta}))\right],
$$
which combines with (\ref{deltat}) to imply that $\bar{\delta}\geq\delta_t$. Consequently,
 by (\ref{eq758}), we have ${\cal SV}_{tk}(\bar{\rho}) < -\half\delta_t^2$, and since we are in the case where (\ref{Rkislarge}) holds and ${\cal SV}$ is non-increasing, we have ${\cal SV}_{tk}(R_k) < -\half\delta_t^2$ as well. Hence, $k\in\L_t(\delta_t)$, therefore $\rho_{tk}<\infty$, and, in addition, $\bar{\rho}>\rho_{tk}$ (since for finite $\rho_{tk}$ we have ${\cal SV}_{tk}(\rho_{tk})=-\half\delta_t^2$, while ${\cal SV}_{tk}(\bar{\rho})<-\half\delta_t^2$ and ${\cal SV}$ is nonincreasing). Thus, we are in the case of $\rho_{tk}\leq \bar{\rho}\leq R_k$, and therefore, by {item (ii)} of Proposition \ref{prop17}, for a feasible signal of shape $k$ and magnitude $\geq\bar{\rho}$, the probability for the inference procedure from Section \ref{sect:refinement} to terminate at time $t$ with the {signal} conclusion is at least $1-\epsilon$, as required.
 \par
 Now, assume that we are in the sub-Gaussian case. By (\ref{Rkislarge}) combined with (\ref{eq757III}),
 we have
 \begin{equation}\label{assume1}
{\cal SV}_{tk}(R_k)<\ln(\bar{\varkappa}).
\end{equation}
 We claim that
\begin{equation}\label{claim1}
\bar{\varkappa}\leq\varkappa_t,
\end{equation}
where $\varkappa_t$ is given by (\ref{varkappat}). Indeed, we have $K_t(\bar{\varkappa})\leq K={\epsilon\epsilon_t/\bar{\varkappa}^2}$ (recall that $\epsilon_t=\epsilon/d$), which combines with (\ref{varkappat}) to imply (\ref{claim1}). Besides this, we have
\[
\ln(\bar{\varkappa}) \leq -\half\ErfInv^2(\epsilon).
\]
Indeed, we have $0<\ErfInv(\epsilon)\leq \sqrt{2\ln(1/\epsilon)}$, whence
\[
-\half\ErfInv^2(\epsilon)\geq \ln(\epsilon)\geq\ln(\sqrt{\epsilon^2/(dK)})=\ln(\bar{\varkappa}).
\]
Invoking (\ref{assume1}), we get ${\cal SV}_{tk}(R_k)<\ln(\varkappa_t)$ by (\ref{claim1}), that is, recalling the construction from Section \ref{impl:constr}, $\rho_{tk}$ is well defined and satisfies
\begin{equation}\label{rhotksatisfies}
\rho_{tk}\in(0,R_k)\;\;\mbox{with}\;\; {\cal SV}_{tk}(\rho_{tk})=\ln(\varkappa_t).
\end{equation}
Because ${\cal SV}_{tk}(\rho)$ is nonincreasing, we conclude from (\ref{eq757III}), (\ref{claim1}) and the second relation in
(\ref{rhotksatisfies}) that $\bar{\rho}\geq\rho_{tk}$. Invoking {item (ii)} of Proposition \ref{prop16}, we conclude that if the input is a feasible signal with activation of
shape $k$ and magnitude at least $\bar{\rho}$, the probability of the inference routine from Section \ref{impl:constr} to terminate at time $t$ with the signal conclusion is at least $1-\epsilon$.\hfill $\Box$

\subsection{Proof of Proposition \ref{concatenation}}\label{AppConcProof}
Observe that when $m,\nu\in\bR^n$ and $S\in\bS^n$, $\|S\|<1$, one has
{\small\begin{equation}\label{newformula}
\begin{array}{l}
\ln\Big(\bE_{\xi\sim\N(\nu,I_n)} \left\{\exp \{ m^T \xi + \frac{1}{2} \xi^T S \xi \}\right\} \Big)\\
\quad= - \frac{1}{2}\ln\Det(I_n-S)
 + m^T \nu +\frac{1}{2}\nu^TS\nu+\frac{1}{2}(m+S\nu)^T[I_n-S]^{-1}(m+S\nu). \\
\end{array}
\end{equation}}
Indeed,
{\small $$
\begin{array}{l}
\bE_{\xi\sim\N(\nu,I_n)} \left\{\exp \{ m^T \xi + \frac{1}{2} \xi^T S \xi \}\right\}\\
=
{1\over (2\pi)^{n/2}}\exp\{m^T\nu+{1\over 2}\nu^TS\nu\}\int
\exp\{d^T\eta-{1\over 2}\eta^T[I_n-S]\eta\}d\eta \\\multicolumn{1}{r}{\hbox{\ [$d=m+S\nu$,\,$\eta=\xi-\nu$]}}\\
={1\over (2\pi)^{n/2}}\exp\{m^T\nu+{1\over 2}\nu^TS\nu\}\int \exp\{-{1\over 2}\zeta^T[I_n-S]\zeta+{1\over 2}d^T[I_n-S]^{-1}d\}d\zeta\\
\multicolumn{1}{r}{\hbox{\ [$\zeta=\eta-[I_n-S]^{-1}d$]}}\\
=\exp\{m^T\nu+{1\over 2}\nu^TS\nu+{1\over 2}d^T[I_n-S]^{-1}d\}\Det^{-1/2}(I_n-S)\left[{1\over (2\pi)^{n/2}}\int\exp\{-{1\over 2}\omega^T\omega\}d\omega\right]\\
\multicolumn{1}{r}{\hbox{\ [$\omega=[I_n-S]^{1/2}\zeta$]}}\\
=\exp\{m^T\nu+{1\over 2}\nu^TS\nu+{1\over 2}d^T[I_n-S]^{-1}d\}\Det^{-1/2}(I_n-S),\\
\end{array}
$$}
which is exactly (\ref{newformula}).

\subsubsection{Proof of Proposition \ref{concatenation}.i}
\paragraph{1$^0$.} Let $b=[0;...;0;1]\in \bR^{n+1}$, so that $B=\left[\begin{array}{c}A\cr
b^T\cr\end{array}\right]$, and let $\A(u)=A[u;1]$. For any $u\in \bR^n$, $h\in\bR^\nu,\;\Theta\in\bS^\nu_+$ and $H\in\bS^\nu$ such that
$-I\prec \Theta^{1/2}H\Theta^{1/2}\prec I$,  we have
\begin{equation}\label{56eq2}
\begin{array}{l}
\Psi(h,H;u,\Theta):=\ln\left(\bE_{\zeta\sim\N(\A(u),\Theta)}\left\{\exp\{h^T\zeta+\half \zeta^TH\zeta\}\right\}\right)\\
=\hbox{\small$\ln\left(\bE_{\xi\sim\N(0,I)}\left\{\exp \{h^T[\A(u)+\Theta^{1/2}\xi]+\half [\A(u)+\Theta^{1/2}\xi]^TH[\A(u)+\Theta^{1/2}\xi] \} \right\}  \right)$}\\
=-\half \ln\Det(I-\Theta^{1/2}H\Theta^{1/2})+h^T\A(u)\\
\quad+\half \A(u)^TH\A(u)+\half [H\A(u)+h]^T\Theta^{1/2}[I-\Theta^{1/2}H\Theta^{1/2}]^{-1}\Theta^{1/2}[H\A(u)+h]
\\
\multicolumn{1}{r}{\hbox{[by (\ref{newformula})]}}\\
=-\half \ln\Det(I-\Theta^{1/2}H\Theta^{1/2})
+\half [u;1]^T\left[bh^TA+A^Thb^T+A^THA\right][u;1]\\
\quad+\half [u;1]^T\left[B^T[H,h]^T\Theta^{1/2}[I-\Theta^{1/2}H\Theta^{1/2}]^{-1}\Theta^{1/2}[H,h]B\right][u;1]\\
\end{array}
\end{equation}
(because $h^T\A(u)=[u;1]^Tbh^TA[u;1]=[u;1]^TA^Thb^T[u;1]$ and $H\A(u)+h=[H,h]B[u;1]$).
\par
 Observe that when $(h,H)\in \H^\gamma$, we have $\Theta^{1/2}[I-\Theta^{1/2}H\Theta^{1/2}]^{-1}\Theta^{1/2}=[\Theta^{-1}-H]^{-1}\preceq
[\Theta_*^{-1}-H]^{-1},$ so that (\ref{56eq2}) implies that for all $u\in \bR^{n},\;\Theta\in\U,$ and $(h,H)\in\H^\gamma$,
\begin{equation}\label{56eq22}
\begin{array}{rcl}
\Psi(h,H;u,\Theta)&\leq&-\half \ln\Det(I-\Theta^{1/2}H\Theta^{1/2})\\
\multicolumn{3}{r}{+\half [u;1]^T\underbrace{\left[bh^TA+A^Thb^T+A^THA+B^T[H,h]^T[\Theta_*^{-1}-H]^{-1}[H,h]B\right]}_{Q[H,h]}[u;1]}\\
&=&-\half \ln\Det(I-\Theta^{1/2}H\Theta^{1/2})+\half \Tr(Q[H,h]Z(u))\\
&\leq&-\half \ln\Det(I-\Theta^{1/2}H\Theta^{1/2})+\Gamma_{\Z}( h,H),\\
\Gamma_{\Z}(h,H)&=&\half \phi_{\Z}(Q[H,h])\\
\end{array}
\end{equation}
(we have taken into account that $Z(u)\in\Z$ when $u\in U$ (premise of the proposition) and therefore $\Tr(Q[H,h]Z(u))\leq\phi_{\Z}(Q[H,h])$).

\paragraph{2$^0$.} We need the following:
\begin{lemma}\label{lemlogdet} Let $\Theta_*$ be a $d\times d$ symmetric positive definite matrix, let $\delta\in[0,2]$, and let $\U$ be a closed convex subset of $\bS^d_+$ such that
\begin{equation}\label{suchthatweq1}
\Theta\in\U\Rightarrow\{\Theta\preceq\Theta_*\}\ \&\ \{\|\Theta^{1/2}\Theta_*^{-1/2}-I\|\leq\delta\}
\end{equation}
(cf. {\rm (\ref{56delta})}). Let also $\H^o:=\{H\in\bS^d:-\Theta_*^{-1}\prec H\prec \Theta_*^{-1}\}$. Then
\begin{equation}\label{itholdsweq1}
\begin{array}{l}
\forall (H,\Theta)\in\H^o\times\U:\\
G(H;\Theta):=-\half \ln\Det(I-\Theta^{1/2}H\Theta^{1/2})\\
\quad \leq G^+(H;\Theta):=-\half\ln\Det(I-\Theta_*^{1/2}H\Theta_*^{1/2}) +\half\Tr([\Theta-\Theta_*]H)\\
\multicolumn{1}{r}{+{\delta(2+\delta)\over 2(1-\|\Theta_*^{1/2}H\Theta_*^{1/2}\|)}\|\Theta_*^{1/2}H\Theta_*^{1/2}\|_F^2,}\\
\end{array}
\end{equation}
where $\|\cdot\|$ is the spectral, and $\|\cdot\|_F$ - the Frobenius norm of a matrix.  In addition, $G^+(H,\Theta)$ is
 {a} continuous function on $\H^o\times\U$ which is convex in $H\in  {\H^o}$ and concave (in fact, affine)
in $\Theta\in \U$.
\end{lemma}
{\bf Proof.}  Let us set
$$
d(H)=\|\Theta_*^{1/2}H\Theta_*^{1/2}\|,
$$
so that $d(H) <1$ for $H\in\H^o$. For $H\in\H^o$ and $\Theta\in\U$ fixed  we have
\begin{equation}\label{56eq1234}
\begin{array}{rcl}
\|\Theta^{1/2}H\Theta^{1/2}\|&=&\|[\Theta^{1/2}\Theta_*^{-1/2}][\Theta_*^{1/2}H\Theta_*^{1/2}][\Theta^{1/2}\Theta_*^{-1/2}]^T\|\\
&\leq& \|\Theta^{1/2}\Theta_*^{-1/2}\|^2\|\Theta_*^{1/2}H\Theta_*^{1/2}\|\leq \|\Theta_*^{1/2}H\Theta_*^{1/2}\|=d(H)\\
\end{array}
\end{equation}
(we have used the fact that $0\preceq\Theta\preceq\Theta_*$ implies $\|\Theta^{1/2}\Theta_*^{-1/2}\|\leq 1$).
Noting that $\|AB\|_F\leq\|A\|\|B\|_F$,  {a} computation completely similar to the one in (\ref{56eq1234}) yields
\begin{equation}\label{56eq1235}
\|\Theta^{1/2}H\Theta^{1/2}\|_F\leq \|\Theta_*^{1/2}H\Theta_*^{1/2}\|_F=:D(H).
\end{equation}
Besides this, setting $F(X)=-\ln\Det(X):\inter\bS^d_+\to\bR$ and equipping $\bS^d$ with the Frobenius inner product, we have $\nabla F(X)=-X^{-1}$, so that with $R_0=\Theta_*^{1/2}H\Theta_*^{1/2}$, $R_1=\Theta^{1/2}H\Theta^{1/2}$, and $\Delta=R_1-R_0$, we have for properly selected $\lambda\in(0,1)$ and $R_\lambda=\lambda R_0+(1-\lambda)R_1$:
\bse
F(I-R_1)&=&F(I-R_0-\Delta)=F(I-R_0)+\langle \nabla F(I-R_\lambda),-\Delta\rangle\\
&=&F(I-R_0)+\langle(I-R_\lambda)^{-1},\Delta\rangle\\
&=&F(I-R_0)+\langle I,\Delta\rangle +\langle (I-R_\lambda)^{-1}-I,\Delta\rangle.
\ese
We conclude that
\begin{equation}\label{56eqtru}
F(I-R_1)\leq F(I-R_0)+\Tr(\Delta)+\|I-(I-R_\lambda)^{-1}\|_F\|\Delta\|_F.
\end{equation}
Denoting by $\mu_i$ the eigenvalues of $R_\lambda$  and noting that $\|R_\lambda\|\leq \max[\|R_0\|,\|R_1\|]=d(H)$ (see (\ref{56eq1234})), we have $|\mu_i|\leq d(H)$, and therefore eigenvalues $\nu_i=1-{1\over 1-\mu_i}=-{\mu_i\over 1-\mu_i}$ of
$I-(I-R_\lambda)^{-1}$ satisfy $|\nu_i|\leq |\mu_i|/(1-\mu_i)\leq |\mu_i|/(1-d(H))$, whence
$$
\|I-(I-R_\lambda)^{-1}\|_F =\sqrt{{\sum}_{i=1}^d \nu_i^2 }   \leq  \frac{\sqrt{\sum_{i=1}^d  \mu_i^2 }}{1 - d(H)} =    \frac{\|R_\lambda\|_F}{1-d(H)}.
$$
Noting that $\|R_\lambda\|_F \leq  \max[\|R_0\|_F,\|R_1\|_F]   {=} D(H)$, see (\ref{56eq1235}), we conclude that $\|I-(I-R_\lambda)^{-1}\|_F\leq D(H)/(1-d(H))$, so that (\ref{56eqtru}) yields
\begin{equation}\label{56eqtru1}
F(I-R_1)\leq F(I-R_0)+\Tr(\Delta)+D(H)\|\Delta\|_F/(1-d(H)).
\end{equation}
Further, by (\ref{56delta}) the matrix $D =\Theta^{1/2}\Theta_*^{-1/2}-I$ satisfies $\|D\|\leq\delta$, whence
{\small$$
\Delta=\underbrace{\Theta^{1/2}H\Theta^{1/2}}_{R_1}-\underbrace{\Theta_*^{1/2}H\Theta_*^{1/2}}_{R_0}=
(I+D)R_0(I+D^T)-R_0=DR_0+R_0D^T+DR_0D^T.
$$}\noindent
Consequently,
$$
\begin{array}{l}
\|\Delta\|_F\leq \|DR_0\|_F+\|R_0D^T\|_F+\|DR_0D^T\|_F\leq [2\|D\|+\|D\|^2]\|R_0\|_F\\
\multicolumn{1}{r}{\leq \delta(2+\delta)\|R_0\|_F
=
\delta(2+\delta)D(H).}
\end{array}
$$
This combines with (\ref{56eqtru1}) and the relation
\[
\Tr(\Delta)=\Tr(\Theta^{1/2}H\Theta^{1/2}-\Theta_*^{1/2}H\Theta_*^{1/2})=\Tr([\Theta-\Theta_*]H)
\] to yield
$$
F(I-R_1)\leq F(I-R_0)+\Tr([\Theta-\Theta_*]H)+{\delta(2+\delta)\over1-d(H)}\|\Theta_*^{1/2}H\Theta_*^{1/2}\|_F^2,
$$
and we arrive at (\ref{itholdsweq1}). It remains to prove that $G^+(H;\Theta)$ is convex-concave and continuous on $\H^o\times\U$. The only component of this claim which is not completely evident is convexity of the
function in $H\in\H^o$.
To see that it is the case, note that $\ln\Det(S)$ is concave on the interior of the semidefinite cone, the function $f(u,v)={u^2\over 1-v}$ is convex and nondecreasing in $u,v$ in the convex domain $\Pi=\{(u,v):u\geq0,v<1\}$, and the function ${\|\Theta_*^{1/2}H\Theta_*^{1/2}\|_F^2\over 1-\|\Theta_*^{1/2}H\Theta_*^{1/2}\|}$ is obtained from $f$ by convex substitution of variables $H\mapsto(\|\Theta_*^{1/2}H\Theta_*^{1/2}\|_F,\|\Theta_*^{1/2}H\Theta_*^{1/2}\|)$ mapping $\H^o$ into $\Pi$.
\hfill
\qed.

\paragraph{3$^0$.} Combining (\ref{itholdsweq1}), (\ref{56eq22}), (\ref{phi}) and the origin of $\Psi$, see (\ref{56eq2}), we arrive at
$$
\begin{array}{l}
\forall ((u,\Theta)\in U\times\U,(h,H)\in\H^\gamma=\H):\\
\multicolumn{1}{c}{\ln\left(\bE_{\zeta\sim\N( {A[u;1]},\Theta)}\left\{\exp\{h^T\zeta+\half \zeta^TH\zeta\}\right\}\right)\leq\Phi_{A,\Z}(h,H;\Theta),}\\
\end{array}
$$
as claimed in (\ref{moments}).
\paragraph{4$^0$.} Now let us check that $\Phi_{A,\Z}(h,H;\Theta):\H\times\U\to\bR$ is continuous and convex-concave. Recalling that the function $G^+(H;\Theta)$ from (\ref{itholdsweq1}) is convex-concave and continuous on $\H^o\times\U$, all we need to verify is that $\Gamma_{\Z}(h,H)$ is convex and continuous on $\H$.
Recalling that $\Z$ is a nonempty compact set, the function $\phi_{\Z}(\cdot):\,\bS^{n+1}
\to\bR$ is continuous, implying the continuity of $\Gamma_{\Z}(h,H)=\half\phi_{\Z}(Q[H,h])$ on $\H=\H^\gamma$ ($Q[H,h]$ is defined in (\ref{56eq22})). To prove the convexity of $\Gamma$, note that $\Z$ is
contained in $\bS^{n+1}_+$, implying that $\phi_{\Z}(\cdot)$ is convex and $\succeq$-monotone. On the other hand, by Schur Complement Lemma, we have
$$
\begin{array}{rcl}
S&:=&\{(h,H,G): G\succeq Q[H,h],(h,H)\in\H^\gamma\}\\
&=&\bigg\{(h,H,G):\left[\begin{array}{c|c}G-[bh^TA+  {A^T h b^T} +  A^THA]&B^T[H,h]^T\cr\hline
[H,h]B&\Theta_*^{-1}-H\cr\end{array}\right]\\
&&\multicolumn{1}{r}{\succeq0,(h,H)\in\H^\gamma\bigg\},}\\
\end{array}
$$
implying that $S$ is convex. Since $\phi_{\Z}(\cdot)$ is $\succeq$-monotone, we have
$$
\begin{array}{l}
\{(h,H,\tau):(h,H)\in \H^\gamma,\;\tau\geq  {\Gamma_\Z}(h,H)\}=\{(h,H,\tau):
\exists G: G\succeq Q[H,h],\\
\multicolumn{1}{r}{2\tau\geq \phi_{\Z}(G),\;(h,H)\in\H^\gamma\},}
\end{array}
$$
and we see that the epigraph of $\Gamma_\Z$ is convex (since the set $S$ and the epigraph of $\phi_{\Z}$ are so), as claimed.

\paragraph{5$^0$.} It remains to prove that $\Phi_{\A,\Z}$ is coercive in $H,h$. Let $\Theta\in\U$ and $(h_i,H_i)\in\H^\gamma$ with $\|(h_i,H_i)\|\to\infty$ as $i\to\infty$, and
let us prove that $\Phi_{\A,\Z}(h_i,H_i;\Theta)\to\infty$. Looking at the expression for $\Phi_{\A,\Z}(h_i,H_i;\Theta)$, it is immediately seen that all
terms in this expression, except for the terms coming from $\phi_{\Z}(\cdot)$, remain bounded as $i$ grows, so that all we need to verify  is
that the $\phi_{\Z}(\cdot)$-term goes to $\infty$ as $i\to\infty$. Observe that $H_i$ are uniformly bounded due
to $(h_i,H_i)\in \H^\gamma$, implying that $\|h_i\|_2\to\infty$ as $i\to\infty$.\\
Denoting  $e=[0;...;0;1]\in\bR^{\nu+1}$ and, as before,
$b=[0;...;0;1]\in \an{\bR^{\nu+1}}{\bR^{n+1}}$, note that, by construction, $B^Te=b$. Now let $W\in\Z$, so that
$W_{n+1,n+1}=1$. Taking into account that the matrices $[\Theta_*^{-1}-H_i]^{-1}$ satisfy $\alpha I_d\preceq [\Theta_*^{-1}-H_i]^{-1}\preceq \beta I_d$  for some positive $\alpha,\beta$ due to $H_i\in\H^\gamma$, we come to
$$
\begin{array}{l}
\underbrace{\left[\left[\begin{array}{c|c}H_i&h_i\cr\hline h_i^T&\end{array}\right]+
\left[H_i,h_i\right]^T[\Theta_*^{-1}-H_i]^{-1}\left[H_i,h_i\right]\right]}_{Q_i}
=
\underbrace{\left[h_i^T[\Theta_*^{-1}-H_i]^{-1}h_i\right]}_{\alpha_i\|h_i\|_2^2}ee^T+R_i,\\
\end{array}
$$
where $\alpha_i\geq\alpha>0$ and $\|R_i\|_F\leq C(1+\|h_i\|_2)$. As a result,
$$
\begin{array}{rcl}
\phi_{\Z}(B^TQ_iB)&\geq&\Tr(WB^TQ_iB)=\Tr(WB^T[\alpha_i\|h_i\|_2^2 ee^T+R_i]B)\\
& {\geq}&\alpha_i\|h_i\|_2^2
\underbrace{\Tr(Wbb^T)}_{=W_{n+1,n+1}=1}-\|BWB^T\|_F\|R_i\|_F\\
&\geq&
\alpha\|h_i\|_2^2-C(1+\|h_i\|_2)\|BWB^T\|_F,\\
\end{array}
$$
and the concluding quantity tends to $\infty$ as $i\to\infty$ due to $\|h_i\|_2\to\infty$, $i\to\infty$. Part (i) is proved.

\subsubsection{Proof of Proposition \ref{concatenation}.ii}
Part (ii) of the proposition is a straightforward combination of part (i) and \cite[Proposition 3.1]{PartI}; for the sake of completeness, here is a simple proof. Since by (i) the function $\Phi(h,H;\Theta_1,\Theta_2)$ is continuous and convex-concave on the domain $\underbrace{(\H_1\cap\H_2)}_{\H}\times \underbrace{(\U_1 \times \U_2)}_{\U}$  and are coercive in $(h,H)$, while $\H$ and $\U$ are closed and convex, and $\U$ in addition is compact, saddle point problem (\ref{SPPLift}) is solvable (Sion-Kakutani Theorem). Now let $(h_*,H_*;\Theta_1^*,\Theta_2^*)$ be a saddle point. To prove
(\ref{riskagain}), let $P\in\G_1$, that is, $P=\N(A_1[u;1],\Theta_1)$ for some $\Theta_1\in \U_1$ and some $u$ with $[u;1][u;1]^T\in \Z_1$. Applying (\ref{moments}) to the first collection of data, with $a$ given by (\ref{iquaddet}), we get the first $\leq$ in the following chain:
$$
\begin{array}{l}
\ln\left(\displaystyle \int {\rm e}^{-\half\omega^TH_*\omega-\omega^Th_*-a}P(d\omega)\right)\\
\leq \Phi_{A_1,\Z_1}(-h_*,-H_*;\Theta_1)-a
\underbrace{\leq}_{
(a)}
 \Phi_{A_1,\Z_1}(-h_*,-H_*;\Theta_1^*)-a
 \underbrace{=}_{(b)}{\cal SV},\\
 \end{array}
 $$
 where $(a)$ is due to the fact that $\Phi_{A_1,\Z_1}(-h_*,-H_*;\Theta_1)+\Phi_{A_2,\Z_2}(h_*,H_*;\Theta_2)$ attains its maximum over $(\Theta_1,\Theta_2)\in \U_1\times\U_2$ at the point $(\Theta_1^*,\Theta_2^*)$, and $(b)$ is due to the origin of $a$ and the relation ${\cal SV}=\half[\Phi_{A_1,\Z_1}(-h_*,-H_*;\Theta_1^*)+\Phi_{A_2,\Z_2}(h_*,H_*;\Theta_2^*)]$. The bound in (\ref{riskagain}.a) is proved. Similarly, let $P\in\G_2$, that is, $P=\N(A_2[u;1],\Theta_2)$ for some $\Theta_2\in \U_2$ and some $u$ with $[u;1][u;1]^T\in \Z_2$. Applying (\ref{moments}) to the second collection of data, with the same $a$ as above, we get the first $\leq$ in the following chain:
$$
\begin{array}{l}
\ln\left(\displaystyle \int {\rm e}^{\half\omega^TH_*\omega+\omega^Th_*+a}P(d\omega)\right)\leq \Phi_{A_2,\Z_2}(h_*,H_*;\Theta_2)+a\\
\underbrace{\leq}_{
(a)}
 \Phi_{A_2,\Z_2}(h_*,H_*;\Theta_2^*)+a\underbrace{=}_{(b)}{\cal SV},\\
 \end{array}
 $$
with exactly the same justification as above of $(a)$ and $(b)$.  The bound in (\ref{riskagain}.b) is proved. \hfill \qed

\subsection{Justification for Remark \ref{remark41}.}\label{just_remark41}

In the easy case, we have  {$B_\chi = I_{\nu+1}$ and therefore}
$$
\begin{array}{rcl}
M_\chi(h,H)&:=&B_\chi^T\left[\hbox{\small$\left[\begin{array}{c|c}H&h\cr\hline h^T&\end{array}\right]+
\left[H,h\right]^T\left[[\Theta_*^{(\chi)}]^{-1}-H\right]^{-1}\left[H,h\right]$}\right]B_\chi\\
&=&\hbox{\small$\left[\begin{array}{c|c}H+H\left[[\Theta_*^{(\chi)}]^{-1}-H\right]^{-1}H&h+ {H}[ {[\Theta_*^{(\chi)}]^{-1}}-H]^{-1}h\cr\hline
h^T+h^T\left[[\Theta_*^{(\chi)}]^{-1}-H\right]^{-1} {H}&h^T\left[[\Theta_*^{(\chi)}]^{-1}-H\right]^{-1}h\cr\end{array}\right]$}\\
\end{array}
$$
and
$$
\begin{array}{rcl}
\phi_{\Z_\chi}(Z)&=&\max\limits_{W}\left\{\Tr(ZW):\,W\succeq0,\,\Tr(WQ^\chi_j)\leq q^\chi_j,\,1\le j\leq J_\chi\right\}\\
&=&\min_\lambda\left\{\sum_jq_j^\chi\lambda_j:\,\lambda\geq0,\,Z\preceq \sum_j\lambda_jQ^\chi_j\right\},\\
\end{array}
$$
where the last equality is due to semidefinite duality.
From the second representation of $\phi_{\Z_\chi}(\cdot)$ and the fact that all $Q^\chi_j$ are diagonal it follows that
$ \phi_{\Z_\chi}( {M_\chi}({0,H}))\leq \phi_{\Z_\chi}( {M_\chi}(h,H))$ (indeed, with diagonal $Q^\chi_j$, this representation clearly says that if $\lambda$ is feasible for the minimization problem participating in the representation when $Z=M_\chi(h,H)$, it remains feasible when $Z$
 is replaced with $M_\chi(0,H)$).
This, in turn, combines straightforwardly with (\ref{phi}) to imply that when replacing  $h_*$ with $0$ in a saddle point $(h_*,H_*;\Theta_1^*,\Theta_2^*)$ of (\ref{SPPLift}), we end up with another saddle
point of (\ref{SPPLift}). In other words, when solving (\ref{SPPLift}), we can from the very beginning set $h$ to $0$,  thus converting (\ref{SPPLift}) into the convex-concave saddle point problem
\begin{equation}\label{spsimpl}
{\cal SV}=\min\limits_{H: (0,H)\in\H_1\cap\H_2}\max\limits_{\Theta_1\in\U_1,\Theta_2\in\U_2} \Phi(0,H;\Theta_1,\Theta_2).
\end{equation}
Taking into account the fact that we are in the case where all matrices from the sets $\U_\chi$, same as the matrices $\Theta_*^{(\chi)}$ and all the matrices $Q^\chi_j$, $\chi=1,2$, are diagonal, it is
immediate to verify that if $E$ is a $\nu\times \nu$ diagonal matrix with diagonal entries $\pm 1$, then $\Phi(0,H;\Theta_1,\Theta_2)=\Phi(0,E H E;\Theta_1,\Theta_2)$.
Due to convexity-concavity of $\Phi$ this implies that (\ref{spsimpl}) admits a saddle point $(0,H_*;\Theta_1^*,\Theta_2^*)$
with $H_*$ invariant w.r.t. transformations $H_*\mapsto EH_*E$ with the above $E$, that is, with diagonal $H_*$, as claimed.\hfill \qed
\subsection{Proof of Proposition \ref{concatenationSG}}\label{AppConcProofSG}
\subsubsection{Preliminaries}
We start with the following result:
\begin{lemma} \label{lemmumu} Let $\bar{\Theta}$ be a positive definite $d\times d$ matrix, and let
$$
u\mapsto\A(u)=A[u;1]
$$
be an affine mapping from $\bR^n$ into $\bR^d$. Finally, let $h\in\bR^d$, $H\in\bS^d$ and $P\in\bS^d$
satisfy the relations
\begin{equation}\label{weq11}
0\preceq P\prec I\ \&\ P\succeq \bar{\Theta}^{1/2}H\bar{\Theta}^{1/2}.
\end{equation}
Then, setting {\scriptsize$B=\left[\begin{array}{cc}A\cr0,...,0,1\cr\end{array}\right]$}, for every $u\in\bR^n$ it holds
\begin{equation}\label{weq111}
\begin{array}{c}
\multicolumn{1}{l}{\zeta\sim \SG(\A(u), {\bar{\Theta}})\Rightarrow
\ln\left(\bE_\zeta\left\{{\rm e}^{h^T\zeta+\half \zeta^TH\zeta}\right\}\right)\leq
-\half\ln\Det(I-P)\qquad\qquad\qquad}\\
\multicolumn{1}{r}{+\half[u;1]^TB^T\left[\Hplus+\left[H,h\right]^T
\bar{\Theta}^{1/2}[I-P]^{-1}\bar{\Theta}^{1/2}\left[H,h\right]\right]B[u;1]}.\\
\end{array}
\end{equation}
Equivalently (set $G=\bar{\Theta}^{-1/2}P\bar{\Theta}^{-1/2}$), whenever $h\in\bR^d$, $H\in\bS^d$, and $G\in\bS^d$
satisfy the relations
\begin{equation}\label{weq11a}
0\preceq G\prec \bar{\Theta}^{-1}\ \&\ G\succeq H,
\end{equation}
one has for every for every $u\in\bR^n$:
\begin{equation}\label{weq111a}
\begin{array}{c}
\multicolumn{1}{l}{\zeta\sim \SG(\A(u),\bar{\Theta})\Rightarrow
\ln\left(\bE_\zeta\left\{{\rm e}^{h^T\zeta+\half \zeta^TH\zeta}\right\}\right)\leq
-\half\ln\Det(I-\bar{\Theta}^{1/2}G\bar{\Theta}^{1/2})}\\
\multicolumn{1}{r}{+\half[u;1]^TB^T\left[\Hplus+\left[H,h\right]^T
[\bar{\Theta}^{-1}-G]^{-1}\left[H,h\right]\right]B[u;1]}.
\end{array}
\end{equation}
\end{lemma}
{\bf Proof.} {\bf 1$^0$.} Let us start with the following observation:\\
\begin{lemma}\label{wlem1} Let $ {\bar{\Theta}}\in\bS^d_+$ and $S\in\bR^{d\times d}$ be such that $S  {\bar \Theta} S^T\prec I_d$. Then for every $\nu\in\bR^d$ one has
\begin{equation}\label{steq1}
\begin{array}{l}
{\ln\bigg(}\bE_{\xi\sim\SG(0, {\bar{\Theta}})}\left\{{\rm e}^{\nu^TS\xi+\half\xi^TS^TS\xi}\right\}{\bigg)}\leq
{\ln\bigg(}\bE_{x\sim\N(\nu,I_d)}\left\{{\rm e}^{\half x^TS {\bar{\Theta}} S^Tx}\right\}{\bigg)}\\
\quad=-\half\ln\Det(I_d-S {\bar{\Theta}} S^T)+\half \nu^T\left[S {\bar{\Theta}} S^T( {I_d}-S {\bar{\Theta}} S^T)^{-1}\right]\nu.\\
\end{array}
\end{equation}
\end{lemma}
{\bf Proof.} Let $\xi\sim \SG(0, {\bar{\Theta}})$ and $x\sim\N(\nu,I_d)$ be independent. We have:
$$
\begin{array}{l}
\bE_\xi\left\{{\rm e}^{\nu^TS\xi+\half\xi^TS^TS\xi}\right\}\underbrace{=}_{a}\bE_\xi\left\{\bE_x\left\{{\rm e}^{[S\xi]^Tx}\right\}\right\}
=\bE_x\left\{\bE_\xi\left\{{\rm e}^{[S^Tx]^T\xi}\right\}\right\}\\
\underbrace{\leq}_{b}\bE_x\left\{{\rm e}^{  {\frac{1}{2}x^TS {\bar \Theta}  S^Tx}}\right\},\\
\end{array}
$$
where $a$ is due to $x\sim\N(0,I_d)$ and $b$ is due to $\xi\sim\SG(0, {\bar{\Theta}})$. We have verified the inequality in (\ref{steq1}); the equality in (\ref{steq1}) is given by direct computation.~\hfill \qed
\medskip\par\noindent {\bf 2$^0$.} Now, in the situation described in Lemma \ref{lemmumu}, given $u\in\bR^n$,
 let us set $\mu=\A(u)=A[u;1]$, $\nu=P^{-1/2}\bar{\Theta}^{1/2}[H\mu+h]$, $S=P^{1/2}\bar{\Theta}^{-1/2}$, so that
 $S\bar{\Theta}S^T=P\prec  {I_d}$  {and $G={\bar \Theta}^{-1/2} P {\bar \Theta}^{-1/2} = S^T S$}.
 Let  $\zeta\sim\SG(\mu,\bar{\Theta})$. Representing $\zeta$ as $\zeta=\mu+\xi$ with $\xi\sim\SG(0,\bar{\Theta})$, we have
$$
\begin{array}{l}
\ln\left(\bE_\zeta\left\{{\rm e}^{h^T\zeta+\half \zeta^TH\zeta}\right\}\right)=
h^T\mu+{1\over 2}\mu^TH\mu+\ln\left(\bE_\xi\left\{{\rm e}^{[h+H\mu]^T\xi+\half \xi^TH\xi}\right\}\right)\\
\leq h^T\mu+{1\over 2}\mu^TH\mu+\ln\left(\bE_\xi\left\{{\rm e}^{[h+H\mu]^T\xi+\half \xi^T  {G}\xi}\right\}\right)\hbox{\ [since  {$G\succeq H$}]}\\
=
h^T\mu+{1\over 2}\mu^TH\mu+\ln\left(\bE_\xi\left\{{\rm e}^{[h+H\mu]^T\xi+\half \xi^T  {S^T S}\xi}\right\}\right)\hbox{\  {[since $G=S^T S$]}}\\
= h^T\mu+{1\over 2}\mu^TH\mu+\ln\left(\bE_\xi\left\{{\rm e}^{\nu^TS\xi+\half \xi^T  {S^T S}\xi}\right\}\right)\hbox{\ [since $S^T\nu=h+H\mu$]}\\
\leq h^T\mu+\half\mu^TH\mu -\half\ln\Det(I_d-S\bar{\Theta} S^T)+\half \nu^T\left[S\bar{\Theta} S^T(I_d-S\bar{\Theta} S^T)^{-1}\right]\nu\\
\quad\hbox{\ [by Lemma \ref{wlem1}]}\\
=h^T\mu+\half\mu^TH\mu-\half\ln\Det(I_d-P)\\
\quad+\half [H\mu+h]^T\bar{\Theta}^{1/2}( {I_d-P})^{-1}\bar{\Theta}^{1/2}[H\mu+h]
\hbox{\ [plugging in $S$ and $\nu$].}\\
\end{array}
$$
It is immediately seen that the concluding quantity in this chain is nothing but the right hand side quantity in (\ref{weq111}). \hfill \qed
\subsubsection{Proof of Proposition \ref{concatenationSG}.i}
\paragraph{1$^0$.} Let us prove (\ref{momentsSG}.$a$). By Lemma \ref{lemmumu} (see (\ref{weq111a})) applied with $\bar{\Theta}=\Theta_*$, setting $\A(u)=A[u;1]$, we have
\begin{equation}\label{whatweused}
\begin{array}{l}
\forall\left((h,H)\in\H, G:0\preceq G\preceq\gamma^+\Theta_*^{-1},G\succeq H,u\in\bR^n:[u;1][u;1]^T\in\Z\right):\\
\ln\left(\bE_{\zeta\sim\SG(\A(u),\Theta_*)}\left\{{\rm e}^{h^T\zeta+\half\zeta^TH\zeta}\right\}\right)\leq
-\half\ln\Det(I-\Theta_*^{1/2}G\Theta_*^{1/2})\\
\multicolumn{1}{r}{+\half[u;1]^TB^T\left[\Hplus+\left[H,h\right]^T
[\Theta_*^{-1}-G]^{-1}\left[H,h\right]\right]B[u;1]}\\
\leq -\half\ln\Det(I-\Theta_*^{1/2}G\Theta_*^{1/2})\\
\multicolumn{1}{r}{+\half\phi_\Z\left(B^T\left[\Hplus+\left[H,h\right]^T
[\Theta_*^{-1}-G]^{-1}\left[H,h\right]\right]B\right)}\\
=\Psi_{A,\Z}(h,H,G),\\
\end{array}
\end{equation}
implying, due to the origin of $\Phi_{A,\Z}$,  that under the premise of (\ref{whatweused}) we have
$$
\ln\left(\bE_{\zeta\sim\SG(\A(u),\Theta_*)}\left\{{\rm e}^{h^T\zeta+\half\zeta^TH\zeta}\right\}\right)\leq \Phi_{A,\Z}(h,H),\,\forall (h,H)\in\H.
$$
Taking into account that when $\zeta\sim\SG(\A(u),\Theta)$ with $\Theta\in\U$, we have also $\zeta\sim\SG(\A(u),\Theta_*)$ due to $\Theta\preceq \Theta_*$, (\ref{momentsSG}.$a$) follows.

\paragraph{2$^0$.} Now let us prove (\ref{momentsSG}.$b$). All we need is to verify the relation
\begin{equation}\label{allweneed}
\begin{array}{l}
\forall\left((h,H)\in\H,G,0\preceq G\preceq\gamma^+\Theta_*^{-1}, {H \preceq G},u\in  {\bR^n}:[u;1][u;1]^T\in\Z,\Theta\in\U\right):\\
\multicolumn{1}{c}{\qquad\qquad\qquad
\ln\left(\bE_{\zeta\sim\SG(\A(u), {\Theta)}}\left\{{\rm e}^{h^T\zeta+\half \zeta^TH\zeta}\right\}\right)\leq\Psi^\delta_{A,\Z}(h,H,G;\Theta);
}\\
\end{array}
\end{equation}
with this relation at our disposal (\ref{momentsSG}.$b$) can be obtained by the same argument as the one we used in item {\textbf{1$^0$}} to derive (\ref{momentsSG}.$a$).\par To establish (\ref{allweneed}), let us fix $h,H,G,u,\Theta$ satisfying the premise of (\ref{allweneed}); note that under the premise of Proposition \ref{concatenationSG}.i, we have  $0\preceq\Theta\preceq \Theta_*$. Now let $\lambda\in(0,1)$, and let $\Theta_\lambda=\Theta+\lambda(\Theta_*-\Theta)$, so that $0\prec\Theta_\lambda\preceq\Theta_*$, and let  $\delta_\lambda=\|\Theta_\lambda^{1/2}\Theta_*^{-1/2}-I\|$, so that $\delta_\lambda\in  [0,2]$. We have $0\preceq G\preceq \gamma^+\Theta_*^{-1}\preceq\gamma^+\Theta_\lambda^{-1}$
that is, $H,G$ satisfy (\ref{weq11a}) w.r.t. $\bar{\Theta}=\Theta_\lambda$. As a result, for our $h,G,H,u$ and the just defined $\bar{\Theta}$, relation (\ref{weq111a}) holds true:
\begin{equation}\label{from}
\begin{array}{l}
\multicolumn{1}{l}{\zeta\sim \SG(\A(u),\Theta_\lambda)\Rightarrow}\\
\ln\left(\bE_\zeta\left\{{\rm e}^{h^T\zeta+\half \zeta^TH\zeta}\right\}\right)\leq
-\half\ln\Det(I-\Theta_\lambda^{1/2}G\Theta_\lambda^{1/2})\\
\multicolumn{1}{r}{+\half[u;1]^TB^T\left[\Hplus+\left[H,h\right]^T
[\Theta_\lambda^{-1}-G]^{-1}\left[H,h\right]\right]B[u;1]}\\
\leq -\half\ln\Det(I-\Theta_\lambda^{1/2}G\Theta_\lambda^{1/2})\\
\multicolumn{1}{r}{+\half
\phi_\Z\left(B^T\left[\Hplus+\left[H,h\right]^T
[\Theta_\lambda^{-1}-G]^{-1}\left[H,h\right]\right]B\right)}\\
\end{array}
\end{equation}
(recall that $[u;1][u;1]^T\in\Z$). As a result,
\begin{equation}\label{eqeq2345}
\begin{array}{l}
\multicolumn{1}{l}{\zeta\sim \SG(\A(u),\Theta)\Rightarrow
\ln\left(\bE_\zeta\left\{{\rm e}^{h^T\zeta+\half \zeta^TH\zeta}\right\}\right)\leq
-\half\ln\Det(I-\Theta_\lambda^{1/2}G\Theta_\lambda^{1/2})}\\
\multicolumn{1}{r}{+\half\phi_\Z\left(B^T\left[\Hplus+\left[H,h\right]^T
[\Theta_*^{-1}-G]^{-1}\left[H,h\right]\right]B\right).}\\
\end{array}
\end{equation}
When deriving (\ref{eqeq2345}) from (\ref{from}), we have used that\\
--- $\Theta\preceq\Theta_\lambda$, so that when $\zeta\sim\SG(\A(u),\Theta)$, we have also $\zeta\sim\SG(\A(u),\Theta_\lambda)$,\\
--- $0\preceq \Theta_\lambda\preceq\Theta_*$   and $G\prec\Theta_*^{-1}$, whence
$[\Theta_\lambda^{-1}-G]^{-1}\preceq [\Theta_*^{-1}-G]^{-1}$,\\
--- $\Z\subset \bS^{n+1}_+$, whence $\phi_\Z$ is $\succeq$-monotone:  $\phi_Z(M)\leq\phi_\Z(N)$ whenever $M\preceq N$.
\par
By Lemma \ref{lemlogdet} applied with $\Theta_\lambda$ in the role of $\Theta$ and $\delta_\lambda$ in the role of $\delta$, we have
$$
\begin{array}{l}
-\half\ln\Det(I-\Theta_\lambda^{1/2}G\Theta_\lambda^{1/2})\\
\quad\leq -\half\ln\Det(I-\Theta_*^{1/2}G\Theta_*^{1/2})\\
\qquad+\half\Tr([\Theta_\lambda-\Theta_*]G)+{\delta_\lambda(2+\delta_\lambda)\over 2(1-\|\Theta_*^{1/2}G\Theta_*^{1/2}\|)}\|\Theta_*^{1/2}G\Theta_*^{1/2}\|_F^2.\\
\end{array}
$$
Consequently, (\ref{eqeq2345}) implies that
$$
\begin{array}{l}
\zeta\sim \SG(\A(u),\Theta)\Rightarrow\\
\ln\left(\bE_\zeta\left\{{\rm e}^{h^T\zeta+\half \zeta^TH\zeta}\right\}\right)\\
\qquad\qquad\quad\leq
-\half\ln\Det(I-\Theta_*^{1/2}G\Theta_*^{1/2})+\half\Tr([\Theta_\lambda-\Theta_*]G)\\
\qquad\qquad\qquad+{\delta_\lambda(2+\delta_\lambda)\over 2(1-\|\Theta_*^{1/2}G\Theta_*^{1/2}\|)}\|\Theta_*^{1/2}G\Theta_*^{1/2}\|_F^2\\
\qquad\qquad\qquad+\half\phi_\Z\left(B^T\left[\Hplus+\left[H,h\right]^T
[\Theta_*^{-1}-G]^{-1}\left[H,h\right]\right]B\right).\\
\end{array}
$$
The resulting inequality holds true for all small positive $\lambda$; taking $\lim\inf$ of the right hand side as $\lambda\to+0$, and recalling that $\Theta_0=\Theta$, we get
$$
\begin{array}{l}
\zeta\sim \SG(\A(u),\Theta)\Rightarrow\\
\ln\left(\bE_\zeta\left\{{\rm e}^{h^T\zeta+\half \zeta^TH\zeta}\right\}\right)\\
\qquad\qquad\quad\leq
-\half\ln\Det(I-\Theta_*^{1/2}G\Theta_*^{1/2})+\half\Tr([ {\Theta}-\Theta_*]G)\\
\qquad\qquad\qquad+{\delta(2+\delta)\over 2(1-\|\Theta_*^{1/2}G\Theta_*^{1/2}\|)}\|\Theta_*^{1/2}G\Theta_*^{1/2}\|_F^2\\
\qquad\qquad\qquad+\half\phi_\Z\left(B^T\left[\Hplus+\left[H,h\right]^T
[\Theta_*^{-1}-G]^{-1}\left[H,h\right]\right]B\right)\\
\end{array}
$$
(note that under the premise of Proposition \ref{concatenationSG}.i  we clearly have $\lim\inf_{\lambda\to+0}\delta_\lambda\leq\delta$). The right hand side of the resulting inequality is nothing
but $\Psi^\delta_{A,\Z}(h,H,G;\Theta)$, see  {\eqref{phiSG}} and we arrive at the inequality required in the conclusion of (\ref{allweneed}).

\paragraph{3$^0$.} To complete the proof of Proposition \ref{concatenationSG}.i, it remains to prove that the functions
$\Phi_{\A,\Z}$, $\Phi^\delta_{A, {\Z}}$ possess the
properties of continuity, convexity-concavity, and coerciveness announced in Proposition \ref{concatenationSG}.
Let us verify that this indeed is so for $\Phi^\delta_{A,\Z}$; reasoning to follow, with evident simplifications, is applicable to $\Phi_{A,\Z}$ as well.
\par
Observe, first, that by exactly the same reasons as in item {\bf 4$^0$} of the proof of Proposition \ref{concatenation}, the function $\Psi^\delta_{A,\Z}(h,H,G;\Theta)$ is real valued, continuous and convex-concave on the domain
$$
\widehat{\H}\times\Z=\{(h,H,G): -\gamma^+\Theta_*^{-1}\preceq H\preceq \gamma^+\Theta_*^{-1},0\preceq G\preceq\gamma^+\Theta_*^{-1},H\preceq G\}\times\Z.
$$
The function $\Phi^\delta_{A,\Z}(h,H;\Theta):\H\times\U\to\bR$ is obtained from $\Psi^\delta(h,H,G;\Theta)$ by the following two operations: we first minimize $\Psi^\delta_{A,\Z}(h,H,G;\Theta)$ over $G$ linked to $(h,H)$ by the convex constraints $0\preceq G\preceq \gamma^+\Theta_*^{-1}$ and $G\succeq H$, thus obtaining a function $$\bar{\Phi}(h,H;\Theta):\underbrace{\{(h,H):-\gamma^+\Theta_*^{-1}\preceq H\preceq\gamma^+\Theta_*^{-1}\}}_{\bar{\H}}\times\U\to\bR\cup\{+\infty\}\cup\{-\infty\}.$$
Second, we restrict the function $\bar{\Phi}(h,H;\Theta)$ from $\bar{\H}\times \U$ onto $\H\times\U$. For $(h,H)\in\bar{\H}$, the set of $G$'s linked to $(h,H)$ by the above convex constraints clearly is a nonempty compact set;
as a result, $\bar{\Phi}$ is a real-valued convex-concave function on $\bar{\H}\times \U$. From continuity of $\Psi^\delta_{A,\Z}$ on its domain it immediately follows that  $\Psi^\delta_{A,\Z}$ is bounded and uniformly continuous on every bounded subset of this domain, implying by evident reasons that
$\bar{\Phi}(h,H;\Theta)$ is bounded in every domain of the form $\bar{B}\times\U$, where $\bar{B}$ is a bounded subset of $\bar{\H}$, and is continuous on $\bar{B}\times \U$ in $\Theta\in \U$ with properly selected modulus of continuity  independent of $(h,H)\in \bar{B}$.
Besides this, by construction, $\H\subset\inter\bar{\H}$, implying that if $B$ is a convex compact subset of $\H$, it belongs to the interior of a properly selected convex compact subset $\bar{B}$ of $\bar{H}$. Since $\bar{\Phi}$ is bounded on $\bar{B}\times \U$ and is convex in $(h,H)$, the function $\bar{\Phi}$ is Lipschitz continuous in $(h,H)\in B$ with Lipschitz constant which can be selected to be independent of $\Theta\in \U$. Taking into account that $\H$ is convex and closed, the bottom line is that $\Phi^\delta_{A,\Z}$ is not just
a real-valued convex-concave function on the domain $\H\times\U$, it is also continuous on this domain.
\par
Coerciveness of $\Phi^\delta_{A,\Z}(h,H;\Theta)$ in $(h,H)$ is proved in exactly the same fashion as the similar property of function (\ref{phi}), see item  {\textbf{5$^0$}} in the proof of Proposition \ref{concatenation}. The proof of  item (i) of Proposition \ref{concatenationSG} is complete.
\par
Item (ii) of Proposition \ref{concatenationSG} can be derived from item  (i) of Proposition  {\ref{concatenationSG}} in exactly the same fashion as Proposition \ref{concatenation}.ii was derived from
 {Proposition} \ref{concatenation}.i.\hfill \qed

\subsection{Proof of Proposition \ref{propquad}}

Let us fix an input $x\in X$, and let $P$ be the corresponding distribution of $y^t$. Assuming that $x$ is a nuisance, at a given step $t$
a signal conclusion can take place only if there were $k$'s such that $\rho_{tk}<\infty$, and for some of these $k$'s it happened that
 $\phi_{tk\rho_{tk}}(y^t)< \alpha$. Invoking (\ref{case1}) with $\rho=\rho_{tk}$ and taking into account that ${\rm e}^{\cal SV}_{tk}(\rho_{tk})=\varkappa$ whenever $\rho_{tk}<\infty$,
 we see that the $P$-probability of the event $\phi_{tk\rho_{tk}}(y^t)< \alpha$ is at most ${\rm e}^\alpha\varkappa={\epsilon\over dK}$. Since there could be at most $K$ values of $k$ such that $\rho_{tk}<\infty$, in the situation under consideration the $P$-probability to terminate with the {signal} conclusion at step $t$ is at most $\epsilon/d$, and thus the $P$-probability of false alarm does not exceed $\epsilon$, as claimed.
  \par
  Now assume that $x\in X$ is a {signal} of shape $k$ and magnitude $\geq\rho$, and that $\rho\geq\rho_{tk}$ for some $k$. The latter may happen only when $\rho_{tk}<\infty$,
  implying that the detector $\phi_{tk\rho_{tk}}(\cdot)$ was used at time $t$. Assuming that the procedure did {\sl not} terminate at step $t$ with the signal conclusion, we should have $\phi_{tk\rho_{tk}}(y^t)\geq\alpha$, and invoking (\ref{case2}), we see that the $P$-probability of the outcome under consideration is at most ${\rm e}^{-\alpha}\exp\{{\cal SV}_{tk}(\rho_{tk})\}={\rm e}^{-\alpha}\varkappa=\sqrt{dK}\varkappa=\epsilon$.\hfill \qed

\end{document}